\documentclass[10pt,a4paper]{article}
\usepackage{amstext}
\usepackage{array}
\usepackage{amsfonts}
\usepackage{amssymb}
\usepackage{graphicx}
\usepackage{epstopdf}
\usepackage{algorithm}     
\usepackage{algpseudocode} 
\usepackage{url}
\usepackage{caption}
\usepackage{subfig}
\usepackage{csquotes}
\usepackage{tikz}
\usetikzlibrary{arrows}
\usetikzlibrary{patterns}
\usepackage{mathtools}
\usepackage{bbm}
\usepackage{amsthm}
\usepackage[textsize=tiny,textwidth=1.05in]{todonotes}

\newcommand{\vertiii}[1]{{\left\vert\kern-0.25ex\left\vert\kern-0.25ex\left\vert #1 
   \right\vert\kern-0.25ex\right\vert\kern-0.25ex\right\vert}}

\usepackage{pgfplots}
\usepackage{float}
\usepgfplotslibrary{external}

\newtheorem{theorem2}{Theorem}[section]

\newtheorem{remark}[theorem2]{Remark}

\DeclareMathOperator{\Tr}{Tr}
\usepackage{geometry}
\graphicspath{ {./fig/} }

\title{Mapping of coherent structures in parameterized flows by learning optimal transportation with Gaussian models}
\author{Angelo Iollo and Tommaso Taddei}
\date{}
\usepackage{fancyhdr}  
  
\fancypagestyle{plain}{
\fancyhead{}
\fancyhead[C]{\hfill Accepted for publication  on Journal of Computational Physics, September 2022}
}

\begin{document}

\maketitle

\begin{abstract}
 We present a general (i.e., independent of the underlying model) interpolation technique based on optimal transportation of Gaussian models for parametric advection-dominated problems. 
 The approach relies on a  scalar testing function  to identify 
the coherent structure we wish to track; a maximum likelihood estimator to identify a Gaussian  model of the coherent structure; and 
a nonlinear interpolation strategy that relies on   optimal transportation maps between Gaussian distributions.
 We  show that well-known self-similar solutions can be recast in the frame of optimal transportation by appropriate rescaling;
 we further   present several numerical examples to motivate our proposal  and  to assess strengths and limitations; finally, we discuss an extension to deal with more complex problems.
 \end{abstract}

\section{Introduction}
\label{sec:intro}

In science and engineering, it is important to identify low-rank approximations valid over a range of configurations (corresponding to different physical properties, different geometries or operational configurations). Low-rank approximations are of paramount importance in parameterized model order reduction (pMOR, \cite{brunton2020machine,lumley1967structure,hesthaven2016certified,quarteroni2015reduced,willcox2002balanced}) to speed up model evaluations in the limit of many queries, but also in optimization and uncertainty quantification to efficiently  generate samples from the solution manifold. In this paper, given two  snapshots of the solution manifold $U_0$, $U_1: \Omega\subset \mathbb{R}^n \to \mathbb{R}^d$, we wish to determine an interpolation $\widehat{U}: [0,1]\times \Omega \to \mathbb{R}^d$ such that 
$\widehat{U}(0,\cdot) = U_0$ and    $\widehat{U}(1,\cdot) = U_1$:
this task is one of the key enablers towards the implementation of    approximation strategies for parameterized systems {and is also relevant in itself}.
Our emphasis is on the development of a 
general (i.e., independent of the underlying parametric model),
interpretable  methodology that allows simple (i.e., non-intrusive) integration with high-fidelity codes and that is robust also  for small datasets.

The vast majority of data compression methods aims to determine linear low-rank approximations. If we denote by $U(x, \mu)$ the solution field, where $x=(x_1,\ldots,x_n)$ denotes the spatial variable and $\mu=(\mu_1,\ldots,\mu_p)$ denotes the vector of parameters, linear  approaches
consider approximations of the form
$$
\widehat{U}(x,\mu) = \sum_{i=1}^r \widehat{\alpha}_i(\mu) \zeta_i(x).
$$
Here, $\widehat{\alpha}_1,\ldots,\widehat{\alpha}_r$ are parameter-dependent coefficients  that can be obtained by solving  a reduced-order model (ROM), while $\zeta_1,\ldots,\zeta_r$ are  a reduced-order basis (ROB) that is computed by exploring the parameter domain. 
Linear models can be interpreted as a generalization  of convex interpolations of two snapshots $U_0, U_1$, that is
\begin{equation}
\label{eq:convex_interpolation}
\widehat{U}^{\rm co}(s, x)
\, =  \,
(1 - s) U_0 (x)
 \, + \,
s U_1(x) 
\quad
s\in [0,1], x\in \mathbb{R}^n.
\end{equation}
The use of  linear methods relies on the assumption that the problem of interest exhibits linear coherent structures.

There exists  a broad class of problems for which linear methods are  effective. To provide concrete examples, the presence  of coherent structures in turbulent flows
provides physical foundations for the use of linear methods in numerous applications  in flow control and design \cite{berkooz1993proper}; 
evanescence of high-frequency modes for diffusion-dominated problems \cite{munjal1987acoustics} is at the foundation of 
component-mode synthesis \cite{craig1968coupling} and more recently of component-based MOR strategies \cite{huynh2013static}. Despite the successes of linear methods, there exists a broad class of problems of interest in engineering for which linear methods are highly inaccurate: this motivates the development of nonlinear methods.

We present a general  approach that relies on optimal transportation
\cite{villani2003topics}
to perform accurate nonlinear interpolations  between solution snapshots.
First, we rely on  a scalar testing function (cf. section \ref{sec:gauss_model_coherent}) to derive a Gaussian model $g[U]$  of the solution field; then, we rely on well-known results for optimal transportation of Gaussian distributions to determine the optimal transport  mapping  $X_g$ from  $g[U_0]$ and  $g[U_1]$ and the 
optimal transport  mapping  $Y_g$ from  $g[U_1]$ and  $g[U_0]$; finally, we define the nonlinear interpolation
\begin{equation}
\label{eq:convex_displacement_interpolation}
\widehat{U}(s, x)
\, =  \,
(1 - s) U_0 \circ 
  W_g(s, x)  
 \, + \,
s U_1 \circ T_g (1-s,x),
\quad
s\in [0,1], x\in \mathbb{R}^n,
\end{equation}
where $W_g(s, x)= (1 - s) x + s Y_g(x)$ and 
$T_g(s, x)= (1 - s) x + s X_g(x)$.
We refer to $\widehat{U}$ as \emph{convex displacement interpolation} (CDI) due to the analogy with displacement interpolation (cf. section \ref{sec:optimal_transportation}) {and the more elementary convex interpolation \eqref{eq:convex_interpolation}}.
We present several numerical examples to motivate our proposal; furthermore, we   show in section \ref{sec:motivating_examples} that well-known self-similar solutions can be recast in the frame of optimal transportation by appropriate rescaling.

The use of optimal transportation theory to devise nonlinear interpolation  has been considered in several works in the pMOR literature,
 \cite{bernard2018reduced,ehrlacher2020nonlinear,iollo2014advection}.
Here, we apply optimal transportation to a model of the solution field: as a result, the parametric  field of interest $U$ does not have to be neither scalar nor positive and does not have to  fulfill   conservation of mass over the parameter domain.

The proposed approach  shares relevant features with Lagrangian or registration-based approaches to pMOR (\cite{mojgani2017arbitrary,ohlberger2013nonlinear,taddei2020registration}) and also with the works on Lagrangian coherent structures (LCS, \cite{haller2015lagrangian,peacock2013lagrangian}) in the field of nonlinear dynamics. In particular, the feature-based Gaussian model of the field is similar in scope to the registration sensor introduced in \cite{taddei2021space} and also to shock-capturing sensors used for high-order schemes (\cite{nicoud1999subgrid,persson2006sub}). On the other hand, we remark that, while our emphasis is on the development of predictive models for parametric systems, LCS literature mainly focuses on the physical understanding --- and subsequently the control --- of chaotic systems.
{Furthermore, while registration methods rely on the introduction of a single reference configuration, CDI is inherently Eulerian.
}

As extensively discussed in the example of section \ref{sec:transonic_flow}, our approach might suffer from (i) the presence of boundaries, and (ii) the presence of multiple coherent structures that we wish to track. In section \ref{sec:extension},  we discuss how to extend the approach to deal with more complex problems: the key ideas are to replace the Gaussian model with a mixture of Gaussian models and to replace the affine-in-$s$ maps $W_g, T_g$ with suitable nonlinear maps.

The outline of the paper is as follows.
In section \ref{sec:optimal_transportation}, we present a short introduction to optimal transportation: we provide a number of references and we introduce relevant notation. {Then, we derive the CDI for cumulative distribution functions and we offer insights about the proposed form}.
In section \ref{sec:motivating_examples}, we illustrate the connection between self-similarity and optimal transportation
{--- and ultimately CDI ---}
through the vehicle of several examples. 
In section \ref{sec:method}, we present our method  in its elementary form and we 
provide numerical investigations.
In section \ref{sec:extension}, we discuss the extension of the elementary approach to deal with more complex problems.
Section \ref{sec:conclusions} concludes the paper by offering a short summary and several perspectives.

\section{Convex displacement interpolation}
\label{sec:optimal_transportation}

\subsection{Optimal transportation}

We  resort to optimal transportation theory \cite{villani2003topics} to model the displacement of coherent structures with respect to parameter variations, as it was already done with different approaches in  \cite{bernard2018reduced,iollo2014advection}.
We introduce the probability measures $\mathbb{P}_0,\mathbb{P}_1$ with probability density functions (pdfs) $\rho_0,\rho_1$ and cumulative distribution functions (cdfs)
$F_0,F_1$,
$$
F_i(x) = \int_{-\infty}^{x_1} \ldots \int_{-\infty}^{x_n} \rho_i(x') \, dx',
\quad i=0,1, \;\; 
x \in \mathbb{R}^n.
$$
We assume that $\rho_0,\rho_1$ have finite second-order moments. We say that $X:\mathbb{R}^n \to \mathbb{R}^n$ transports $\mathbb{P}_0$ to $\mathbb{P}_1$ if 
$\mathbb{P}_1(B)
= \mathbb{P}_0(X^{-1}(B) )$ for all
$\mathbb{P}_1$-measurable  sets  
$B$, with $X^{-1}(B) := \{\xi \in \mathbb{R}^n : X(\xi) \in B    \}$, and we use notation $\mathbb{P}_1 = X_{\#}\mathbb{P}_0$. Note that the latter implies  local mass conservation
\begin{equation}
\label{eq:local_mass_conservation}
\rho_0(\xi) = \rho_1(X(\xi))  \big|  {\rm det} \nabla_{\xi} X(\xi) \big|, \quad \forall \, \xi \in \mathbb{R}^n,
\end{equation}
or equivalently
\begin{equation}
\label{eq:local_mass_conservation_cdf}
F_0(\xi) = F_1(X(\xi)), \quad \forall \, \xi \in \mathbb{R}^n.
\end{equation}

With this notation, we can introduce the Monge's optimal transport problem as follows:
find $X:\mathbb{R}^n\to \mathbb{R}^n$ to minimize
\begin{equation}
\label{eq:monge_problem}
I(X;  \rho_0, \rho_1) \, = \, \int_{\mathbb{R}^n} \| X(\xi) - \xi \|_2^2 \, \rho_0(\xi) \, d\xi,
\quad
{\rm subject \; to} \; \eqref{eq:local_mass_conservation}.
\end{equation}
Under the assumptions stated above, 
there exists a unique convex potential $\Psi$: $\mathbb{R}^n \xrightarrow[]{}\mathbb{R}$,  such that the mapping $X =\nabla_{\xi} \Psi$ minimizes \eqref{eq:monge_problem}:
this result follows from a duality principle in convex optimization introduced in \cite{kantorovich1942translocation} that is linked to the polar factorization and monotone rearrangement of vector-valued functions \cite{brenier1987polar};  the full proof is based on the convexity of the potential minimizing $I(X; \rho_0,\rho_1)$ and the existence and uniqueness of monotone measure-preserving maps \cite{mccann1995existence}.
If we denote by $J$ the minimum of $I(X; \rho_0,\rho_1)$
over all transport maps $X$, existence and uniqueness of the potential $\Psi$  imply that
$$
W_2(\rho_0, \rho_1)=\displaystyle \sqrt{J}
$$ 
is a distance function  between probability measures; $W_2$ is 
known as the Wasserstein metric.

{If we denote by $\mathbb{W}_2(\mathbb{R}^n)$ the space of probabilities in $\mathbb{R}^n$ endowed with the distance $W_2$, it is possible to show that $\mathbb{W}_2(\mathbb{R}^n)$ is a geodesic space and 
\begin{equation}
\label{eq:mccann_interpolation_plus}
\mathbb{P}^+(s)
:=\left(
(1-s) Id + s X
\right)_{\#} \mathbb{P}_0,
\quad
s\in [0,1]
\end{equation}
is the unique geodesic curve that connects 
$\mathbb{P}_0$ and $\mathbb{P}_1$, and has constant speed 
\cite[Chapter 5]{santambrogio2015optimal}.
Eq. \eqref{eq:mccann_interpolation_plus} provides a nonlinear interpolation --- dubbed as displacement or McCann interpolation  --- that is very different from the linear interpolation $(1-s) \mathbb{P}_0 + s \mathbb{P}_1$ and is useful in numerous applications. Eq. \eqref{eq:mccann_interpolation_plus} also provides the foundations for the CDI introduced in the next section.
}

\subsection{Interpolation procedure}

{We denote by $Y:\mathbb{R}^n \to \mathbb{R}^n$ the optimal transport map that connects 
 $\mathbb{P}_1$ to  $\mathbb{P}_0$ and we introduce the reverse interpolation:
 \begin{equation}
\label{eq:mccann_interpolation_minus}
\mathbb{P}^-(s)
:=\left(
(1-s) Id + s Y
\right)_{\#} \mathbb{P}_1,
\quad
s\in [0,1]
\end{equation}
Provided that $X$ is invertible, it is possible to prove that 
 $Y=X^{-1}$; furthermore, since $\mathbb{P}^-,\mathbb{P}^+$ have constant speed and $\mathbb{W}_2(\mathbb{R}^n)$ is a geodesic space, we must have
\begin{equation}
\label{eq:key_fact}
\mathbb{P}^+(s) = \mathbb{P}^-(1-s).
\end{equation}
If we introduce the maps
\begin{subequations}
\label{eq:exact_interpolation}
\begin{equation}
\label{eq:exact_interpolation_a} 
T(s,\xi) = (1-s) \xi + s X(\xi), \quad
W(s,x) = (1-s) x+ s Y(x),
\end{equation}
the identity \eqref{eq:key_fact} implies that the cdf $U(s,\cdot)$  associated with $\mathbb{P}^+$ satisfies
\begin{equation}
\label{eq:exact_interpolation_b}
U(s,y) 
=F_0\left( T^{-1}(s,y)  \right)
=F_1\left( W^{-1}(1-s,y)  \right)
\quad
\forall \, s\in [0,1] \;\;  y\in \mathbb{R}^n.
\end{equation}
\end{subequations}
Then, if we linearize $T^{-1}(s,\cdot)$ and $W^{-1}(s,\cdot)$,
\begin{subequations}
\label{eq:CDI_interpolation_section2}
\begin{equation}
\label{eq:CDI_interpolation_section2_a}
\begin{array}{l}
\displaystyle{
T^{-1}(s,y) \approx 
(1- s) T^{-1}(0,\xi) 
+
s T^{-1}(1,\xi) 
=
(1- s) y + s Y(y)
=
W(s,y),
} \\[3mm]
\displaystyle{
W^{-1}(s,y) \approx 
(1- s) W^{-1}(0,\xi) 
+
s W^{-1}(1,\xi) 
=
(1- s) y + s X(y)
=
T(s,y),
} \\ 
\end{array}
\end{equation}
we obtain 
\begin{equation}
\label{eq:CDI_interpolation_section2_b}
U(s,y)  \approx 
(1-s)
 F_0\left( W(s,y)  \right) \, 
+ \, s \, F_1\left( T(1-s,y)  \right)
=: \widehat{U}(s,y)
\quad
\forall \, s\in [0,1] \;\;  y\in \mathbb{R}^n, 
\end{equation}
which is the CDI  \eqref{eq:convex_displacement_interpolation} for the cdfs $F_0,F_1$.
\end{subequations}

The linearization \eqref{eq:CDI_interpolation_section2_a} introduces an approximation in the identities \eqref{eq:exact_interpolation_b}: the combination of direct (cf. \eqref{eq:mccann_interpolation_plus}) and
reverse 
(cf. \eqref{eq:mccann_interpolation_minus}) displacement interpolations ensures exact interpolation at end points
---
that is, 
$\widehat{U}(s,\cdot)= {U}(s,\cdot)$ for $s=0,1$
---
and symmetry
---
that is, if we build $\widehat{U}'$ by interchanging $\rho_0$ with $\rho_1$, we obtain 
$\widehat{U}'(s,\cdot)= \widehat{U}(1-s,\cdot)$.
We observe that the CDI stems from the optimal transport theory but it can also be  applied to arbitrary parametric fields through the vehicle of suitable probabilistic models of the coherent structures of the solution: in section \ref{sec:method}, we show how to derive Gaussian models for arbitrary fields using user-defined scalar testing functions.
In this section, we have considered interpolations over the whole space, i.e.,  $\Omega=\mathbb{R}^n$. If $\Omega$ is not convex, the maps \eqref{eq:exact_interpolation_a} might not map $\Omega$ in itself and so \eqref{eq:CDI_interpolation_section2_b} might not be well-defined. In section \ref{sec:extension}, we investigate a boundary-aware generalization of the CDI
that relies on nonlinear-in-$s$ maps $T,W$ 
 to deal with this case.
}

\begin{remark}
\textbf{Lagrangian interpolation.}
Convex displacement interpolation is related to the following nonlinear interpolation:
\begin{equation}
\label{eq:Lagrangian_displacement}
\widehat{U}(s, x)
\, =  \, \left(
(1 - s) U_0  \, + \,
s U_1 \circ X  
\right)
\circ W (s,x),
\quad
s\in [0,1], x\in \mathbb{R}^n.
\end{equation}
The latter performs linear interpolation of the mapped field $\{ \widetilde{U}(s,\cdot) = U(s,  T(s,\cdot) ) \, : \,s\in [0,1] \}$ and can thus be referred to  as Lagrangian interpolation: it is indeed consistent with Lagrangian (or registration-based) approaches presented in the MOR literature, (e.g. \cite{iollo2014advection,taddei2020registration,taddei2021space}).
.
Similarly to \eqref{eq:convex_displacement_interpolation}, the nonlinear interpolation \eqref{eq:Lagrangian_displacement}  satisfies  $\widehat{U}(0, \cdot) = U_0$ and 
$\widehat{U}(1, \cdot) = U_1$;
however, $\widehat{U}$  is not symmetric:
interchanging $U_0$ with  $U_1$ leads to a different nonlinear  interpolant.
{More fundamentally, \eqref{eq:Lagrangian_displacement} relies on the introduction of a reference configuration --- in this case $U_0$ --- while CDI is inherently Eulerian}
 
 Lagrangian interpolation can   be generalized to the case of multiple snapshots: the strategy  is equivalent to the one discussed in 
\cite{taddei2020registration,taddei2021space} and is here outlined for completeness. Given snapshots $\{  (s_k, U^k:=U(s_k,\cdot)    \}_{k=0}^K$,  
 (i) we compute the forward maps $X^k$ between $U_0$ and $U^k$ for $k=0,\ldots,K$;
 (ii) we infer the parametric map $s \in [0,1] \mapsto \widehat{T}(s,\cdot)$;
 (iii) we define the mapped snapshots
 $\widetilde{U}^k : = U^k \circ \widehat{T}(s_k,\cdot)$ for $k=0,\ldots,K$;
 (iv) we define the Lagrangian interpolant
\begin{equation}
\label{eq:Lagrangian_displacement_plus}
\widehat{U}(s, x)
\, =  \, \left(
\sum_{k=0}^K \alpha_k(s) \,\widetilde{U}^k
\right)
\circ \widehat{T}^{-1}(s,x),
\quad
s\in [0,1], x\in \mathbb{R}^n,
\end{equation}
 where the coefficients 
$\alpha_0,\ldots,\alpha_K:[0,1]\to \mathbb{R}$ should also be learned based on the available training data, or based on a mathematical physical model.
Note that inference of the parametric map
$\widehat{T}$ at step (ii) should preserve bijectivity: we refer to the above-mentioned literature for a discussion on this issue.
\end{remark}

\section{Motivating examples}
\label{sec:motivating_examples}

Self-similarity plays a fundamental role in physics and we show in the next examples that well-known self-similar solutions to  PDEs can be recast in the frame of optimal transportation by appropriate rescaling. 
{More precisely, in all examples below
we identify a function 
$F:[t_0,t_1]\times \mathbb{R}^n \to \mathbb{R}$ of the solution such that
\begin{equation}
\label{eq:mccann_cdf}
F(t, x) = F(t_0, T^{-1}(s(t), x)), \quad
{\rm with} \; T(s, \xi) = (1-s)\xi + X(\xi),
\end{equation}
where $X$ solves a suitable optimal transportation problem and
$s:[t_0,t_1]\to [0,1]$ is a suitable rescaling function. Note that \eqref{eq:mccann_cdf} reads as the McCann interpolation 
\eqref{eq:mccann_interpolation_plus}
for cumulative distribution functions.
}
In the remainder, we repeatedly use the Brenier's theorem (\cite{brenier1987polar}), which states that
given two densities $\rho_0$ and $\rho_1$
 there exists a unique  function $X$ that is the  gradient of a convex function and transports 
 $\rho_0$  onto $\rho_1$.

\subsection{The heat kernel}
\label{sec:heat_equation}

\subsubsection{Convex potential}

Let us consider the heat kernel $K:  \mathbb{R}_+ \times \mathbb{R}^n \rightarrow \mathbb{R}_+$, $(t,x)\mapsto K(t,x)$,
\begin{equation}
\label{eq:heat_kernel}
K(t,x)={\frac{1}{(4\pi t)^{n/2}}} \, e^
{-\frac{\| x  \|_2^2}{4 t}} \,
\end{equation}
that satisfies the initial value problem:
\begin{equation}
\label{eq:heat_equation}
\left\{
\begin{array}{ll}
{\displaystyle\frac {\partial K}{\partial t}}=\Delta K  &
{\rm in} \; \mathbb{R}^+ \times \mathbb{R}^n  ,  \\[3mm]
K(0,x) = \delta_0(x) & 
{\rm in} \;  \mathbb{R}^n ,
\\
\end{array}
\right.
\end{equation}
where $\delta_0$ is the 
Dirac mass concentrated in $x=0$. Of course, we have for all $ t_0, t_1 \in \mathbb{R}^+$:
$$
\int_{\mathbb{R}^n}  K(t_0,\xi) d\xi=\int_{\mathbb{R}^n}  K(t_1,x) dx = 1.
$$
Define now the convex potential $\Psi: \mathbb{R}_+ \times \mathbb{R}_+ \times  \mathbb{R}^n \rightarrow \mathbb{R}_+$, 
$$
\Psi(t_0,t_1,\xi)=\frac{1}{2} \sqrt{\frac{t_1}{t_0}} \, 
\| \xi \|_2^2.
$$ 
Note that $\nabla_{\xi} \Psi =\sqrt{  \displaystyle\frac{t_1}{t_0}} \xi$ and 
${\rm det} (  \nabla_{\xi}^2 \Psi     )
= \left(\displaystyle \frac{t_1}{t_0}  \right)^{n/2}$.  It is thus easy to verify that 
$$
K(t_0,\cdot )=K(t_1,\nabla_\xi \Psi)\det{\left(\nabla^2_\xi \Psi\right)}
\quad
\forall \; t_0,t_1\in \mathbb{R}_+,
$$
which corresponds to  \eqref{eq:local_mass_conservation} for $\rho_0 = K(t_0,\cdot)$ and 
$\rho_1 = K(t_1,\cdot)$. Since 
$\Psi$ is a convex function whose gradient  satisfies \eqref{eq:local_mass_conservation},   $X = \nabla_{\xi} \Psi$ must be the unique solution to \eqref{eq:monge_problem}.
 
\subsubsection{Displacement interpolation}

{We introduce  the forward mapping and the corresponding backward map:
$$
 T(s,\xi)= \, (1 - s)   \, \xi+ \, s \,  \sqrt{\frac{t_1}{t_0}} \, \xi
 = \left(  
\frac{\sqrt{t_0}  +  (\sqrt{t_1} - \sqrt{t_0}    ) s}{\sqrt{t_0} }  \right) \xi,
\quad
T^{-1}(s, x)=  \,  \left(  
\frac{\sqrt{t_0}  +  (\sqrt{t_1} - \sqrt{t_0}    ) s}{\sqrt{t_0} }  \right)^{-1} x.
$$
If we define the parameter re-scaling
$s:[t_0,t_1]\to [0,1]$ such that
$$
s(t)=\frac{\sqrt{t}-\sqrt{t_0}}{\sqrt{t_1}-\sqrt{t_0}},
$$
we find $T^{-1}(s(t), x) = \sqrt{ \frac{t}{t_0} } x$ and then
$$
K(t, x) = \frac{1}{(4\pi t)^{n/2}} e^{ - \frac{\|x\|_2^2}{4 t}   }
=
\left( \frac{t}{t_0} \right)^{n/2}
K\left(t_0,  T^{-1}(s(t), x) \right).
$$
The latter implies that  the rescaled solution 
$F:[t_0,t_1]\times \mathbb{R}^n \to \mathbb{R}$ such that
$F(t, x) = t^{n/2} K(t, x)$ satisfies \eqref{eq:mccann_cdf}: we conclude that the displacement interpolation --- with suitable rescaling --- is an exact solution to the heat equation for all $t\in [t_0,t_1]$ up to  a multiplicative factor.
}


 
\subsection{Nonlinear diffusion}
\label{sec:no-lin}

\subsubsection{Convex potential}

A suitable model for diffusion of heat in hot plasma, very intense thermal waves or diffusion in porous media \cite{barenblatt1996scaling} is the following nonlinear diffusion equation:
\begin{equation}
    \label{non-linheat}
\frac{\partial \theta}{\partial t}=\Delta  \theta^m
\end{equation}
where $\theta: \mathbb{R}_+\times\mathbb{R}^n \rightarrow \mathbb{R}_+$, $(t,x)\mapsto \theta(t,x)$, and $m \in \mathbb{N}$ with $m>1$. For example, for an instantaneous release of heat at time $t=0$ and concentrated at the origin, this equation admits the so-called ZKB \cite{barenblatt1996scaling} solution profile:
\begin{equation}
\label{zkb}
B(t,x)=t^{-\alpha} \left(\left(  C-k \| x \|_2^2 t^{-2\beta}\right)^+\right)^{1/(m-1)}
\end{equation}
where $\alpha=\displaystyle\frac{n}{n(m-1)+2}$, $\beta=\displaystyle \frac{\alpha}{n}$, $k=\displaystyle \frac{(m-1)\alpha}{2  m n}$, $C$ is a positive constant and $z^+=\max(z,0)$. The $L^1$ norm of the ZKB profile, $\displaystyle \int_{\mathbb{R}^n} B(t,x) \, dx$, is time invariant and equal to the ``initial heat'' released.

We introduce the convex potential    $\Psi: \mathbb{R}_+ \times \mathbb{R}_+ \times \mathbb{R}^n \to \mathbb{R}_+$,  
$$
\Psi(t_0,t_1,\xi)= \frac{1}{2}\left( \frac{t_1}{t_0}\right)^{\beta}\left\|\xi\right\|_2^2,
$$
Note that  
$\nabla_\xi \Psi=
\left( t_1/t_0 \right) ^{\beta}\,\xi$ and  
$\det{(\nabla^2_\xi \Psi)=(t_1/t_0)^\alpha}$: we thus obtain   
$$
B(t_0,\xi)=B\left(t_1,\xi \left(\frac{t_1}{t_0}\right)^{\beta} \right)\,\left(\frac{t_1}{t_0}\right)^{\alpha},
$$
which corresponds to \eqref{eq:local_mass_conservation} for 
$\rho_0(\xi)=B(t_0,\xi)$ and $\rho_1(x)=B(t_1,x)$.
As in the previous case, since $\Psi$ is convex and satisfies  \eqref{eq:local_mass_conservation}, we must have that 
$\nabla_\xi \Psi$ is the unique optimal forward mass transportation  between the solutions at $t_0$ and $t_1$.
 
\subsubsection{Displacement interpolation}
As in the heat kernel case, {we introduce the forward   mapping
and the corresponding backward mapping:
$$
    T(s,\xi)=(1-s)  \, \xi+s \left(\frac{t_1}{t_0}\right)^{\beta} \, \xi
=  \left( 1+\left(\left(\displaystyle\frac{t_1}{t_0}\right)^{\beta}-1\right)\,s  \right) \xi,
\quad
T^{-1}(s,x)  =
\frac{x}{1+\left(\left(\displaystyle\frac{t_1}{t_0}\right)^{\beta}-1\right)\,s}.
$$
Given the rescaling $s(t)=\frac{t^\beta-t_0^\beta}{t_1^\beta-t_0^\beta}$,  
we obtain 
$$
1+\left(\left(\displaystyle\frac{t_1}{t_0}\right)^{\beta}-1\right)\,s(t) = \left( \frac{t}{t_0} \right)^{\beta} = \left( \frac{t}{t_0} \right)^{\alpha/n},
$$ 
which implies $T^{-1}(s(t),x) =\left( \frac{t}{t_0} \right)^{\alpha/n} x $
and then
$$
B(t,x) 
=
 \left( \frac{t}{t_0} \right)^{-\alpha}
 B\left(t_0,
T^{-1}(s(t),x)
\right).
$$
The latter implies that  the rescaled solution 
$F:[t_0,t_1]\times \mathbb{R}^n \to \mathbb{R}$ such that
$F(t, x) = t^{\alpha} B(t, x)$ satisfies \eqref{eq:mccann_cdf}.}
Again, in model-order reduction vocabulary, we have that given two snapshots of the solution for different parameter values, their  displacement  interpolation  is an exact solution to the PDE up to a multiplicative constant, provided that the parameter is appropriately rescaled.

\subsection{Conservation laws}
\label{sec:example_euler}
The Euler equations for an  inviscid compressible ideal fluid flow in one dimension are given by: 
$$
\frac{\partial U}{\partial t}+\frac{\partial F}{\partial x}=0
$$
where $U=\left(\rho, \rho u,E\right)$ is the vector of conserved variables, $F=\left(\rho u,\rho u^2+p,(E+p) u\right)$ is the flux,  $\rho$ is the density, $p$ is the pressure, $E$ is the total internal energy per unit volume and $u$ is the velocity.
We denote by $\gamma>0$ the ratio of specific heats and we denote by $a = \displaystyle\sqrt{\frac{\gamma p}{\rho}}$ the speed of sound.

\subsubsection{Displacement interpolation of simple wave solutions}
\label{sec:simple_wave}

We assume that the flow is isentropic and does not contain shock waves for all $t\in (0, t^{\star}) = A$: under this assumption, the equation of state reduces to $p = C \rho^{\gamma}$; the Cauchy problem is well-posed in $A\times \mathbb{R}$ and all physical quantities can be traced back to the value at time $t=0$. It is possible to show that the Euler system admits the two Riemann invariants $R^{\pm} =
\displaystyle
\frac{2}{\gamma-1} a \pm u$ that satisfy the equations:
\begin{equation}
\label{eq:riemann_invariants}
\frac{\partial R^{\pm}}{\partial t}
+ (u \pm a) \, \frac{\partial R^{\pm}}{\partial x}
= 0 \quad {\rm in} \; A \times \mathbb{R}.
\end{equation}
In this section, we consider flows for which one of the two invariants, say $R^-$, is constant and equal to $c\in \mathbb{R}$ at time $t=0$; the corresponding solution to the Euler equations is known as simple wave (\cite{landau2013course}).

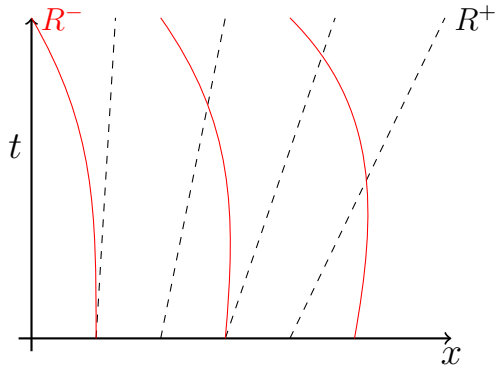
\begin{figure}[h]
\centering
\begin{tikzpicture}[scale=0.85]
\linethickness{0.3 mm}
\linethickness{0.3 mm}

\draw[thick, ->]  (-0.2,0)--(6.5,0);
\coordinate [label={below:  {\Large {$x$}}}] (E) at (6.5,0) ;
\draw[thick, ->]  (0,-0.2)--(0,5);
\coordinate [label={left:  {\Large {$t$}}}] (E) at (0,3) ;
			
\draw[dashed]  (1,0)--(1.3,5);
\draw[dashed]  (2,0)--(3,5);		
\draw[dashed]  (3,0)--(4.7,5);
\draw[dashed]  (4,0)--(6.4,5);		
					
%
\draw [,red] (1,0) to [out=90,in=-60] (0,5);
\draw [-,red] (3,0) to [out=85,in=-55] (2,5);
\draw [-,red] (5,0) to [out=80,in=-45] (4,5);
\coordinate [label={right:  {\large {$R^+$}}}] (E) at (6.4, 5) ;
\coordinate [label={right:  {\large {${\color{red} R^-}$}}}] (E) at (0, 5) ;

\end{tikzpicture}

\caption{Simple wave solution: the right-going characteristics are straight lines. $X^+(t,\xi)=\xi+(u_0(\xi)+a_0(\xi))\, t $ where $u_0$ and $a_0$ are the initial condition for the state variables.}
\label{sing_w}
\end{figure}

We define the characteristics $X^{\pm}: A \times \mathbb{R} \to \mathbb{R}$ such that
\begin{equation}
\label{eq:characteristics_simplewave}
\left\{
\begin{array}{ll}
\displaystyle{
\frac{d X^{\pm}}{d t} (t, \xi) = u(t,X^{\pm}(t,\xi) ) \,  \pm \, a(t,X^{\pm}(t,\xi) )
} & \displaystyle{ t\in A } \\[3mm]
X^{\pm} (0, \xi) = \xi &
\\
\end{array}
\right.
\end{equation}
Combining \eqref{eq:riemann_invariants} with \eqref{eq:characteristics_simplewave}, we find that 
$\displaystyle \frac{d}{dt}   R^{\pm}( t,  X^{\pm} (t,\xi)  )  = 0$, which implies that $ R^+$ (resp. $R^-$) is constant on the characteristic $X^+$ (resp. $X^-$). Since $R^-$ is constant at $t=0$, we must have that $R^-(t,x) = c$ for all $(t,x)\in A\times \mathbb{R}$. This implies that $u(t,x) = c+\displaystyle\frac{2}{\gamma-1} a(t,x)$ and then 
$$
\begin{array}{l}
\displaystyle{
a(t, X^+(t,\xi)) =  R^{+}( t,  X^{+} (t,\xi)  )  - u(t, X^+(t,\xi))  =
 R_0^{+}( \xi )  - c - \frac{2}{\gamma-1} a(t,X^+(t,\xi)) 
} \\
\displaystyle{
\hfill
 \Rightarrow
a(t, X^+(t,\xi))  = (\gamma-1)\frac{ R_0^{+}( \xi ) - c      }{\gamma+1} ,
}
\\
\end{array}
$$
and in particular $a$ (and thus all state variables) are constant on the $X^+$ characteristic.   Furthermore, since $u+a$ in  constant on $X^+$, we must have that
$\displaystyle\frac{d X^+}{dt}$ is constant and thus the $X^+$ characteristics  are straight lines
(cf. Figure \ref{sing_w}) and satisfy
$X^+(t,\xi)=\xi+ (u_0(\xi)+a_0(\xi))\, t$.  Due to  the assumption on the smoothness of the flow, characteristics do not coalesce for $t\in A$:
as a result, we find that $X^+(t,\cdot)$ is bijective in $\mathbb{R}$ for all $t\in A$.

We denote by $Y^+$ the inverse map of $X^+$ such that $X^+(t, Y^+(t,x)) = x$ for all
$(t,x) \in A \times \mathbb{R}$
($Y^+$ is the mapping that associates $(t,x)$ to the foot of the corresponding $X^+$ characteristic); we have that the pressure $p$ (and any other state variable) satisfies
\begin{equation}
\label{eq:simple_wave_cdf}
p(t,x) = p_0( Y^+(t,x)   ) = p_0(\xi).
\end{equation}
{We observe that \eqref{eq:simple_wave_cdf} is of the form 
\eqref{eq:mccann_cdf}:  to conclude, it thus suffices  to verify  that 
$X^{+}(t,\cdot)$  solves a suitable optimal transportation problem.

Towards this end, we define 
the scalar field 
$\Psi: A \times \mathbb{R} \to \mathbb{R}$ such that
$$
\Psi(t,\xi) \, = \, \int_0^{\xi} X^+(t,\xi') \, d \xi'.
$$
Since $X^+(t, \cdot)$ is bijective in $\mathbb{R}$ for all  $t\in A$, we must have $\displaystyle\frac{\partial X^+}{\partial \xi} > 0$: as a result, $\Psi$ is convex in the second argument; furthermore, 
by differentiating \eqref{eq:simple_wave_cdf}, we find
$$
\frac{\partial p}{\partial x} = \frac{\partial p_0}{\partial \xi}\frac{\partial Y^+}{\partial x} ,
$$
which  corresponds to \eqref{eq:local_mass_conservation} for 
$\rho_0 =  \frac{\partial p_0}{\partial \xi}$ and
$\rho_1 =  \frac{\partial p(t,\cdot)}{\partial x}$, provided that $\frac{\partial p_0}{\partial \xi}>0$\footnote{
The result can be trivially extended to the more general case by taking the absolute value of $\frac{\partial p_0}{\partial \xi}$.
}.
In conclusion, for any $t\in A$,  $ X^+(t,\cdot)$ is  the unique optimal forward mass transportation between  $ \displaystyle \frac{\partial p_0}{\partial \xi}  $ and  $ \displaystyle \frac{\partial p(t,\cdot)}{\partial x}$.
}
Application of optimal transportation thus  detects the appropriate self-similarity transformation associated with the problem.

\subsubsection{Displacement interpolation for Riemann problems }
\begin{figure}
    \centering
\subfloat[] { 
  \includegraphics[width=.45\linewidth]{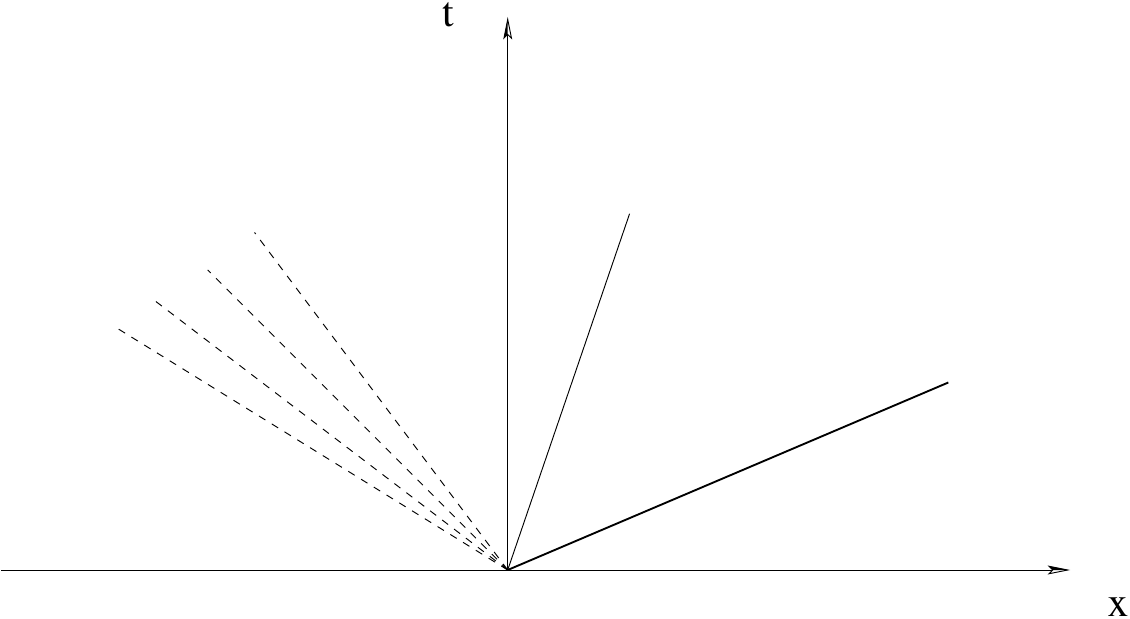} 
  }
  ~~
 \subfloat[] {  
  \includegraphics[width=.45\linewidth]{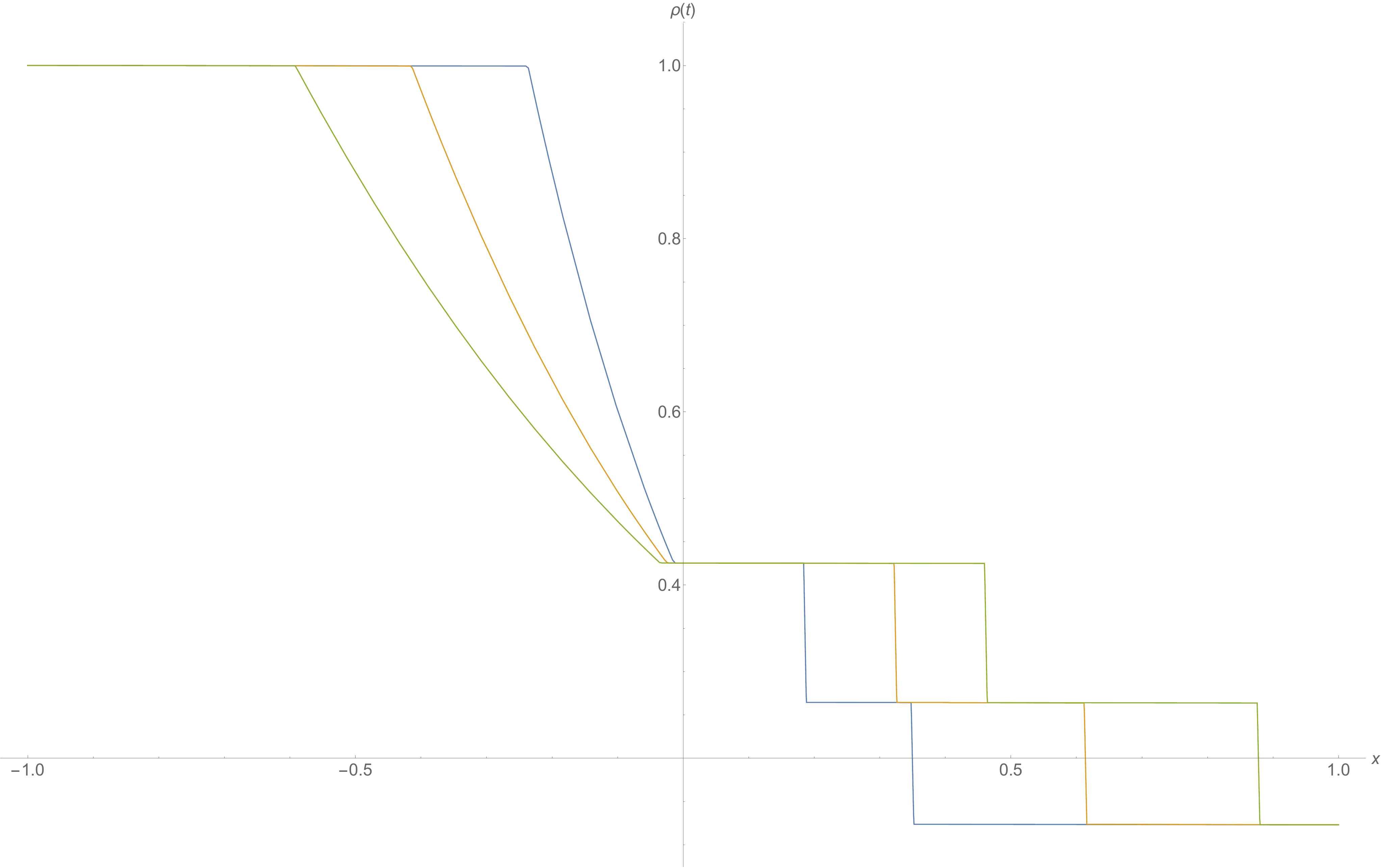} 
  }
  
\caption{Sod problem.
(a) shock tube solution pattern: a receding expansion fan and two forward waves corresponding to a contact discontinuity and a shock. The solution is self-similar with respect to the variable $\eta=x/t.$
(b) shock tube solution: density at $t=t_0$ (blue), density at $t=t_1$ (green) and density displacement interpolant at $t=(t_0+t_1)/2$ (yellow). The displacement interpolant coincides with the the exact solution.
}
\label{sod}
\end{figure}

The Riemann problem for conservation laws is a Cauchy problem with piece-wise constant initial data where a single discontinuity is placed at $x=0$ in the domain of interest. The problem is essentially one-dimensional and the solution is self-similar with respect to the self-similarity variable $\eta=x/t$ (see, e.g., \cite[Chapter 3]{toro2013riemann}). 
To fix the ideas, we consider the Sod shock tube problem, which corresponds to impose that the ideal gas is at rest ($u=0$) 
with high pressure and density for $x<0$ and   low pressure and density for $x>0$.  There are three waves emerging from the initial discontinuity: a receding expansion fan and two forward waves corresponding to a contact discontinuity and a shock, see Figures \ref{sod}(a) and \ref{sod}(b).

Density is a monotonically-decreasing function in $x$ and we have that $\rho(t,x)=\hat{\rho}(x/t)=\hat{\rho}(\eta)$, $\forall t>0$. As before, we can write:
$$
\frac{\partial \rho}{\partial x} = \frac{\partial \hat{\rho}}{\partial \eta}\frac{\partial \eta}{\partial x}.
$$
This means that,  up to a  change of sign, the density spatial derivative satisfies optimal transportation between any two times $t_0$ and $t_1$ (see \eqref{eq:local_mass_conservation})     with 
$\Psi(t_0,t_1,\xi)  = \frac{1}{2} \xi^2 \frac{t_1}{t_0}$. 
If we consider a linear scaling $s(t) = \frac{t-t_0}{t_1-t_0}$, we obtain 
$T(s(t),\xi) = \frac{t}{t_0} \xi$, which is the appropriate self-similarity transformation associated with the problem.

\section{Nonlinear interpolation based on Gaussian models}
\label{sec:method}

Although for other notable cases (such as steady boundary layers) the exact displacement interpolant coincides with the exact solution to  the physical model, in general it is not possible to readily identify an extensive scalar physical quantity that is representative of the whole solution field and for which the optimal transport map is available in closed form.

In the last few decades, there has been a growing interest in determining effective algorithms to approximate the solution to optimal transportation problems for arbitrary choices of the densities $\rho_0$ and $\rho_1$ in \eqref{eq:monge_problem},  
 \cite{peyre2019computational}. 
In this work, we pursue a different approach: first, we identify a Gaussian model $g[U]$ of the solution $U$ to the PDE; then, we exploit the knowledge in closed-form of the forward mapping $T$ between two standard Gaussians to define the  interpolation operator.
To illustrate the many features of our approach and its limitations  that motivate the extensions  of section \ref{sec:extension}, we present extensive numerical investigations for  one-dimensional and two-dimensional test problems with exact or numerical solutions.

\subsection{Methodology}

\subsubsection{Optimal transportation of multivariate normal density distributions}
\label{sec_opt_mult}
We briefly review the solution to \eqref{eq:monge_problem} for multivariate Gaussian densities; we refer to \cite{mccann1997convexity} for the proofs.
We define the normal density distribution $\phi$ with mean $\mu\in \mathbb{R}^n$ and symmetric positive definite covariance ${ {\Sigma }}\in\mathbb{R}^{n \times n}$: 
\begin{equation}
\label{eq:gaussian_distribution}
{\displaystyle \phi(x; \, \mu,\Sigma)={\frac {1}{(2\pi )^{n/2}\left|{ {\Sigma }}\right|^{1/2}}}\;e^{ -{\frac {1}{2}}\left({ {x}}-{ {\mu }}\right)^{\top }{ {\Sigma }}^{-1}\left({ {x}}-{{\mu }}\right)}}.
\end{equation}
Given the densities $\rho_0 = \phi(\cdot; \, \mu_0,\Sigma_0)$ and 
$\rho_1= \phi(\cdot; \, \mu_1,\Sigma_1)$,  we find that the  displacement interpolant $\widehat{\phi}_s$ is Gaussian with mean and covariance given by 
\begin{subequations}
\label{gaussian_transport}
\begin{equation}
\label{eq:mean_covariance}
\mu_s=(1-s)\,\mu_0+s\, \mu_1,
\quad
\Sigma_s=\Sigma_0^{-1/2}\left( \left(1-s\right) \Sigma_0 + s \left( \Sigma_0^{1/2} \Sigma_1 \Sigma_0^{1/2}\right)^{1/2} \right)^2     \Sigma_0^{-1/2},
\end{equation}
for all $s\in[0,1]$. The forward mapping $T$ is also available in closed form:
\begin{equation}
\label{forwaedm_gaussian}
  T(s,\xi)=(1-s)  \, \xi+s \, \left( \mu_1+ \Sigma_0^{-1/2}\left(  \Sigma_0^{1/2} \Sigma_1 \Sigma_0^{1/2}\right)^{1/2}     \Sigma_0^{-1/2} \left( \xi-\mu_0\right) \right).
\end{equation}
Finally, the Wasserstein distance between Gaussian density distributions is given by:
\begin{equation}
\label{wdist}
W_2\left( \phi\left(\mu_0,\Sigma_0\right),\phi\left(\mu_1,\Sigma_1\right)\right)=\sqrt{ \| \mu_1-\mu_0 \|_2^2 \, + \,    \Tr\left(\Sigma_0+\Sigma_1-2
\,     
 \left( \Sigma_0^{1/2} \Sigma_1 \Sigma_0^{1/2}\right)^{1/2}    \right)}.
\end{equation}
\end{subequations}

As a final remark, we note that the optimal mapping between Gaussian distributions is always well-defined, affine and can be obtained at negligible computational cost. 

\subsubsection{Gaussian models of coherent structures}
\label{sec:gauss_model_coherent}

Given the field $U: \mathbb{R}^n \to \mathbb{R}^d$, we define the scalar testing function $\mathcal{T}(\cdot; U): \mathbb{R}^n \to \mathbb{R}$ and the set 
\begin{equation}
\label{eq:coherent_set}
\mathcal{C}_{\mathcal{T}}(U) : =
\left\{
x \in \mathbb{R}^n \, : \, \mathcal{T}(x; U) > 0
\right\},
\end{equation}
which  identifies the coherent structure associated with the criterion $\mathcal{T}$. To provide a concrete example,
if $U$ is the velocity field, we might define 
 $\mathcal{T}(x; U) = \| \nabla \times U(x)  \|_2 - \tau$ with $\tau>0$: in this case,  
$\mathcal{C}_{\mathcal{T}}(U)$ identifies the region of the domain where the 
enstrophy exceeds a user-defined threshold.

 In order to fit a Gaussian model to 
 $\mathcal{C}_{\mathcal{T}}(U)$ in \eqref{eq:coherent_set}, we define a finite-dimensional discretization of the domain of interest $P_{\rm hf} = \{  x_i  \}_{i=1}^{N_{\rm hf}}$ and we define
\begin{equation}
\label{eq:discrete_coherent_set}
 P_{\rm hf}^+
 : = \left\{
 x\in P_{\rm hf}  \,: \,
 \mathcal{T}(x; U) > 0
  \right\} 
=
\left\{
y_j
\right\}_{j=1}^{N_{\rm hf}^+ }
\end{equation}
Then, the statistical parametric model of the coherent structure is obtained by assuming that $\left\{
y_j \right\}_j$ are independent identically distributed (iid) realizations of a multivariate Gaussian distribution, and then  resorting to maximum likelihood estimation (MLE, see, e.g., \cite[Chapter 8]{rice2006mathematical}) to estimate mean and variance:
\begin{equation}
\label{eq:gauss_model}
g(x; U) : = \phi \left( x; \mu_{\rm mle}[U],\Sigma_{\rm mle}[U] \right),
\quad {\rm where} \;
\left\{
\begin{array}{l}
\displaystyle{
\mu_{\rm mle}[U] = \frac{1}{N_{\rm hf}^+} \;
\sum_{j=1}^{N_{\rm hf}^+}   y_j,}
\\[3mm]
\displaystyle{
\Sigma_{\rm mle}[U] = \frac{1}{N_{\rm hf}^+ } \;
\sum_{j=1}^{N_{\rm hf}^+}   \,  ( y_j - \mu_{\rm mle}[U]) \, ( y_j - \mu_{\rm mle}[U])^T. }
\\
\end{array}
\right.
\end{equation}
{Note that the value of the scalar testing function is not used to weight the estimates of mean and variance in \eqref{eq:gauss_model}.
In all the examples considered, the scalar testing function  includes gradients of the solution field, which is   discontinuous, and is thus highly irregular and oscillatory. For this reason, the decision  not to  weight the points $\{y_j \}_j$ based on the values of $\mathcal{T}(\cdot; U)$ improves the robustness of the parameter estimation procedure and is less sensitive to the mesh size.
}

The  scalar testing function $\mathcal{T}$ identifies flow features  that we wish to track. From an approximation standpoint, it is natural to identify and then track high-gradient regions of the flow, which correspond to shock waves or contact discontinuities. In this respect, the testing function  $\mathcal{T}$ is related in scope to shock-capturing sensors that are used in high-order methods to activate numerical dissipation where needed, and also to error indicators used for mesh adaptation and refinement.  
In the framework of model reduction, we observe that 
we might also interpret the Gaussian model 
$U\mapsto g(\cdot ; U)$ as a  registration sensor: similarly to  \cite{taddei2021space}, $g(\cdot ; U)$ is indeed used to learn a suitable parametric  mapping $T_g$ that is ultimately used to approximate the parametric field of interest.

The use of Gaussian distributions allows to readily define the optimal mapping $T_g$ and ultimately the displacement interpolant, possibly at the price of inaccurate representations of the coherent structure of interest. The choice of the distribution model should be a compromise between  \emph{learnability} and \emph{expressivity}:
here, learnability can be measured in terms of the degree of difficulty of solving the subsequent optimal transportation problem, while expressivity is related to the difference in performance between the displacement interpolant based on Gaussian models and the displacement interpolant obtained by transporting the indicator function of $\mathcal{C}_{\mathcal{T}}(U)$. Note also that, for practical high-fidelity data, estimates of $\mathcal{C}_{\mathcal{T}}(U) $ might be noisy: our Gaussian model might thus also filter raw data and ultimately prevent over-fitting.

\subsubsection{Convex displacement interpolation}

Given the parametric field $U:[0,1]\times \mathbb{R}^n \to \mathbb{R}^d$, we consider the problem of constructing (nonlinear) interpolations between $U_0 = U(0,\cdot)$ and 
$U_1 = U(1,\cdot)$. Towards 
this end, we use the procedure in section
\ref{sec:gauss_model_coherent} to generate the Gaussian models $g_0,g_1$ and we use \eqref{forwaedm_gaussian} to compute the forward map $X_g$ and its inverse $Y_g = X_g^{-1}$ --- the latter is simply obtained by interchanging $U_0$ with $U_1$.
Then,  we define the CDI
$\widehat{U}:[0,1]\times \mathbb{R}^n\to \mathbb{R}^d$
of the form 
\eqref{eq:convex_displacement_interpolation}
 such that 
$$
\widehat{U}(s, x)
\, =  \,
(1 - s) U_0 \circ 
  W_g(s, x)  
 \, + \,
s U_1 \circ T_g (1-s,x),
\quad
s\in [0,1], x\in \mathbb{R}^n,
$$
where $W_g(s, x)= (1 - s) x + s Y_g(x)$.

 If multiple snapshots of $U$ are available for $0=t_0<\ldots<t_K=1$, we can improve the accuracy of \eqref{eq:convex_displacement_interpolation} by considering piecewise approximations in the intervals $A_k := ( t_k,t_{k+1})$, for $k=0,\ldots,K-1$, or by learning a more accurate rescaling function $s:[0,1]\to [0,1]$. {Regarding} the latter, we might (i) compute $s_1,\ldots,s_K$ such that
\begin{subequations}
\label{eq:convex_displacement_interpolation_multisnapshots}
\begin{equation}
\label{eq:learn_rescale}
s_k \in  {\rm arg} \min_{s\in [0,1]} \| \widehat{U}(s,\cdot) - U(t_k,\cdot) \|_{\star},
\quad
{\rm such \; that \;}
0 = s_0<s_1\ldots<s_K=1,
\end{equation}
where $\| \cdot \|_{\star}$ is a functional norm of interest;
(ii) compute a bijective  rescaling $\widehat{s}:[0,1]\to [0,1]$
based on the dataset $\{(t_k,s_k)\}_{k=0}^K$ using a standard regression algorithm; and
(iii)  define the interpolant:
\begin{equation}
\label{eq:convex_displacement_interpolation_plus}
\widehat{U}'(t, x)
\, =  \,
\widehat{U}(\widehat{s}(t), x),
\quad
t\in [0,1], \; x\in \mathbb{R}^n.
\end{equation}
\end{subequations}
In the numerical examples, we show that optimizing the rescaling function $s$ might have a significant impact on performance; furthermore, it might  unveil relevant features of the coherent structure of interest.

\subsection{Numerical examples}
\label{sec:numerics_examples}

\subsubsection{Simple wave field}
\label{sec:numerics_simple_wave}
We study the problem described in section \ref{sec:simple_wave}. In this case, the parametric evolution of the solution is considered with respect to time. As discussed in the previous sections, optimal transportation exactly maps the initial condition to subsequent solution profiles. Here, we compare the empirical similarity transform determined based on two solution snapshots at $t_0$ and $t_1$, to the exact time-dependent solution. 

\begin{figure}[H]
    \centering
    \subfloat[ ]{
    \includegraphics[width=.45\textwidth]{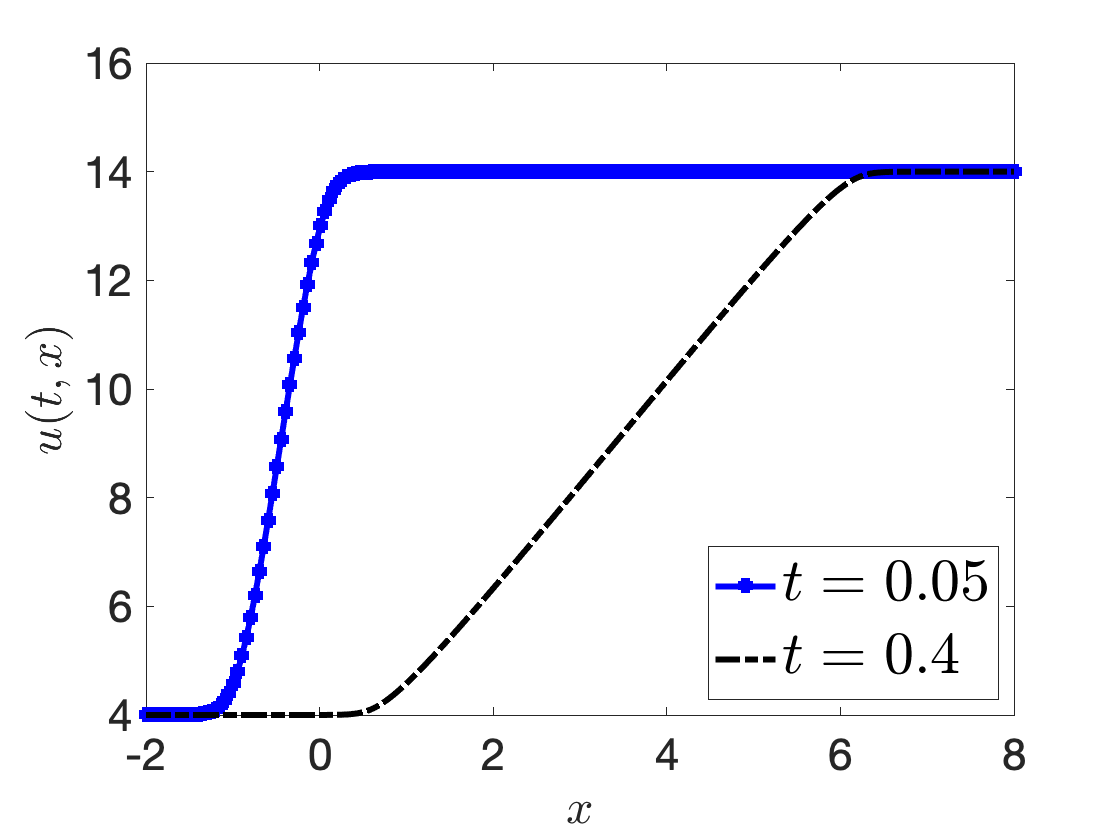}}
    ~~
    \subfloat[ ]{
    \includegraphics[width=.45\textwidth]{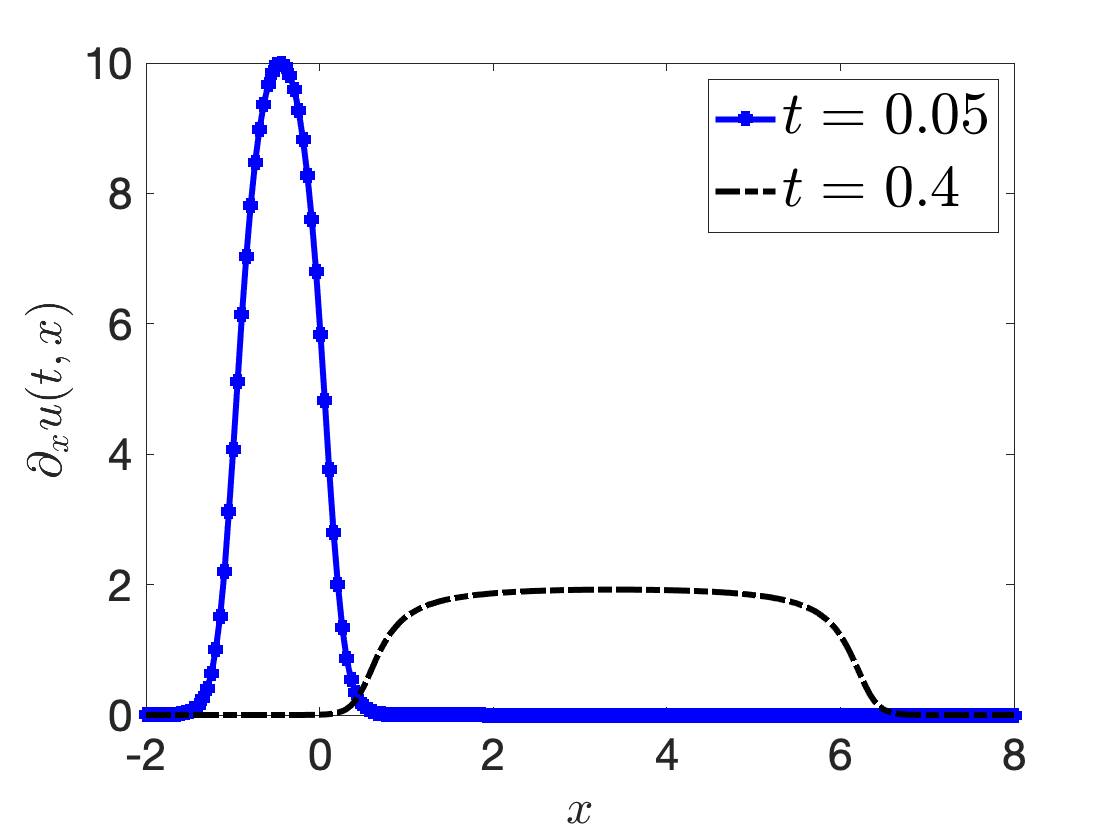}
    }
    \caption{
simple wave field.    
(a) solution velocity snapshots $u(t,x)$ at time $t_0=0.05$ (blue) and $t_1=0.4$ (dashed black). 
(b) space derivative of the velocity at time $t_0=0.05$ (blue) and $t_1=0.4$ (dashed black).}
    \label{fig_sw1}
\end{figure}

Let the initial condition for the speed of sound be  $a(0,x)=2 + \tanh{((x + 1)/0.2)}$, the initial condition for the left-going Riemann invariant $R^-=1$ and $\gamma=7/5$. The right-going characteristics are hence straight lines and the solution is an expansion fan traveling rightward. As an example, 
in Figure \ref{fig_sw1}
we show two snapshots of the velocity field $u(t,x)$ at times $t_0=0.05$ and $t_1=0.4$.

\begin{figure}[h!]
    \centering
    \subfloat[$t_0=0.05$]{
    \includegraphics[width=.45\textwidth]{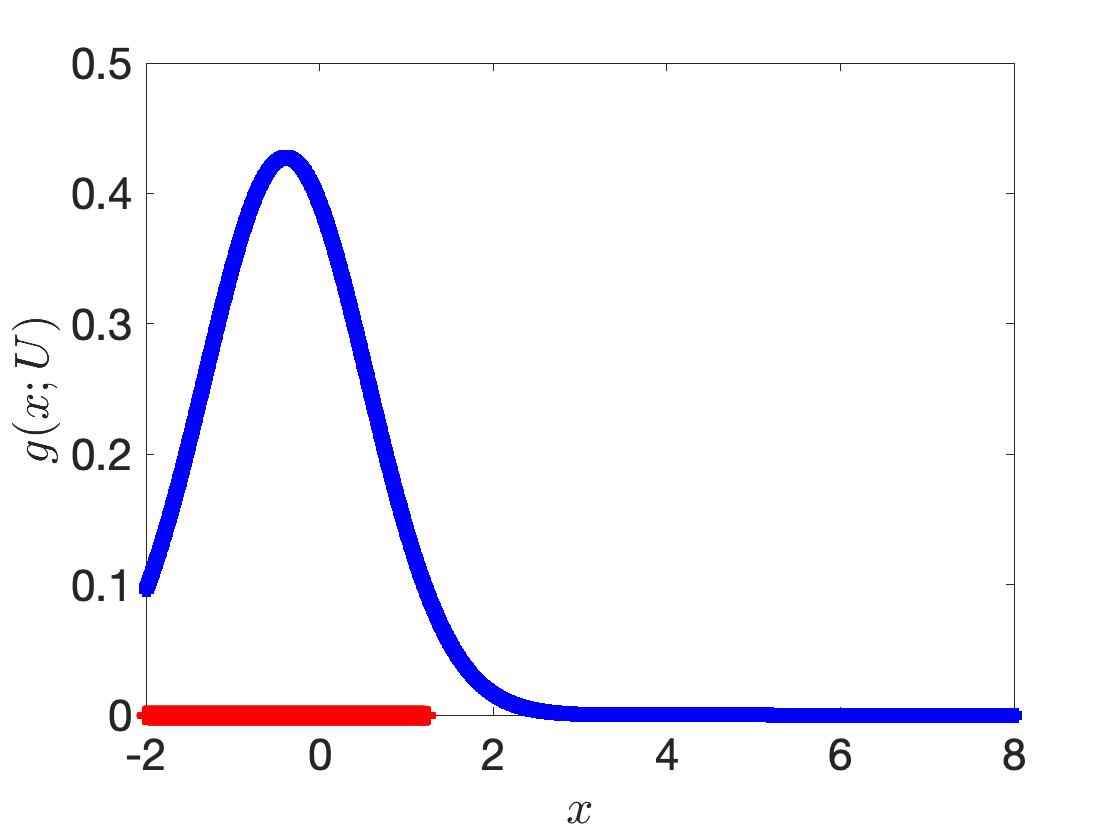}
    }
    ~~
    \subfloat[$t_1=0.4$]{
    \includegraphics[width=.45\textwidth]{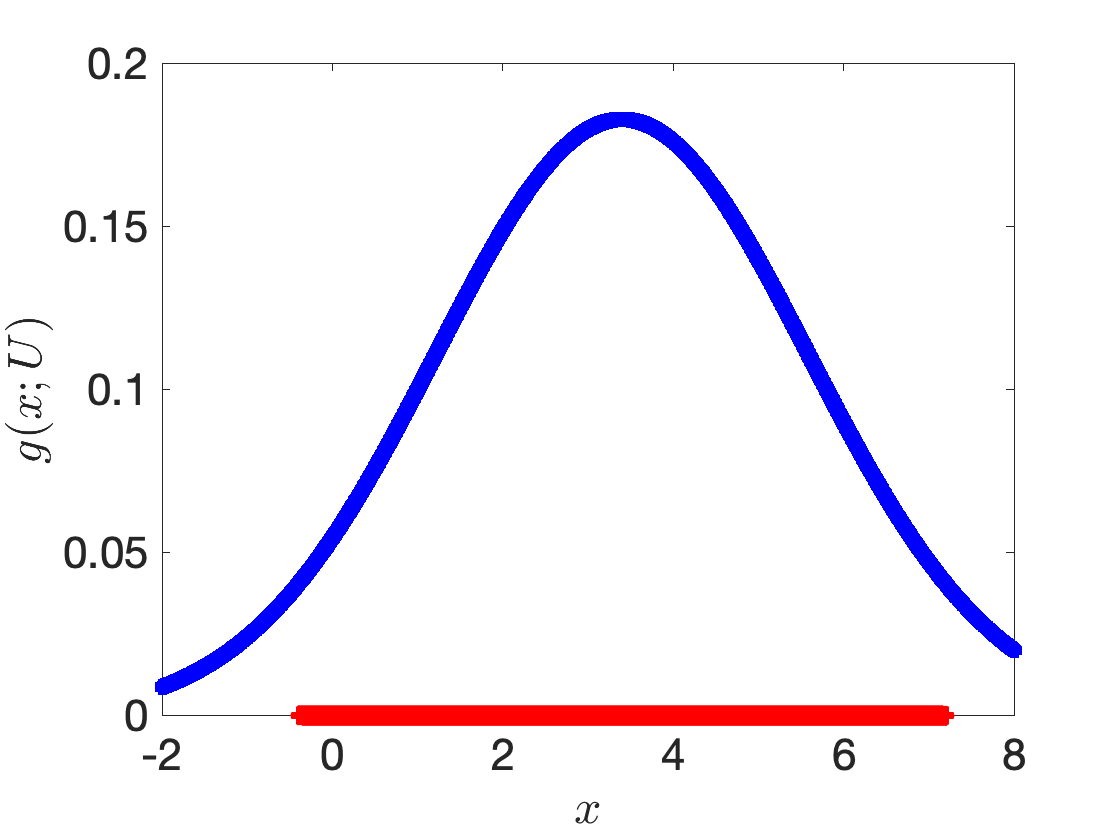}
    }
    \caption{
simple wave field. 
(a)-(b) MLE Gaussian density distributions (blue) 
at times  $t_0=0.05$ and $t_1=0.4$. Red points indicate the elements of $P_{\rm hf}^+$.}
    \label{fig_sw2}
\end{figure}

We consider the scalar testing function
\begin{equation}
\label{eq:test_function_sw}
\mathcal{T}(x; U) : = 
\left| \frac{\partial u}{\partial x} (x)\right|
- \epsilon,
\end{equation}
with $\epsilon=10^{-4}$:   a point $x \in \mathbb{R}$ belongs to $C_{\cal T}$ if the absolute value of the space derivative of the velocity field $u(x)$ is larger than $\epsilon$. Figure \ref{fig_sw2} shows the MLE
Gaussian density distributions (blue) 
at times  $t_0=0.05$ and $t_1=0.4$; the  red points indicate the elements of the sets $P_{\rm hf}^+$ \eqref{eq:discrete_coherent_set}.

We compare the exact velocity solutions at $t\in[0.05,0.4]$ to their  $L^2$ projection in the manifold of  displacement interpolants $\{ \widehat{u}_s : s\in[0,1] \}$, cf.  eq.
\eqref{eq:convex_displacement_interpolation}.  More precisely, given  $\alpha\in[0,1]$,  we define $s_{\alpha}$ such  that
\begin{equation}
\label{eq:s_alpha}
s_{\alpha} : = 
\text{arg}\,\min_{s}\,   \int_{\mathbb{R}} 
\left( \widehat{u}_s(x) - u((1-\alpha)t_0+\alpha t_1, x) \right)^2\,dx,
\end{equation}
and the projection $\widehat{u}_{{\alpha}}^{\rm proj} = \widehat{u}_{s_{\alpha}}$.
Similarly, we project the exact solution in the convex set spanned by the two solutions at time $t_0$ and $t_1$ and we define
\begin{equation}
\label{eq:s_alpha_co}
\begin{array}{l}
\displaystyle{
s_{\alpha}^{\rm co} : = 
\text{arg}\,\min_{s}\,   \int_{\mathbb{R}} 
\left( \widehat{u}_s^{\rm co}(x) - u((1-\alpha)t_0+\alpha t_1 , x)\right)^2\,dx,
}
\\[3mm]
\displaystyle{
{\rm with}  \; \; 
\widehat{u}_s^{\rm co}(x)  = (1-s) u(t_0,x) + s\, u( t_1, x). 
}\\
\end{array}
\end{equation}

\begin{figure}[h!]
\centering
\subfloat[]{
\includegraphics[width=.495\textwidth]{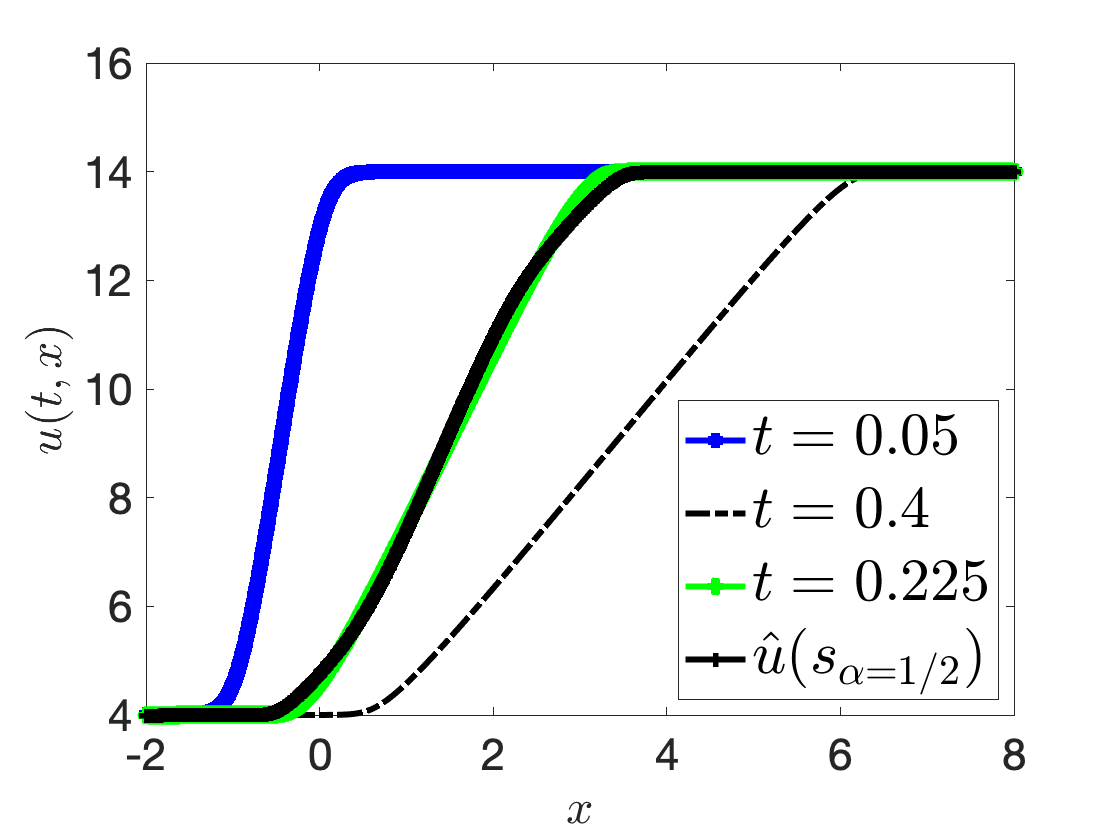}}
~~~
\subfloat[]{
\includegraphics[width=.495\textwidth]{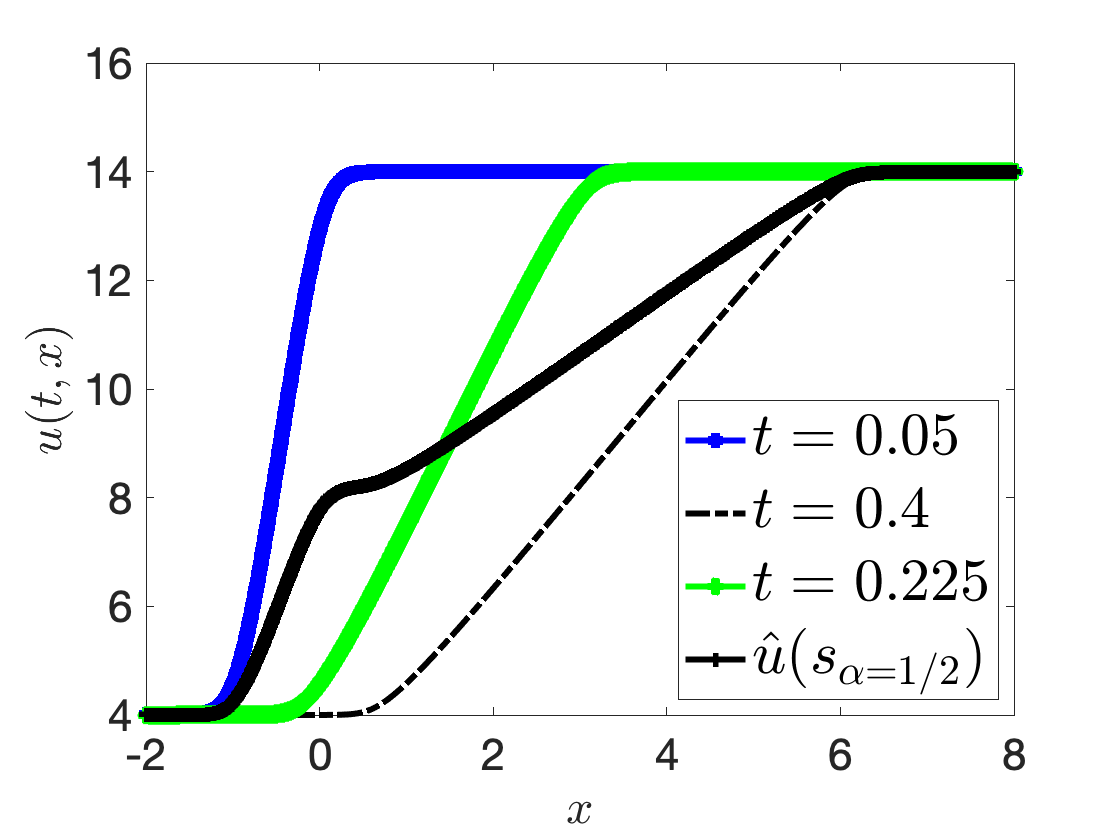}
}
\caption{simple wave field.
 (a)  velocity profile at $t_0$ (blue), $t_1$ (dashed black), exact solution for $\frac{t_0+t_1}{2}$ (green) and  convex displacement interpolant  $\widehat{u}_s$ for    $s=0.7$ (black dots). 
(b)  
convex interpolant  $\widehat{u}_s^{\rm co}$ for    $s=0.58$ (black  dots).}
    \label{fig_sw3}
\end{figure}

 In  Figure \ref{fig_sw3},  we compare the exact velocity profile for $\alpha=0.5$ (i.e., $t = 0.225$) to the optimal  displacement interpolant $\widehat{u}_{s_\alpha}$ and  to the convex interpolation $\widehat{u}_{s_\alpha^{\rm co}}^{\rm co}$. We find  that $s_\alpha=0.70$ and $s_\alpha^{\rm co}=0.58$. The displacement interpolant captures the essential features of the solution while the convex projection in the convex set of the initial and final snapshots is completely inaccurate.
 
\begin{figure}[h!]
    \centering
            \subfloat[]{
    \includegraphics[width=.495\textwidth]{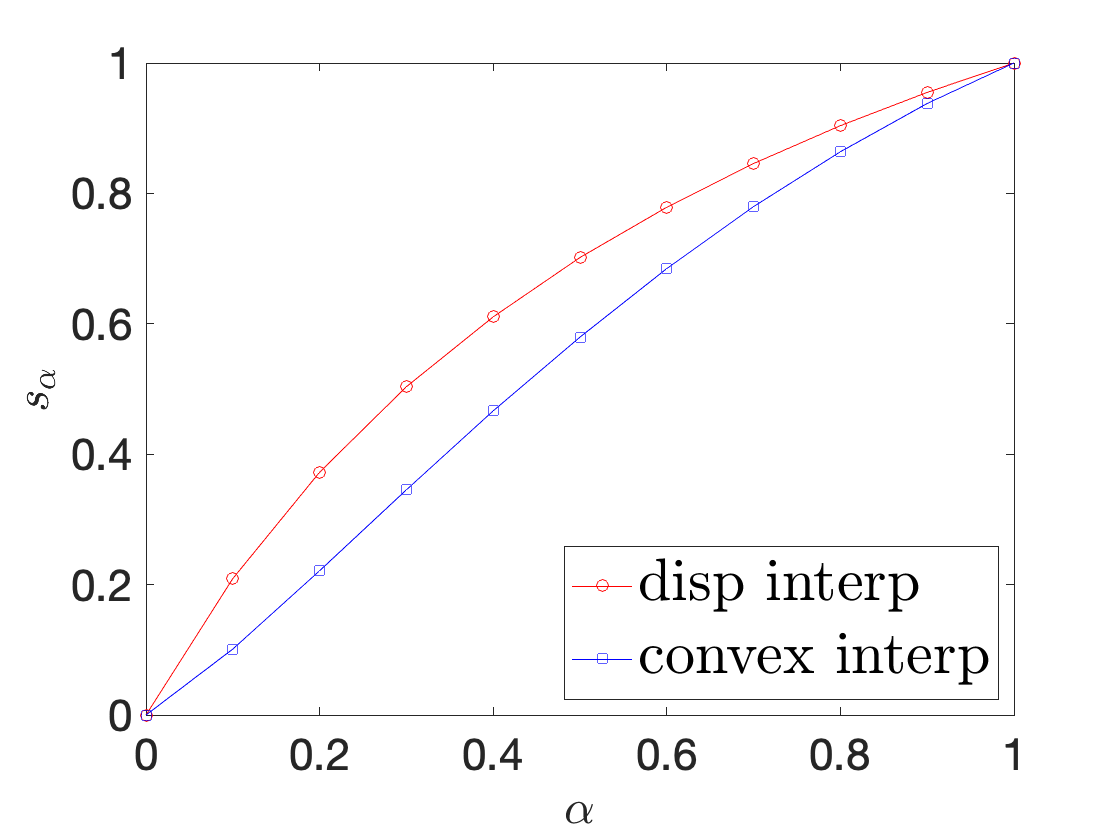}
    }
    ~~
            \subfloat[]{
    \includegraphics[width=.495\textwidth]{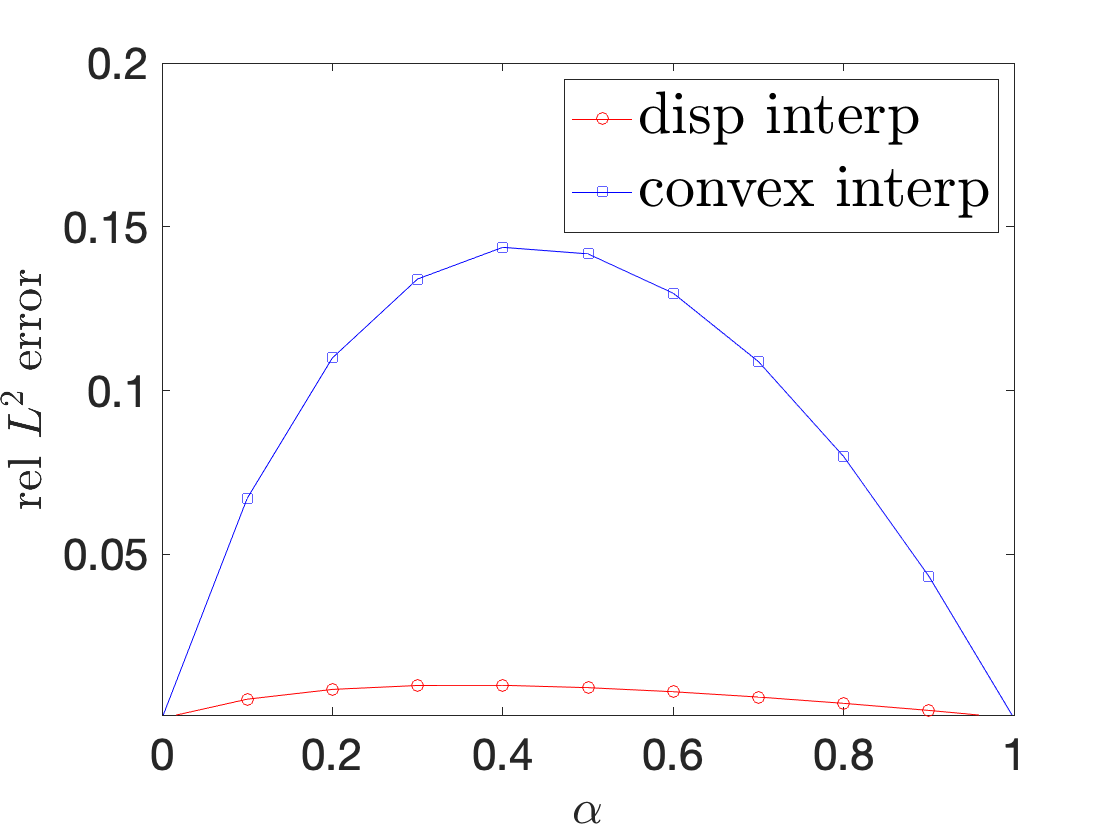}
    }
\caption{
simple wave field.
(a) behavior of $s_\alpha$ (disp interp) and  ${s}_\alpha^{\rm co}$   (convex interp)  with respect to  $\alpha$.
(b) behavior of the relative $L^2$ error in $(-2,8)$.
}
\label{fig_sw4}
\end{figure}

In Figure \ref{fig_sw4},  we show the behaviors  of $s_\alpha$ and  ${s}_\alpha^{\rm co}$  with respect to  $\alpha$. As expected from the motivating examples in section \ref{sec:motivating_examples}, $s_\alpha$ is not necessarily linear with respect to $\alpha$. In the same figure,  we show the relative $L^2$ projection error of the exact solution with respect to the displacement interpolant $\widehat{u}$ and with respect to  the convex interpolant $\widehat{u}^{\rm co}$, for several values of $\alpha$. These results show that even for a smooth solution with a non-compact support, convex displacement interpolation  systematically improves the approximations with respect to the convex projection in the space of the snapshots.

\subsubsection{Supersonic flow past a wedge}
\label{sec:wedge}
\begin{figure}[h!]
\centering
\begin{tikzpicture}[scale=1.1]
\draw [black,thick] (0,0) -- (3,0) -- (5,1) -- (5,3) --(0,3)--(0,0);
 \draw [->,black,ultra thick] (-0.3,0.5) -- (0,0.5);
  \draw [->,black,ultra thick] (-0.3,1.5) -- (0,1.5);
   \draw [->,black,ultra thick] (-0.3,2.5) -- (0,2.5);
   
 \draw [->,black,ultra thick] (0,-0.1) -- (1,-0.1);

\coordinate [label={right:  {\Large {$\Omega_{\rm f}$}}}] (E) at (5,2.5) ;
\draw [->,black,ultra thick] (0,-0.1) -- (1,-0.1);
\draw [black,ultra thick] (0,-0.2) -- (0,0);

\draw[white,pattern=north west lines,pattern color=red] 
(0,0) -- (3,0)--(3,-0.25) --(0,-0.25) -- (0,0);

  \draw[white,pattern=north west lines,pattern color=red] 
(3,0) -- (5,1)--(5,0.75) --(3,-0.25) -- (3,0);

\draw[ultra thick] (4,0) arc (0:26:1);
\draw [black,ultra thick] (3,0) -- (4,0);

\draw [blue,ultra thick,dashed] (3,0) -- (5,2);
\draw[blue,ultra thick] (4.5,0) arc (0:45:1.56);
\coordinate [label={right:  {\Large {$\theta$}}}] (E) at (4.5,0.2) ;

\coordinate [label={right:  {\Large {$\delta$}}}] (E) at (4,0) ;
\coordinate [label={left:   {\Large {$M^{\rm u}$}}}] (E) at (0,2) ;
\coordinate [label={right:  {\Large {$M^{\rm d}$}}}] (E) at (5,1) ;

\end{tikzpicture}
\caption{supersonic flow past a wedge, problem setting.
$M^{\rm u}$ is the free-stream Mach number, $\delta$ is the wedge angle, $\theta$ is the shock angle. Wall conditions are imposed on the bottom boundary.}
\label{fig:supersonic_wedge_explanation}
\end{figure}
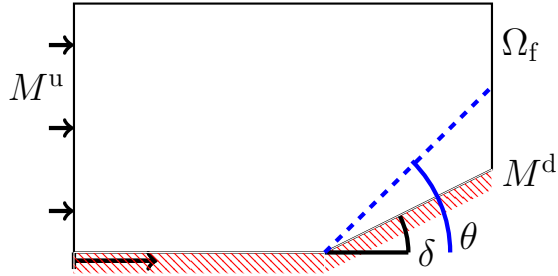

We consider a two-dimensional compressible Euler flow of air ($\gamma=7/5$) past a wedge. The upstream flow is supersonic and it induces a steady attached shock wave that develops from the leading edge if the upstream Mach number is within a given range, which depends on the wedge angle. In this test case we let the solution vary with respect to the upstream Mach number $M^{\rm u}$ and  the wedge angle $\delta$. 
We compare below the convex displacement interpolant obtained by the empirical similarity transform to the exact solution. As an example, we study the interpolation between $M_0^{\rm u}=5, \delta_0=28.275$  and
$M_1^{\rm u}=8, \delta_1=22.80$.
{Figure \ref{fig:supersonic_wedge_explanation} illustrates the system configuration.}

We denote by $\Omega_{\rm f}$
the physical domain, and we denote by
$(x_1,x_2) \mapsto M(x_1,x_2)$ the Mach number.
In view of the discussion, we introduce the reference domain $\Omega_{\rm r} = (-0.5,1)\times (0,1)$ and the geometric transformation 
$\Lambda: \Omega_{\rm r} \times (-\pi/2, \pi/2) \rightarrow \Omega_{\rm f} \subset \mathbb{R}^2$, $(x_1,x_2,\delta) \mapsto \Lambda(x_1,x_2,\delta)$,  such that
\begin{subequations}
\begin{equation}
\label{eq:Lambda_wedge}
\Lambda(x_1,x_2,\delta)= 
\left\{
\begin{array}{ll}
(x_1,x_2) &  x_1<0 \\[3mm]
\left(x_1, \displaystyle 
x_1 \tan (\delta) +     (1 - x_1 \tan (\delta) ) x_2    
\right)
        & x_1 \geq 0.\\
    \end{array}
\right.
\end{equation}
and its inverse
\begin{equation}
\label{eq:Theta_wedge}
\Theta(x_1,x_2,\delta)= 
\left\{
    \begin{array}{ll}
        (x_1,x_2) &  x_1<0 \\[3mm]
        \left(x_1, \displaystyle 
\frac{ x_2-x_1 \, \tan(\delta) }{
1 - x_1 \tan (\delta) }
\right)
        & x_1 \geq 0.
    \end{array}
\right.
\end{equation}
Finally, we define the mapping
$(x_1,x_2,\delta) \mapsto \Phi(x_1,x_2,\delta)$ such that
\begin{equation}
\label{eq:Phi_wedge}
\Phi(x_1,x_2,\delta)= 
\Lambda \left(
\Theta \left( x_1,x_2,\bar{\delta}\right),
\, \delta
\right),
\qquad
\bar{\delta} = \frac{\delta_0 +\delta_1}{2},
\end{equation}
which maps $\Omega_{\rm f} (  \bar{\delta})$ into $\Omega_{\rm f}(\delta)$. 
\end{subequations}

It is possible to show that the exact solution is piecewise-constant and exhibits a straight  shock discontinuity that is generated at the wedge  leading edge. If  we denote by $\theta$ the shock angle  and by 
$M^{\rm d}$ the downstream Mach number, we obtain the expression for $M$: 
\begin{subequations}
\label{eq:wedge_exact}
\begin{equation}
\label{eq:wedge_exact_a}
M(x_1,x_2) = \;
\left\{
\begin{array}{ll}
M^{\rm u} & {\rm if} \; x_2 > x_1 \tan(\theta) \\[3mm]
M^{\rm d} & {\rm if} \; x_2 < x_1 \tan(\theta) \\
\end{array}
\right.
\end{equation}
Given $M^{\rm u}$ and $\delta$, we can employ the relationships (cf. 
\cite{staff1953equations})
\begin{equation}
\label{eq:wedge_exact_b}
\begin{array}{l}
\displaystyle{
\cot \left(   \delta \right) \, = \, 
\tan \left(   \theta \right) \,
 \left(   
 \frac{(\gamma+1) (  M^{\rm u}   )^2  }{ 2 ( M^{\rm u} \sin (\theta) )^2 - 1       } \, - \,
 1   \right),
}
\\[3mm]
\displaystyle{
 \left( M^{\rm d} \sin ( \delta -\theta) \right)^2
 \, = \,
 \frac{ (\gamma - 1)   \left( M^{\rm u} \sin ( \theta) \right)^2 + 2}{
 2 \gamma  \left( M^{\rm u} \sin ( \theta) \right)^2 - (\gamma - 1)
 },
}
\\
\end{array}
\end{equation}
\end{subequations}
to find the downstream Mach number $M^{\rm d}$ and the shock angle $\theta$. 

In order to deal with geometry variations, we   pursue two different strategies.
\begin{enumerate}
\item
\emph{Extension:} we  extend the Mach number $M$ to $\mathbb{R}^2$ for all parameters.
\item
\emph{Geometric registration:}
we apply the interpolation procedure to the mapped field $\widetilde{M}(x_1,x_2; M^{\rm u}, \delta) = 
{M}(  \Phi(x_1,x_2; \delta) ; M^{\rm u}, \delta)$, which is defined in $\Omega_{\rm f}(\bar{\delta})$ for all values of the parameters $M^{\rm u}, \delta$. 
\end{enumerate}
Note that since the proposed displacement interpolation strategy does not preserve boundaries, extension outside $\Omega_{\rm f}$ is  necessary for both techniques; note also that, since the solution is piecewise-constant and the shock curve is linear,  the extension   is straightforward.

\begin{figure}[H]
\centering
\subfloat[ ]{
\includegraphics[width=.45\textwidth]{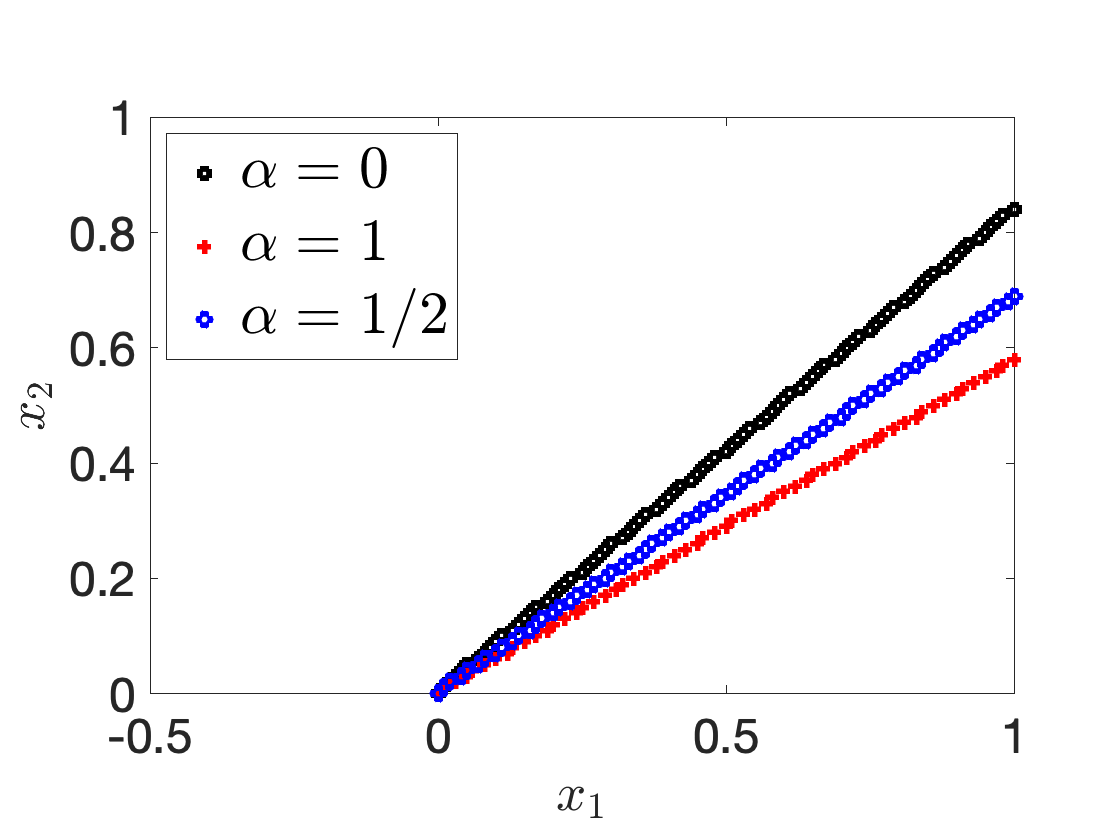}}
~~
\subfloat[ ]{
\includegraphics[width=.45\textwidth]{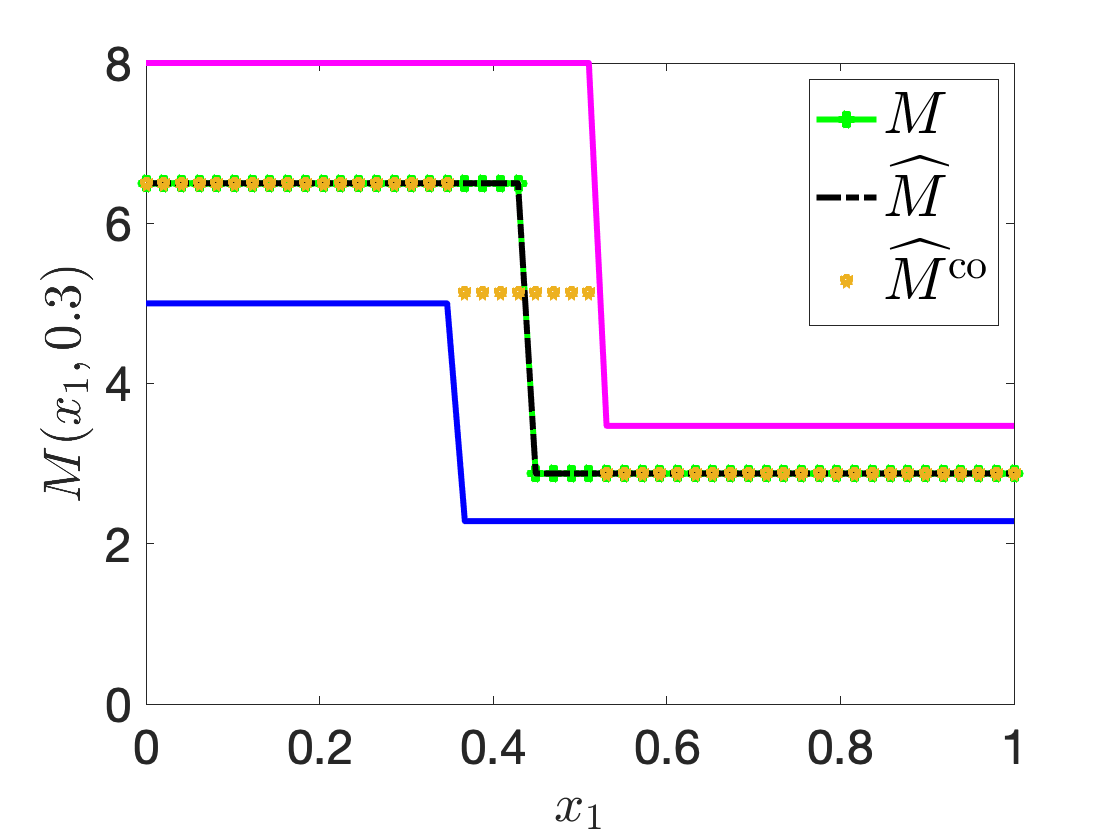}
}
    
\subfloat[ ]{
\includegraphics[width=.45\textwidth]{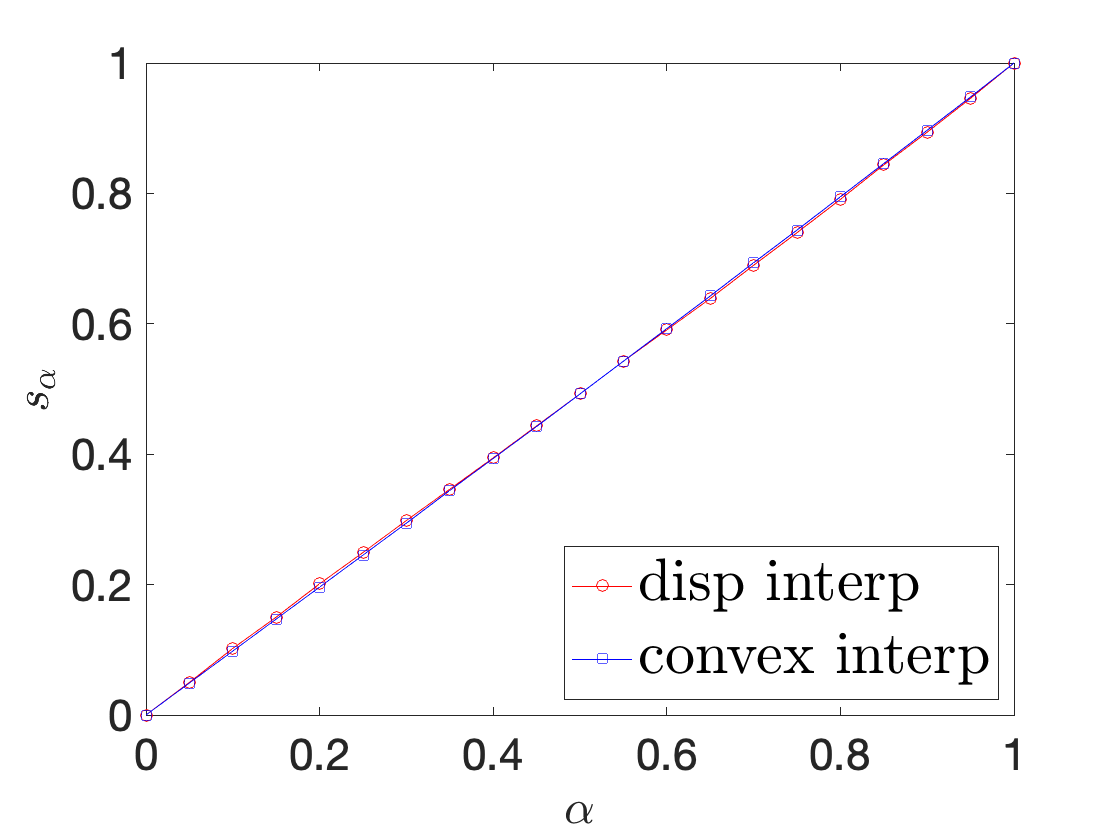}}
~~
\subfloat[ ]{
    \includegraphics[width=.45\textwidth]{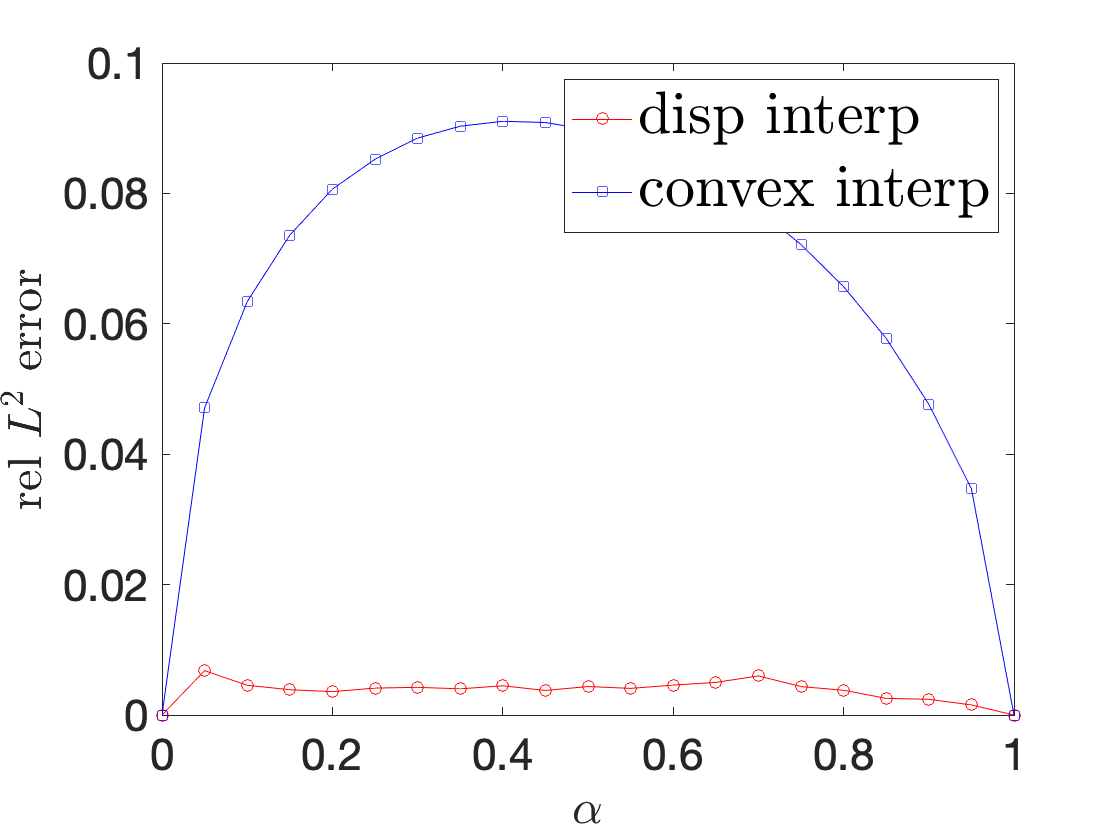}
    }
    
\caption{supersonic flow past a wedge; approach  based on extension.
(a) selected points $P_{\rm hf}^+$ for three choices of the parameter pair   $\mu_{\alpha}$;
(b) Mach profile for $x_2=0.3$ at  $\mu_0$ (blue), 
$\mu_1$ (violet), $\mu_{1/2}$ (green),  and
CDI  $\widehat{M}$ and 
convex interpolant $\widehat{M}^{\rm co}$ for $s=1/2$;
(c) behavior of $s_{\alpha}$  and  $s_{\alpha}^{\rm co}$  in \eqref{eq:s_alpha} and  \eqref{eq:s_alpha_co};
(d) behavior of the relative $L^2$ projection error.}
    \label{fig:approach1_wedge}
\end{figure}

We pursue the first approach based on extension.  Towards this end, we consider a regular $151 \times 101$ grid 
 in the rectangle $\Omega_{\rm r}$ and we consider the scalar testing function
\begin{equation}
\label{eq:scalar_testing_fun_wedge}
\mathcal{T}(x; U) =
\frac{1}{\epsilon} \big|M(x_1 + \epsilon, x_2) 
-  M(x_1, x_2)  \big| 
 \,  - \,  1,
\end{equation}
 where $\epsilon=10^{-2}$ is equal to the size of the grid. 
  Figure \ref{fig:approach1_wedge} shows the 
  results:
in Figure  \ref{fig:approach1_wedge}(a), we show  the selected points for three choices of the parameter pair   $\mu_{\alpha} : =  (1 - \alpha) (M_1^{\rm u}, \delta_1)+\alpha 
(M_1^{\rm u}, \delta_1)$, 
$\alpha \in \{ 0, 1/2, 1  \}$;
in  Figure  \ref{fig:approach1_wedge}(b), we compare the exact
Mach profile for $x_2=0.3$ and $x_1\in (0,1)$ for  $\mu_{\alpha=1/2}$ with the optimal CDI $\widehat{M}_{s_{\alpha}}$ and the convex interpolation $\widehat{M}_{s_{\alpha}}^{\rm co}$; 
in Figure  \ref{fig:approach1_wedge}(c), we show the behavior of $s_{\alpha}$  and  $s_{\alpha}^{\rm co}$ defined as  in \eqref{eq:s_alpha} and  \eqref{eq:s_alpha_co}, respectively;
in Figure  \ref{fig:approach1_wedge}(d), we show the behavior of the $L^2$ relative projection error.  Note that the optimal value of $s$ is a linear function of $\alpha$; note also that displacement interpolation offers extremely accurate results compared to the  more standard convex interpolation.
 
\begin{figure}[H]
\centering
\subfloat[ ]{
\includegraphics[width=.45\textwidth]{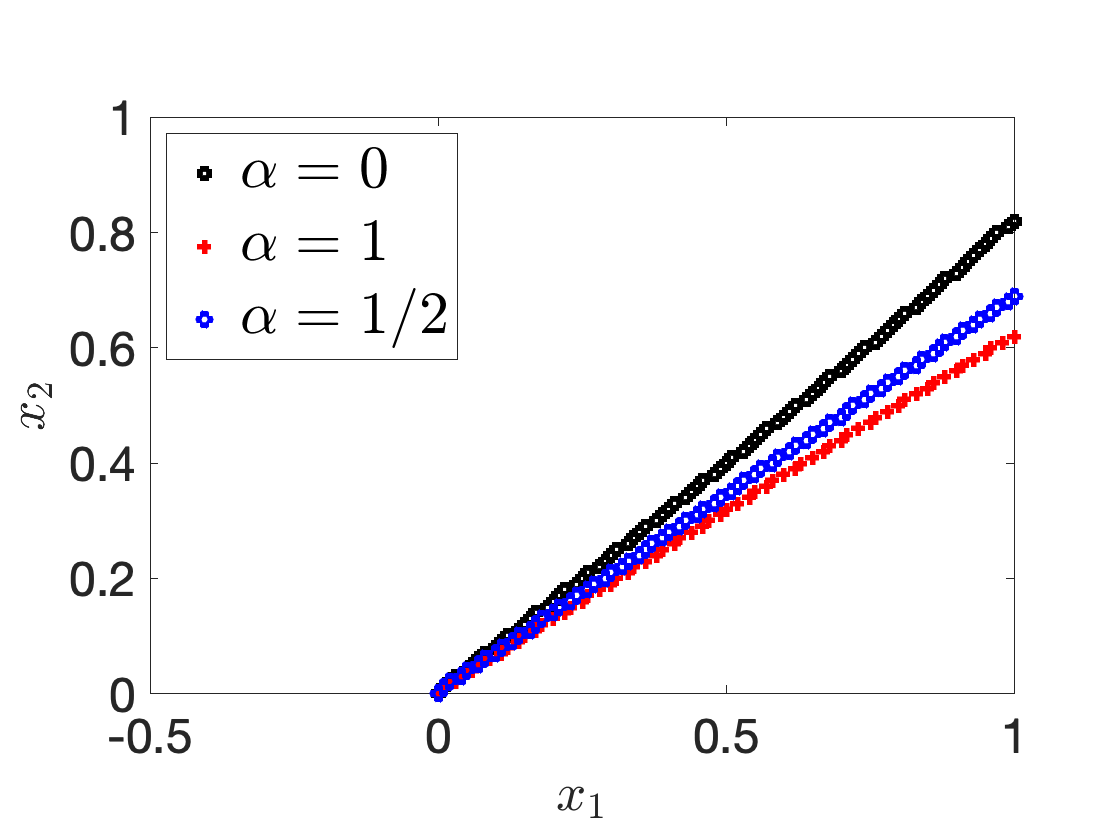}}
~~
\subfloat[ ]{
\includegraphics[width=.45\textwidth]{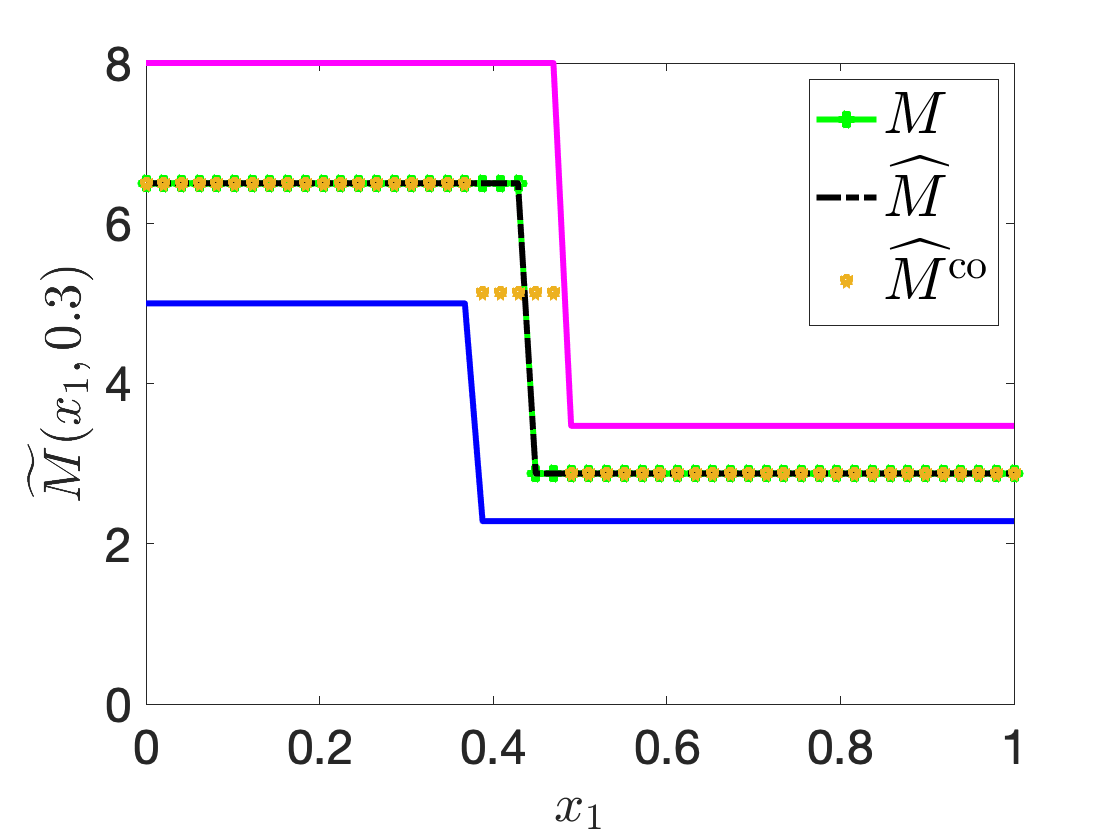}
}
    
\subfloat[ ]{
\includegraphics[width=.45\textwidth]{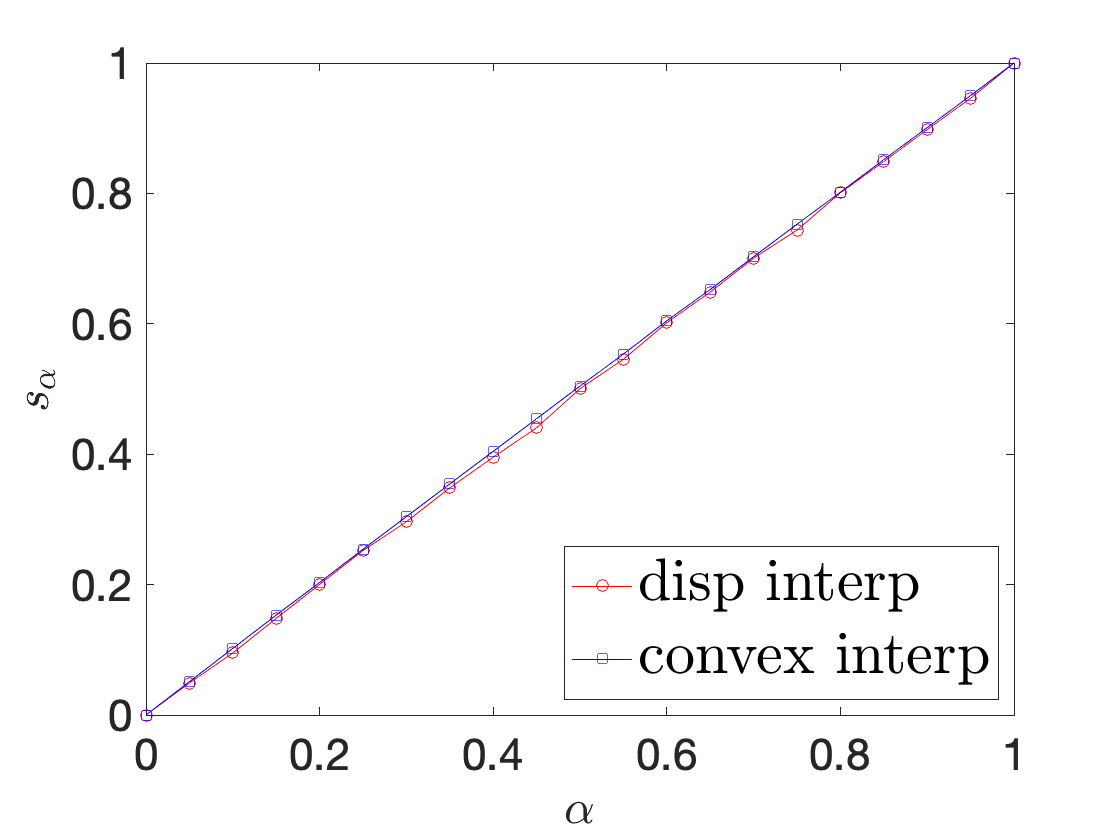}}
    ~~
\subfloat[ ]{
\includegraphics[width=.45\textwidth]{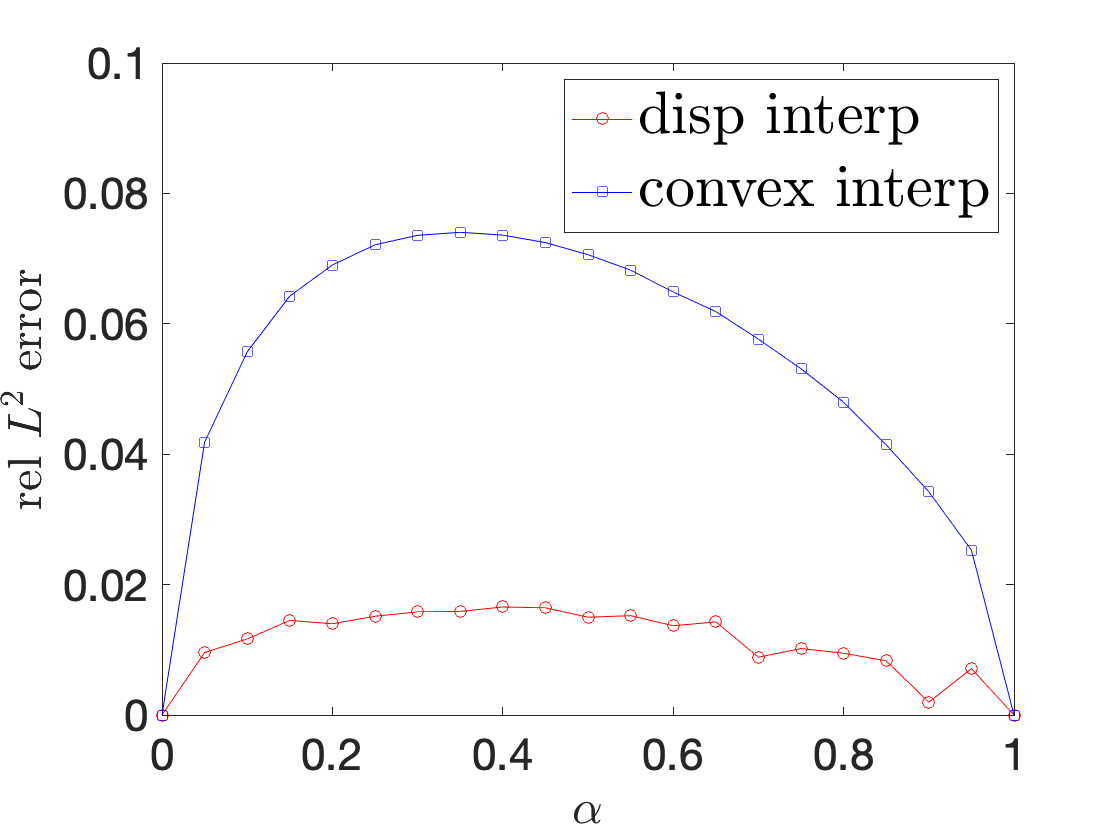}
}
    
\caption{
supersonic flow past a wedge; approach based on registration.
(a) selected points $P_{\rm hf}^+$ for three choices of the parameter pair   $\mu_{\alpha}$;
(b) Mach profile for $x_2=0.3$ at  $\mu_0$ (blue), 
$\mu_1$ (violet), $\mu_{1/2}$ (green),  and CDI  $\widehat{M}$ and 
convex interpolant $\widehat{M}^{\rm co}$ for $s=1/2$;
(c) behavior of $s_{\alpha}$  and  $s_{\alpha}^{\rm co}$  in \eqref{eq:s_alpha} and  \eqref{eq:s_alpha_co};
(d) behavior of the relative $L^2$ projection error.}
    \label{fig:approach2_wedge}
\end{figure}

We also pursue the second approach based on geometric registration.  Towards this end, we consider the same  regular $151 \times 101$ grid
 in the rectangle $\Omega_{\rm r}$,  but we  discard points outside $\Omega_{\rm f}(\bar{\delta})$; then, we consider the scalar testing function
\begin{equation}
\label{eq:scalar_testing_fun_wedge_2}
\mathcal{T}'(x; U) =
\frac{1}{\epsilon} \big|\widetilde{M}(x_1 + \epsilon, x_2) 
-  \widetilde{M} (x_1, x_2)  \big| 
 \,  - \,  1,
 \quad
 \epsilon=10^{-2}.
\end{equation}
Figure \ref{fig:approach2_wedge} replicates the same tests considered for the other strategy: as for the  previous approach, displacement interpolation significantly outperforms linear interpolation. We note, however, that the geometric mapping has a beneficial effect on the performance of the linear approach, while it is slightly detrimental for displacement interpolation.

\subsubsection{Transonic flow past an airfoil}
 \label{sec:transonic_flow}
 
\begin{figure}[H]
\centering
\subfloat[$M=0.77$]{
\includegraphics[width=.45\textwidth]{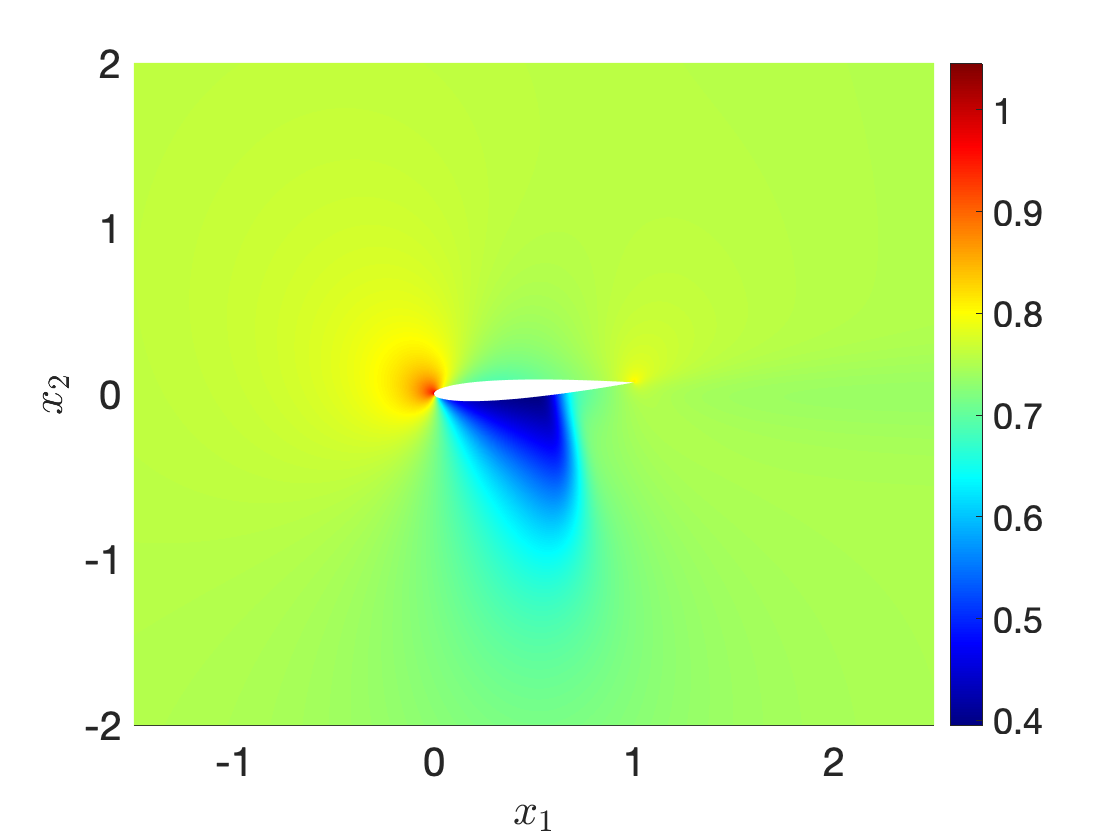}}
~~
\subfloat[$M=0.83$]{
\includegraphics[width=.45\textwidth]{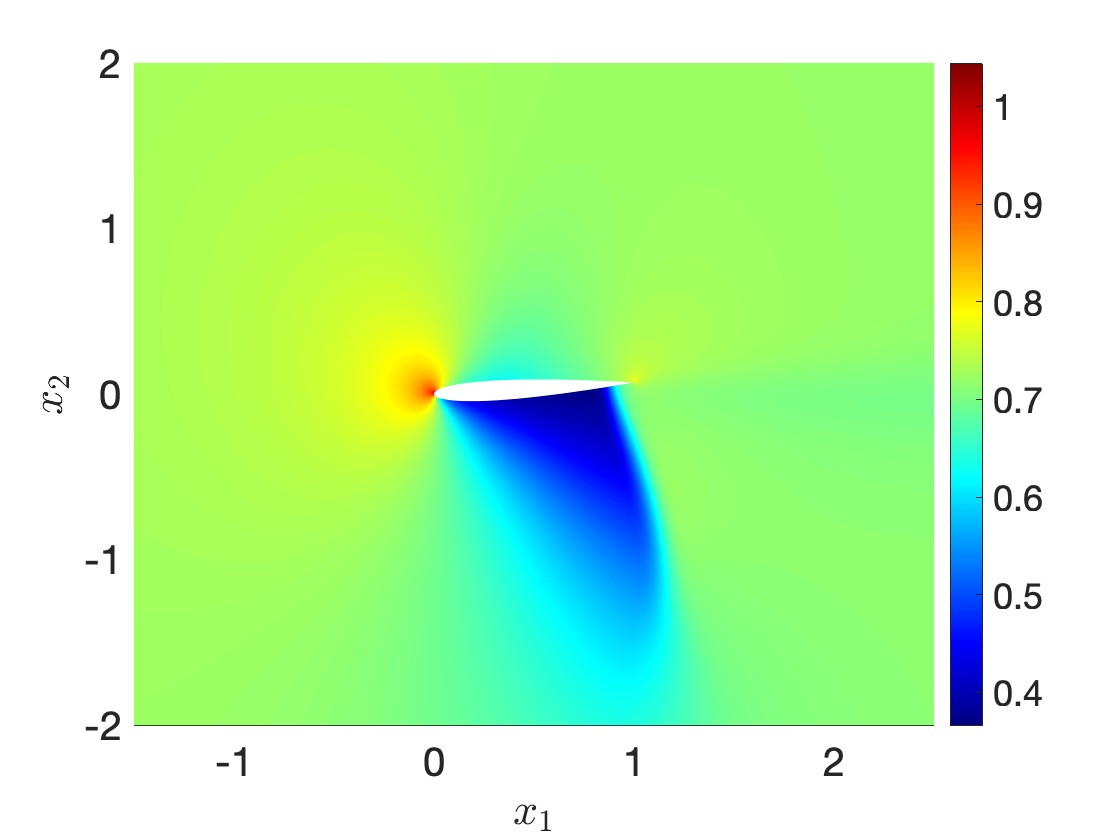}
}
        
\caption{transonic flow past an airfoil. Flow density for two values of the Mach number.}
\label{fig:airfoil_density}
\end{figure}

We consider a two-dimensional transonic flow past a NACA 0012 airfoil at angle of attack $\alpha=-4^o$; we let the solution vary  with respect to the free-stream Mach number $M \in  [0.77,0.83]$; we study the interpolation between $M_0=0.77$ and $M_1=0.83$. 
Related examples are considered in \cite{riffaud2021dgdd,taddei2021registration}.
We resort to a discontinuous Galerkin (DG) discretization with artificial viscosity to estimate the solution field; computations are performed in the domain 
$\Omega = (-4,10) \times (-10, 10) \setminus \Omega_{\rm naca}$ where $\Omega_{\rm naca}$ is the domain associated with the airfoil. Figure \ref{fig:airfoil_density} shows the behavior of the flow density for $M=0.77$ and 
 $M=0.83$: note that the solution 
develops a shock on the lower boundary of the airfoil that 
 is extremely sensitive to the value of the Mach number.

\begin{figure}[H]
\centering
\subfloat[ ]{
\includegraphics[width=.45\textwidth]{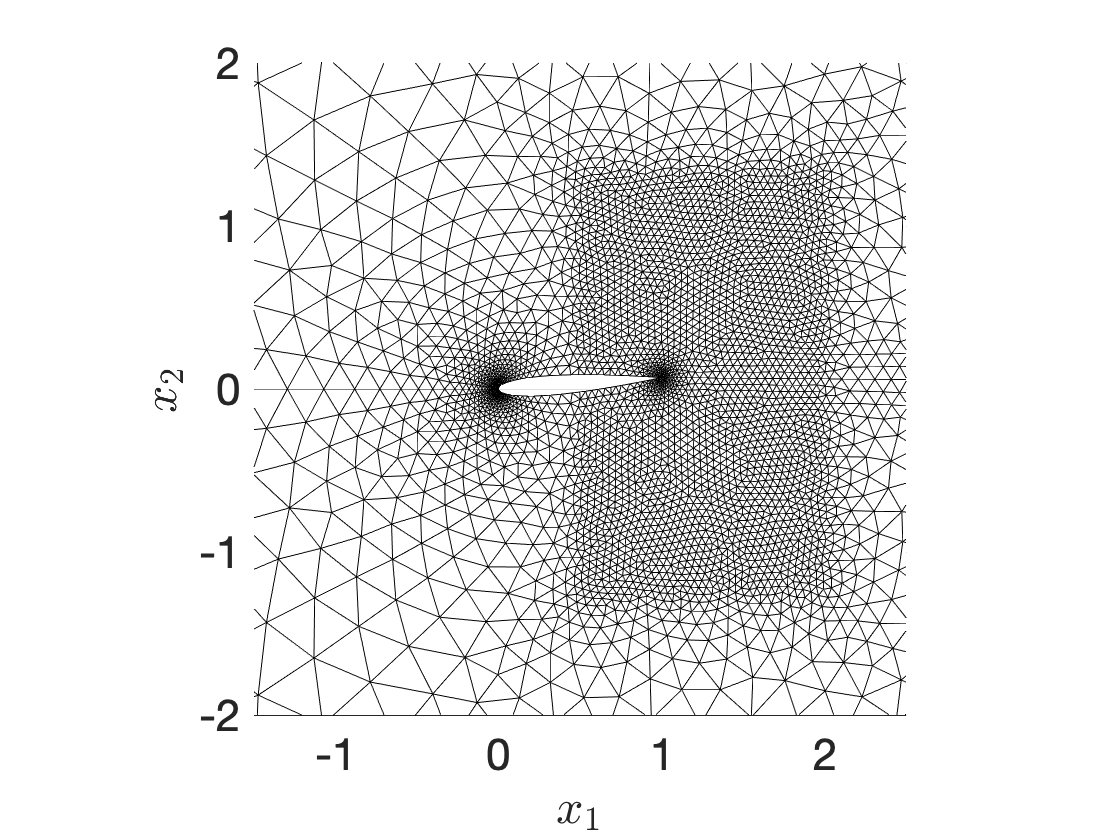}}
~~
\subfloat[ ]{
\includegraphics[width=.45\textwidth]{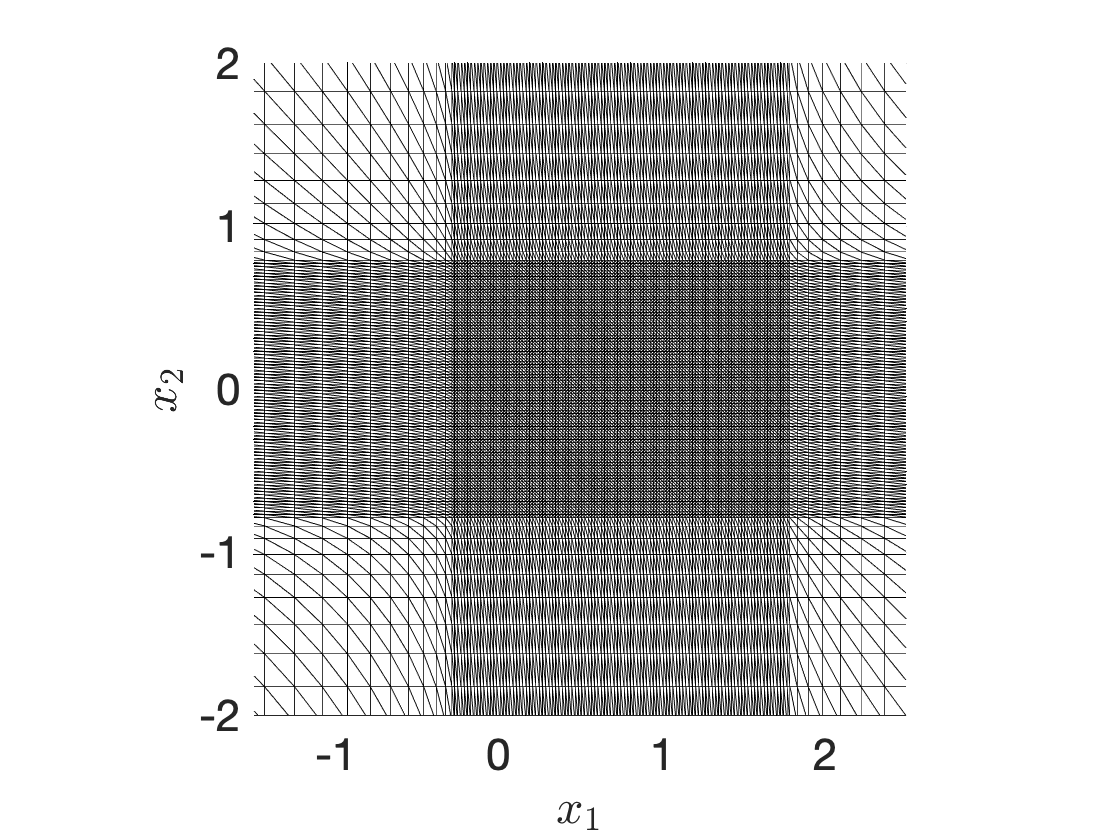}
}
        
\caption{transonic flow past an airfoil.
(a) mesh used for DG calculations; 
(b) structured mesh used for interpolation.}
\label{fig:airfoil_mesh}
\end{figure}

Figure  \ref{fig:airfoil_mesh}(a) shows the computational mesh used for DG calculations. For simplicity of implementation, we here apply our interpolation procedure in the mapped domain corresponding to angle of attack  
 $\alpha = 0^o$; furthermore, for efficiency reasons, interpolation is performed on the structured mesh in Figure \ref{fig:airfoil_mesh}(b). Since the proposed approach does not ensure bijectivity in the domain $\Omega$, it is necessary to extend the solution field inside the airfoil: we here build the extension based on the solution to a Laplace problem in the interior of the airfoil.
 
\begin{figure}[H]
\centering
\subfloat[$M=0.77$]{
\includegraphics[width=.48\textwidth]{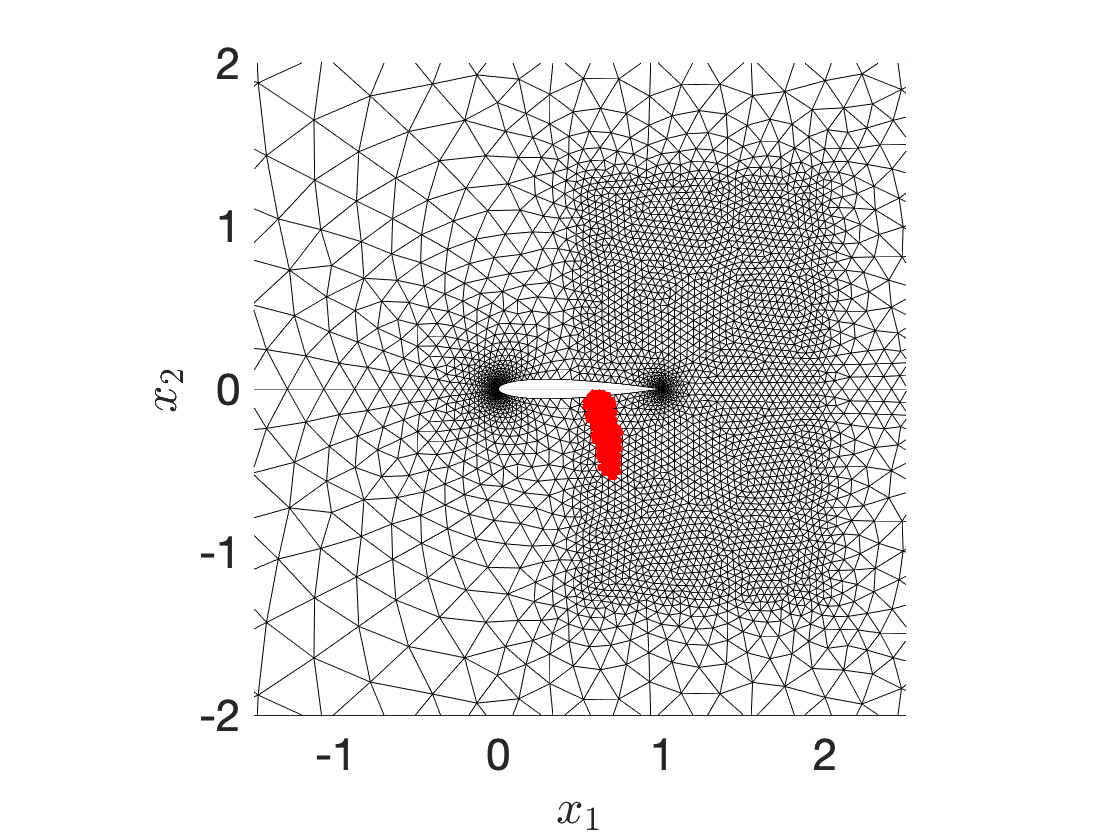}}
~~
\subfloat[ $M=0.83$]{
\includegraphics[width=.48\textwidth]{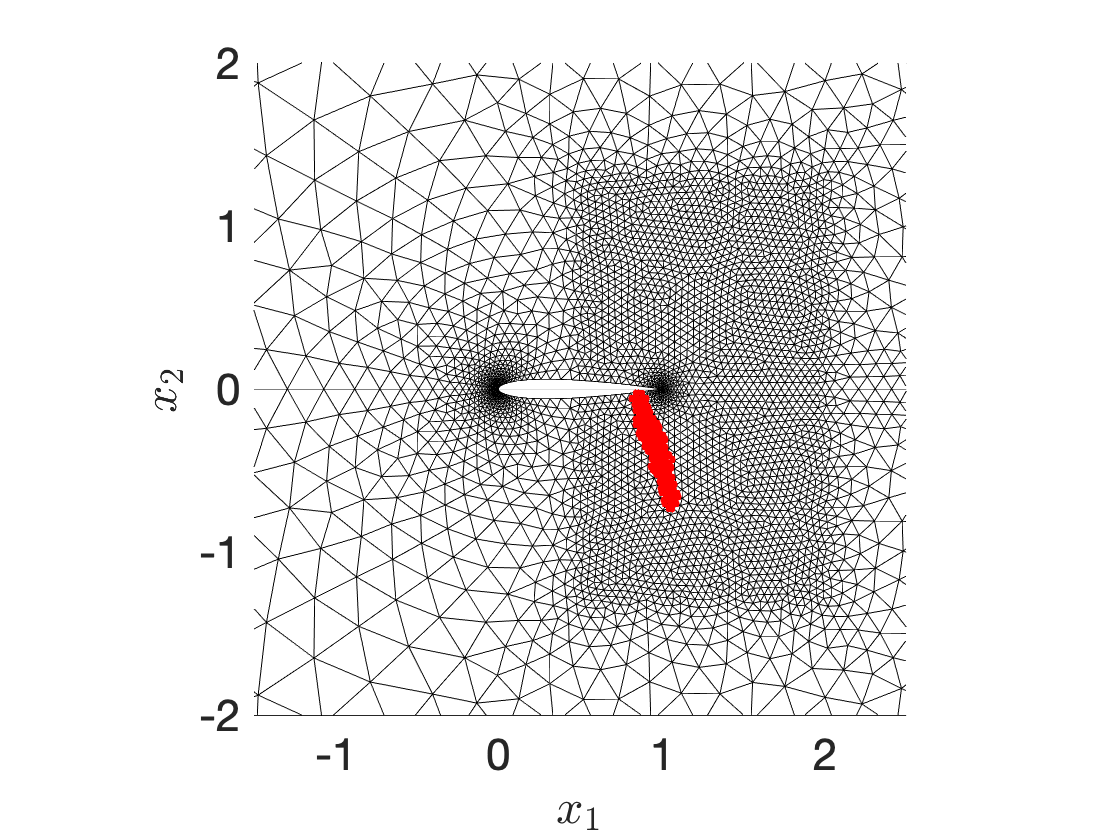}
}
        
\caption{transonic flow past an airfoil.
Elements of $P_{\rm hf}^+$ for two values of the Mach number.}
\label{fig:airfoil_selectpoints}
\end{figure}

The definition of an effective scalar testing function that detects the presence of shock discontinuities is challenging due to  numerical dissipation. We here proceed as follows: first, we define $P_{\rm hf}$ as the set of elements' centers and we compute the indicator
\begin{equation}
\label{eq:modified_ducros}
d_k \, = \, 
\max_{x \in \texttt{D}_k} \, \big| \phi(x; U)    \big|,
\quad
\phi(x; U)  : =
\frac{ \left( -\nabla \cdot u  \right)^+}{
\sqrt{  ( \nabla \cdot u )^2 + \| \nabla \times u  \|_2^2 + a^2    }}
\, 
\frac{\| \nabla p  \|_2}{p +\epsilon} \| u \|_2
\end{equation}
where $\texttt{D}_k$ denotes the $k$-th element of the DG mesh,  $k=1,\ldots,N_{\rm e}$, and $\epsilon=10^{-4}$; then, we define 
$P_{\rm hf}^+$ as the set associated with the largest $0.5\%$  values of the indicator 
$\{ d_k  \}_{k=1}^{N_{\rm e}}$. 
We observe that the first term in $\phi(x; U) $ is a modified Ducros sensor (see \cite{modesti2017low,nicoud1999subgrid}) that identifies strong compressions of the flow, the second term identifies regions characterized by large pressure gradients and the third term is intended to discard regions where the  velocity is small --- such as the leading edge.
We further remark that the indicator \eqref{eq:modified_ducros} is used in 
\cite{ferrero2020hybrid} to define the artificial viscosity for high-order DG  discretizations of inviscid flows.
Figure \ref{fig:airfoil_selectpoints} shows the selected points for two values of the Mach number.

\begin{figure}[H]
\centering
\subfloat[ ]{
\includegraphics[width=.45\textwidth]{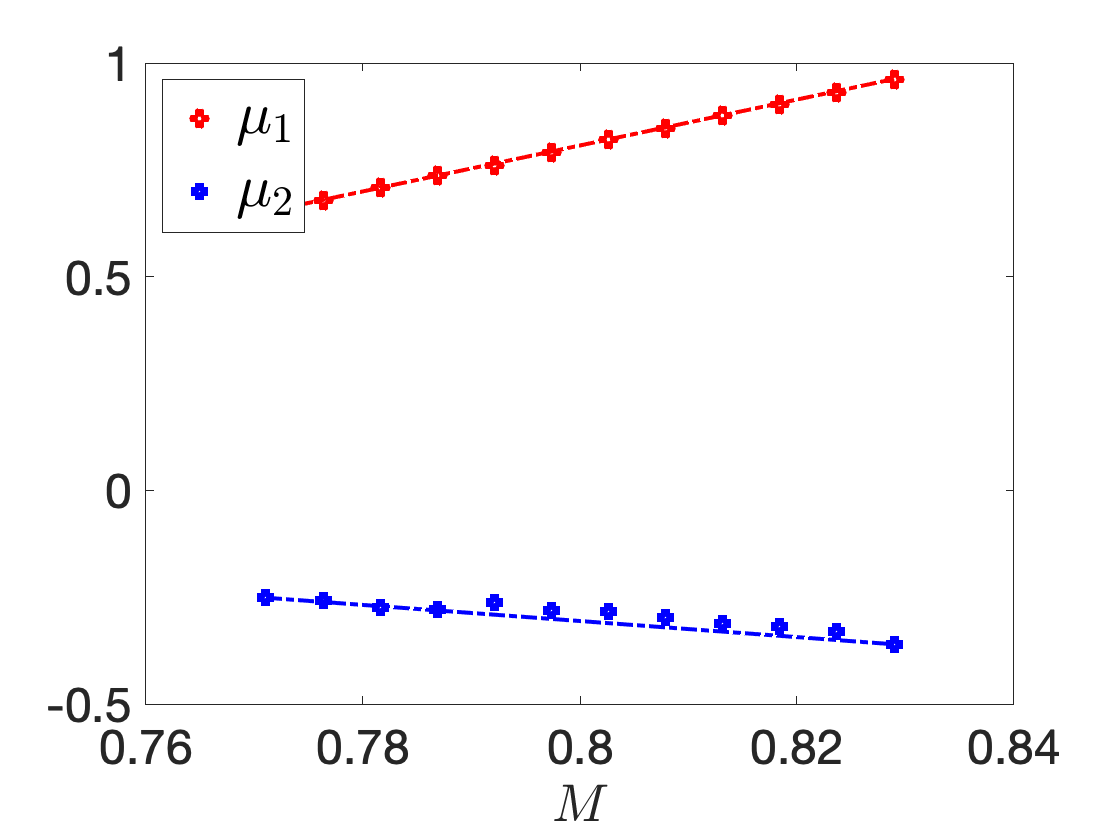}}
~~
\subfloat[ ]{
\includegraphics[width=.45\textwidth]{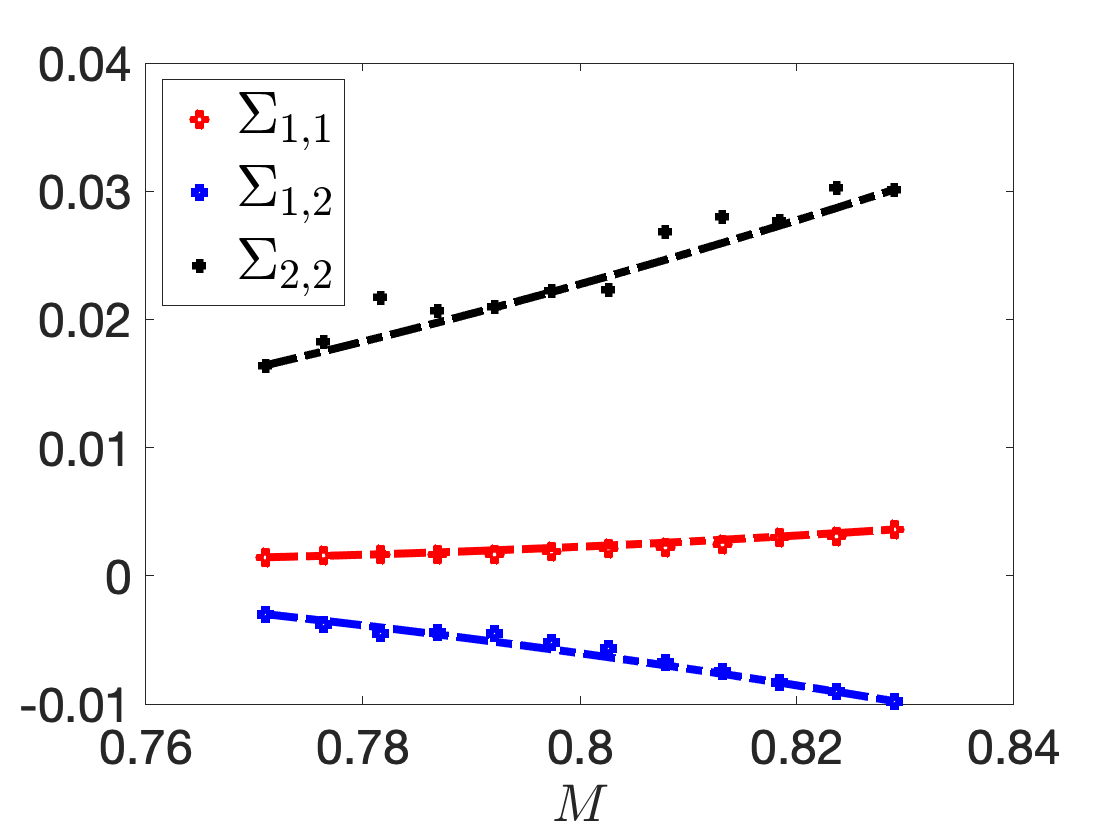}
}
        
\caption{transonic flow past an airfoil.
Behavior of the MLE estimates $\mu_{\rm mle}$ and $\Sigma_{\rm mle}$  \eqref{eq:gauss_model} with respect to the Mach number.
}
\label{fig:airfoil_mle}
\end{figure}

\begin{figure}[H]
\centering
\subfloat[]{
\includegraphics[width=.45\textwidth]{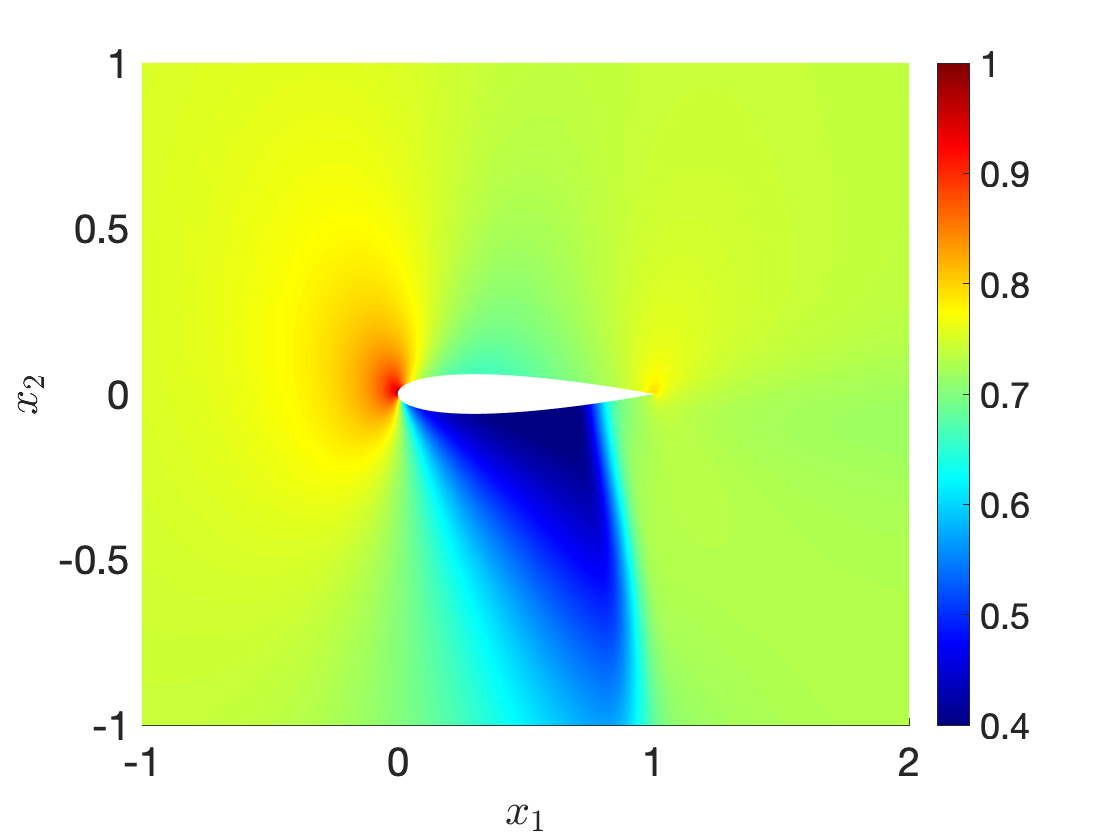}}
~~
\subfloat[ ]{
\includegraphics[width=.45\textwidth]{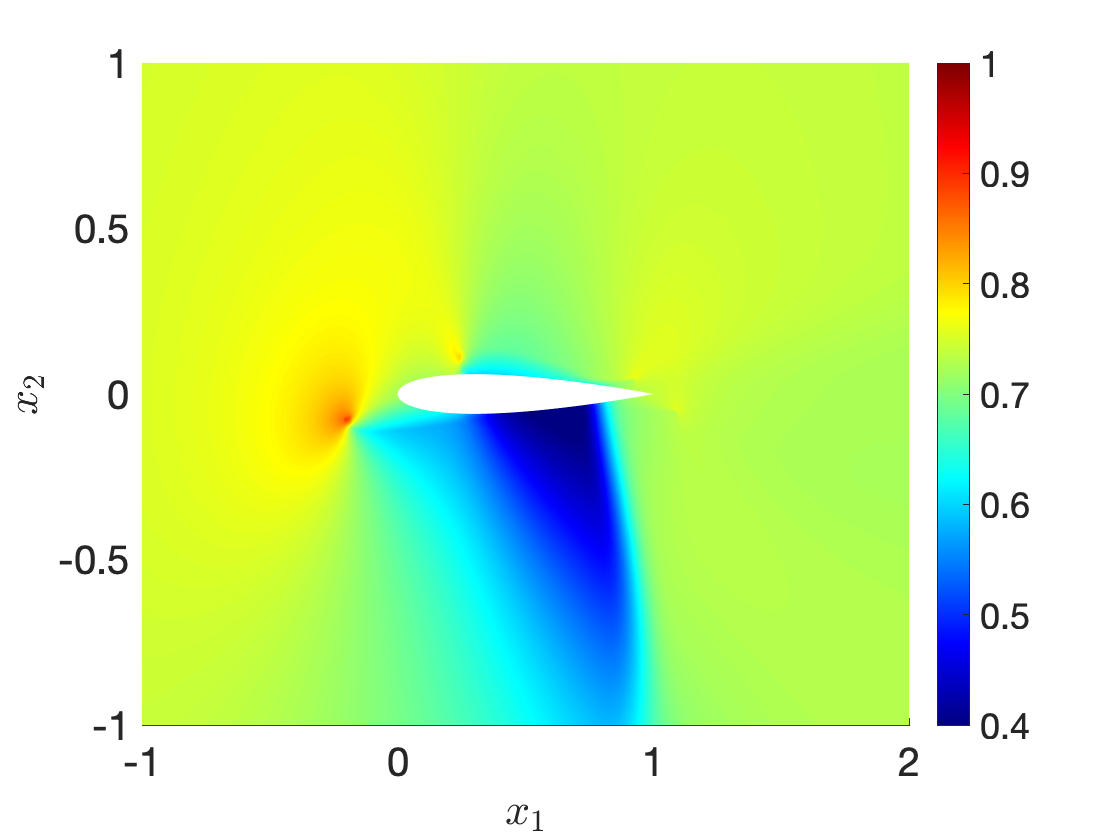}
}
        
\subfloat[ $x_2=-0.1$]{
\includegraphics[width=.45\textwidth]{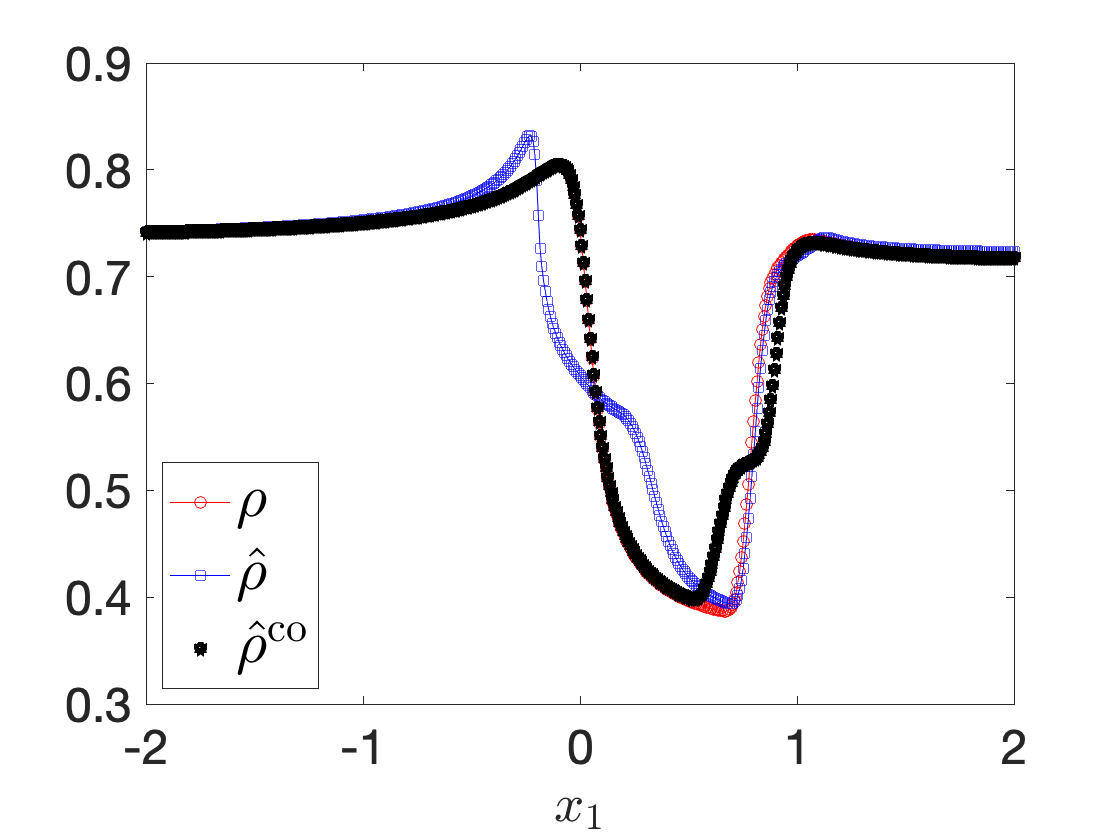}}
~~
\subfloat[$x_2=-0.5$ ]{
\includegraphics[width=.45\textwidth]{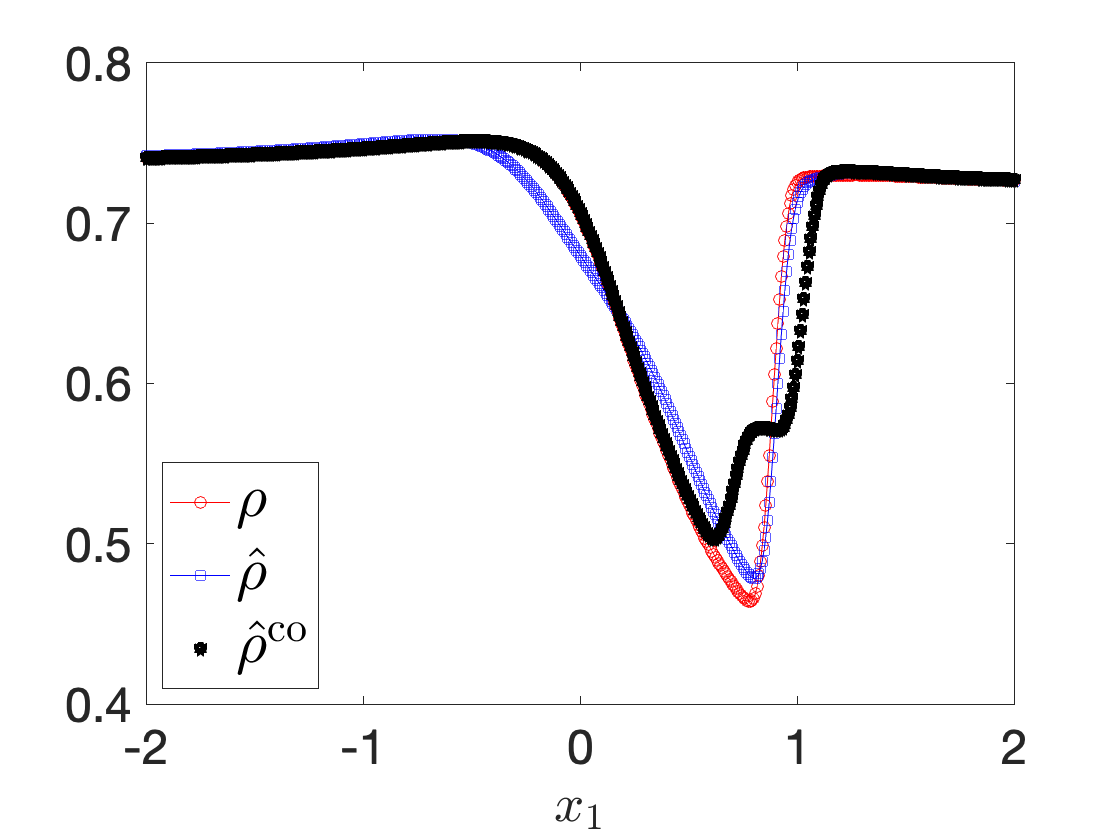}
}              
        
\caption{transonic flow past an airfoil.
(a)-(b)  behavior  of the density field
for $M=0.8$
 and of the CDI 
 \eqref{eq:convex_displacement_interpolation}
 for $s=1/2$.
 (c)-(d) two horizontal slices of truth and predicted density profiles.
}
\label{fig:airfoil_vis}
\end{figure}

In Figure \ref{fig:airfoil_mle}, we investigate the behavior of the MLE estimates $\mu_{\rm mle}$ and $\Sigma_{\rm mle}$ in \eqref{eq:gauss_model} with respect to the Mach number. We observe that both mean and variance are smooth functions of the parameter. In Figure \ref{fig:airfoil_vis}, we compare the density field for $M=0.8$ ($\alpha=0.5$) with the 
CDI  \eqref{eq:convex_displacement_interpolation}
with $s=\alpha$; we further provide horizontal slices of the DG solution (in red), the CDI (in blue) and   the convex interpolant (in black) for two values of $x_2$.
We observe that the CDI  is extremely accurate in the proximity of the shock, while it is highly inaccurate far from the shock,  especially in the proximity of the airfoil.

\section{Extension: boundary-aware transportation of Gaussian models}
\label{sec:extension}

The examples of the previous section show that CDI  \eqref{eq:convex_displacement_interpolation}  based on  optimal transportation of Gaussian models is effective if 
(i) boundaries are not present 
 (cf. section \ref{sec:numerics_simple_wave})
or the extension of the solution outside the domain is trivial (cf. section \ref{sec:wedge}), and
(ii) the solution field presents a single coherent  structure that   is well-approximated by an ellipsoid.
Inaccuracy of the displacement interpolation for the example in section \ref{sec:transonic_flow} is the consequence of two factors.
First, since displacement interpolation does not preserve the boundaries of the domain,  interpolation might be  highly inaccurate in the neighborhood of the airfoil
(cf. Figure \ref{fig:airfoil_vis}(c)) {and highly depends on the choice of the extension operator, which is  typically very difficult to construct, particularly for slender bodies}. Second,
the  Gaussian model considered   is not able to take into account the coherent structures that develop at leading and trailing edges: it is thus a too simplistic representation of the solution field.

Based on these considerations, we propose here an  extension of the approach in section \ref{sec:method}: CDI  based on boundary-aware (BA) transportation of multiple Gaussian models. We investigate performance of our  approach for the transonic flow test case introduced in section \ref{sec:transonic_flow}.
For simplicity, in the remainder we assume that the domain $\Omega$ is a Lipschitz two-dimensional domain.

\subsection{Methodology}
\label{sec:method_extension}

Given $M\in \mathbb{N}$, we introduce the approximation map 
$\mathcal{N}: \Omega \times \mathbb{R}^M \to \mathbb{R}^2$ and the set $\mathcal{A}_{\rm bj}\subset \mathbb{R}^M$ such  that
$\mathcal{N}(\cdot,\mathbf{a})$ is a bijection in $\Omega$ for all $\mathbf{a} \in \mathcal{A}_{\rm bj}$. We pursue the approach in 
\cite{taddei2021registration}
to define  
$\mathcal{N}$: 
we refer to   \ref{sec:registration} for further  details.
We denote by  
$P_{\rm hf}^{+,0} = \{ y_j^0  \}_{j=1}^{N_{\rm hf,0}^+}$
 and $P_{\rm hf}^{+,1} = 
 \{  y_j^1  \}_{j=1}^{N_{\rm hf,0}^+}$  the selected points for $U_0$ and $U_1$, and we denote by $X_g$ and $Y_g$ the optimal maps obtained using   
\eqref{forwaedm_gaussian}.

Then, 
we define $\widehat{\mathbf{a}}_{0,1}  \in \mathbb{R}^M$ to minimize 
\begin{subequations}
\label{eq:registration_task}
\begin{equation}
\sum_{j=1}^{N_{\rm hf,0}^+}
\| X_g(  y_j^0  ) -  \mathcal{N}( y_j^0,\mathbf{a})   \|_2^2
\, + \,
\mathfrak{P}(\mathbf{a}),
\quad {\rm subject \; to} \;\;
\mathfrak{C}(\mathbf{a}) \leq 0,
\end{equation}
where $\mathfrak{P}$ is a suitable regularization that penalizes the $H^2$ seminorm of the mapping and
$\mathfrak{C}$ is a bijectivity constraint that, combined with 
$\mathfrak{P}$,  enforces that $\widehat{\mathbf{a}}_{0,1}$ belongs to $\mathcal{A}_{\rm bj}$: we refer to \cite{taddei2021registration} for the details. Similarly, 
we define 
$\widehat{\mathbf{a}}_{1,0} \in \mathbb{R}^M$ to minimize
\begin{equation}
\sum_{j=1}^{N_{\rm hf,1}^+}
\| Y_g(  y_j^1  ) -  \mathcal{N}( y_j^1, \mathbf{a})   \|_2^2
\, + \,
\mathfrak{P}(\mathbf{a}),
\quad {\rm subject \; to} \;\;
\mathfrak{C}(\mathbf{a}) \leq 0.
\end{equation}
\end{subequations}
In conclusion, we introduce the boundary-aware (BA) CDI  as
\begin{equation}
\label{eq:convex_displacement_interpolation_bnd_aware}
\widehat{U}(s, x)
\, =  \,
(1 - s) U_0 \circ 
  \widetilde{W}_g(s, x)  
 \, + \,
s U_1 \circ \widetilde{T}_g (1-s,x),
\quad
s\in [0,1],  \;  x\in \Omega,
\end{equation}
where 
$ \widetilde{T}_g(s, x)= \mathcal{N}(x, \, s \cdot  \widehat{\mathbf{a}}_{0,1}) $ and
$ \widetilde{W}_g(s, x)= \mathcal{N}(x, \, s \cdot  \widehat{\mathbf{a}}_{1,0}) $.
In the next section, we investigate performance of \eqref{eq:convex_displacement_interpolation_bnd_aware} for the transonic flow test case.

Registration provides   sub-optimal  --- in the sense of optimal transportation --- 
bijective-in-$\Omega$
approximations of the Gaussian maps obtained using \eqref{forwaedm_gaussian}: 
{we can thus  interpret $\widetilde{T}_g$ and $\widetilde{W}_g$ as approximate projections of the actions on marked points of the  optimal transport maps 
${T}_g$ and ${W}_g$ onto the space of bijective maps in $\Omega$.}

We cannot in general guarantee that $\widetilde{T}_g$ and $\widetilde{W}_g$ are bijections in  $\Omega$ for all $s\in [0,1]$: nevertheless, in our experience, provided that the distance --- in the sense of Wasserstein --- between the Gaussian models associated with  $P_{\rm hf}^{+,0}$ and $P_{\rm hf}^{+,1}$ is moderate, the approach leads to bijective maps for all $s\in [0,1]$.
{For large deformations, we might set $\widetilde{T}_g(s,\cdot)= \mathcal{N}(\cdot; \widehat{\mathbf{a}}_{0,1}(s))$ where $\widehat{\mathbf{a}}_{0,1}(s)$ minimizes
$$
\sum_{j=1}^{N_{\rm hf,1}^+}
\| T_g(s,   y_j^0  ) -  \mathcal{N}( y_j^0, \mathbf{a})   \|_2^2
\, + \,
\mathfrak{P}(\mathbf{a}),
\quad {\rm subject \; to} \;\;
\mathfrak{C}(\mathbf{a}) \leq 0.
$$
Similarly, we set $\widetilde{W}_g(s, \cdot) = \mathcal{N}(\cdot; \widehat{\mathbf{a}}_{1,0}(s))$ where 
$\widehat{\mathbf{a}}_{1,0}(s)$ is defined based on 
$\{ W_g(s,   y_j^1  ) \}_j$.
This choice of the mapping ensures bijectivity of 
$\widetilde{T}_g,\widetilde{W}_g$ for all $s\in [0,1]$ at the price of additional offline costs. 
}

Following \cite{chen2018optimal}, it is straightforward to extend \eqref{eq:registration_task}-\eqref{eq:convex_displacement_interpolation_bnd_aware} to track multiple structures.
First, given the points
$$
 \left\{   y_j^{0,k} \,: \, j=1,\ldots, 
N_{\rm hf,0,k}^+, k=1,\ldots,N_{\rm g} \right\},
\quad
 \left\{   y_j^{1,k} \,: \, j=1,\ldots, 
N_{\rm hf,1,k}^+, k=1,\ldots,N_{\rm g} \right\}
$$
we first compute the mappings $\{T_{\rm g}^{(k,k')} \}_{k,k'}$ and the associated Wasserstein distances 
$\{W_2^{(k,k')} \}_{k,k'}$ using the identities in section \ref{sec:gauss_model_coherent}. Then, we identify the permutation $I$ of $\{1,\ldots,N_{\rm g}\}$ that minimizes 
$$
\sum_{k=1}^{N_{\rm g}} \, W_2^{k,I_k}
$$
over all possible permutations. Finally, we compute 
$\widehat{\mathbf{a}}_{0,1}  \in \mathbb{R}^M$ to minimize 
\begin{subequations}
\label{eq:registration_task_plus}
\begin{equation}
\sum_{k=1}^{N_{\rm g}}
\sum_{j=1}^{N_{\rm hf,0,k}^+}
\| T_g^{k,I_k}(  y_j^{0,k}  ) -  \mathcal{N}( y_j^{0,k},\mathbf{a})   \|_2^2
\, + \,
\mathfrak{P}(\mathbf{a}),
\quad {\rm subject \; to} \;\;
\mathfrak{C}(\mathbf{a}) \leq 0,
\end{equation}
and we compute  $\widehat{\mathbf{a}}_{1,0}  \in \mathbb{R}^M$ to minimize 
\begin{equation}
\sum_{k=1}^{N_{\rm g}}
\sum_{j=1}^{N_{\rm hf,1,k}^+}
\| R_g^{k,I_k}(  y_j^{1,k}  ) -  \mathcal{N}( y_j^{1,k},\mathbf{a})   \|_2^2
\, + \,
\mathfrak{P}(\mathbf{a}),
\quad {\rm subject \; to} \;\;
\mathfrak{C}(\mathbf{a}) \leq 0,
\end{equation}
which are of the same form as  \eqref{eq:registration_task}.  
We validate this approach  through  the vehicle of a transonic flow with two shocks.
\end{subequations}

\subsection{Numerical results}

\subsubsection{Transonic flow at angle of attack $4^o$}

\begin{figure}[H]
\centering
\subfloat[$s=0.25$]{
\includegraphics[width=.45\textwidth]{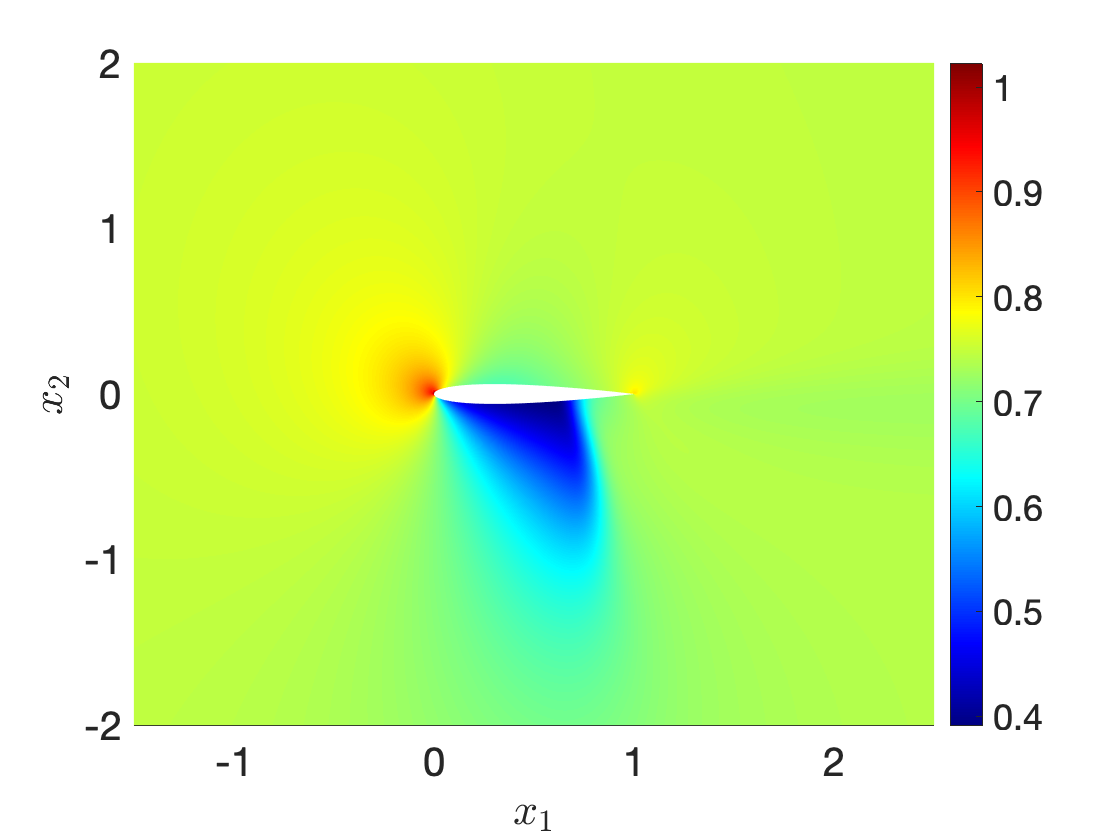}}
~~
\subfloat[$s=0.75$]{
\includegraphics[width=.45\textwidth]{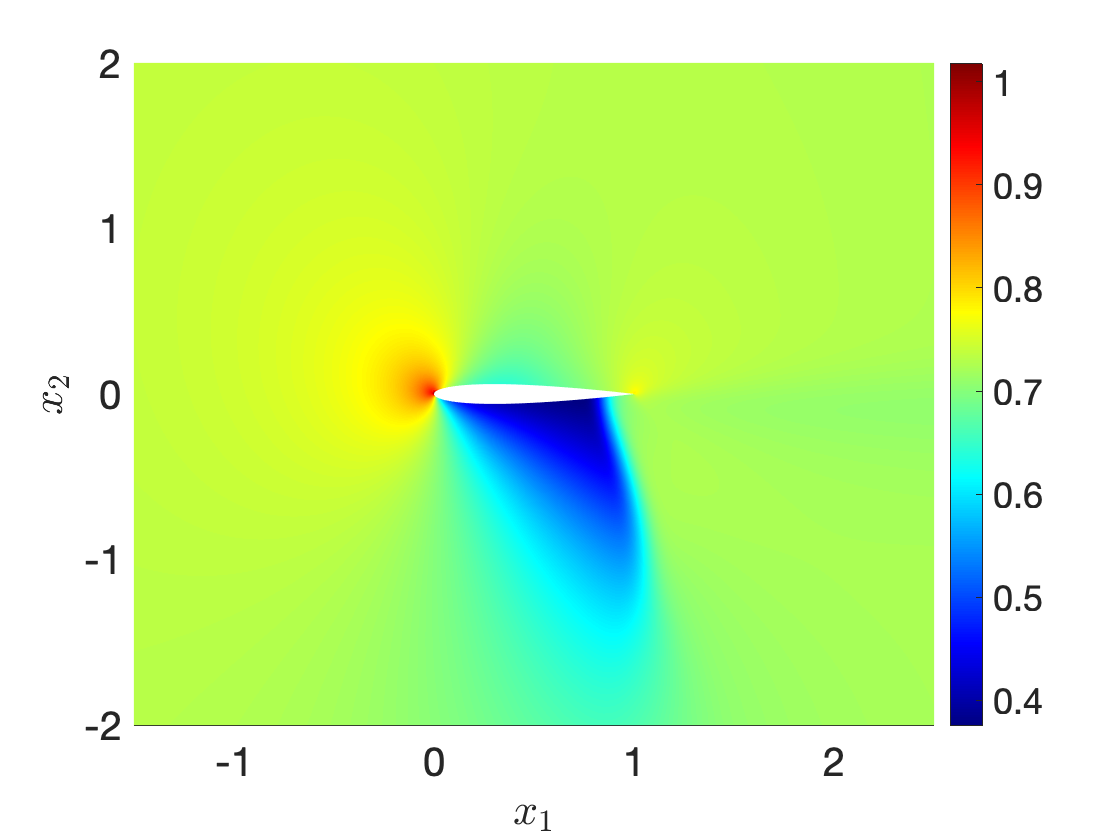}}
%
%
%
        
\caption{boundary-aware transportation of Gaussian models (angle of attack $4^o$). BA CDI   
\eqref{eq:convex_displacement_interpolation_bnd_aware}
for two values of $s$. }
\label{fig:airfoil_density_bnd_aware}
\end{figure}

Figure \ref{fig:airfoil_density_bnd_aware} shows the behavior of the BA CDI  \eqref{eq:convex_displacement_interpolation_bnd_aware}
for two values of $s\in [0,1]$: we observe that the proposed interpolation preserves the structures at leading and trailing edges and is able to smoothly deform the shock attached to the airfoil.

\begin{figure}[H]
\centering
\subfloat[]{
\includegraphics[width=.45\textwidth]{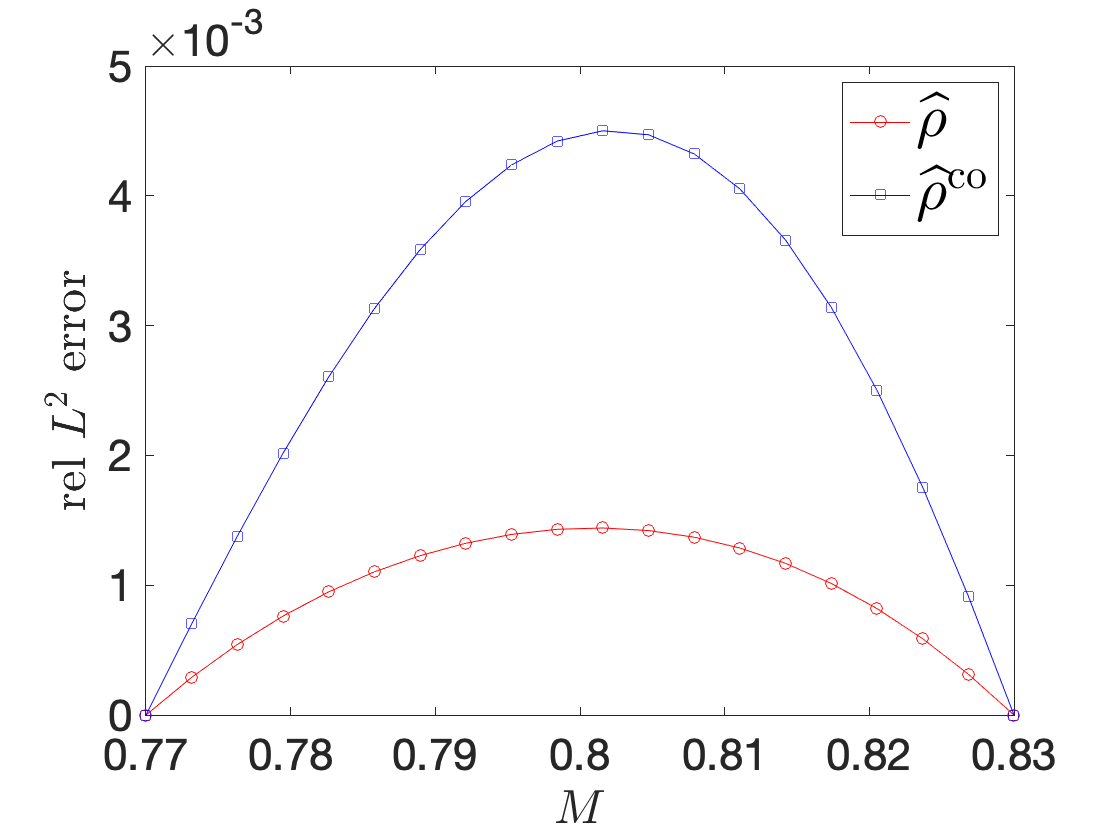}}

\subfloat[]{
\includegraphics[width=.45\textwidth]{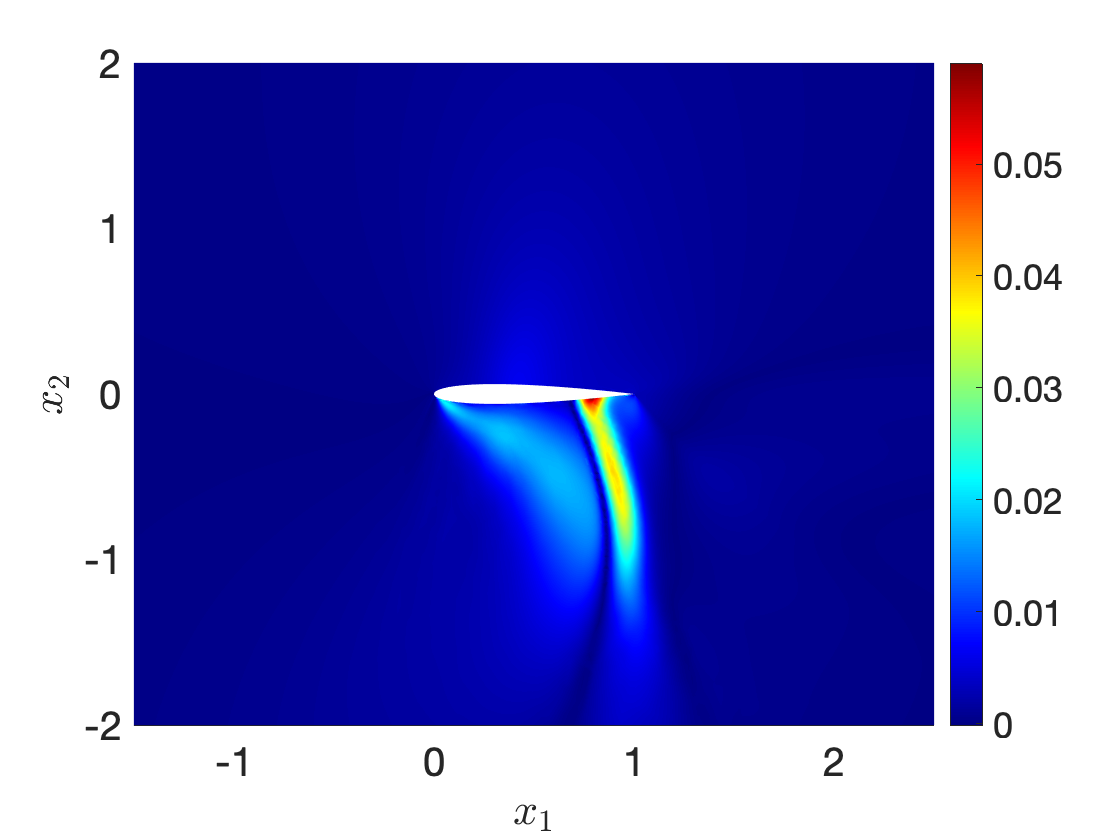}}
~~
\subfloat[]{
\includegraphics[width=.45\textwidth]{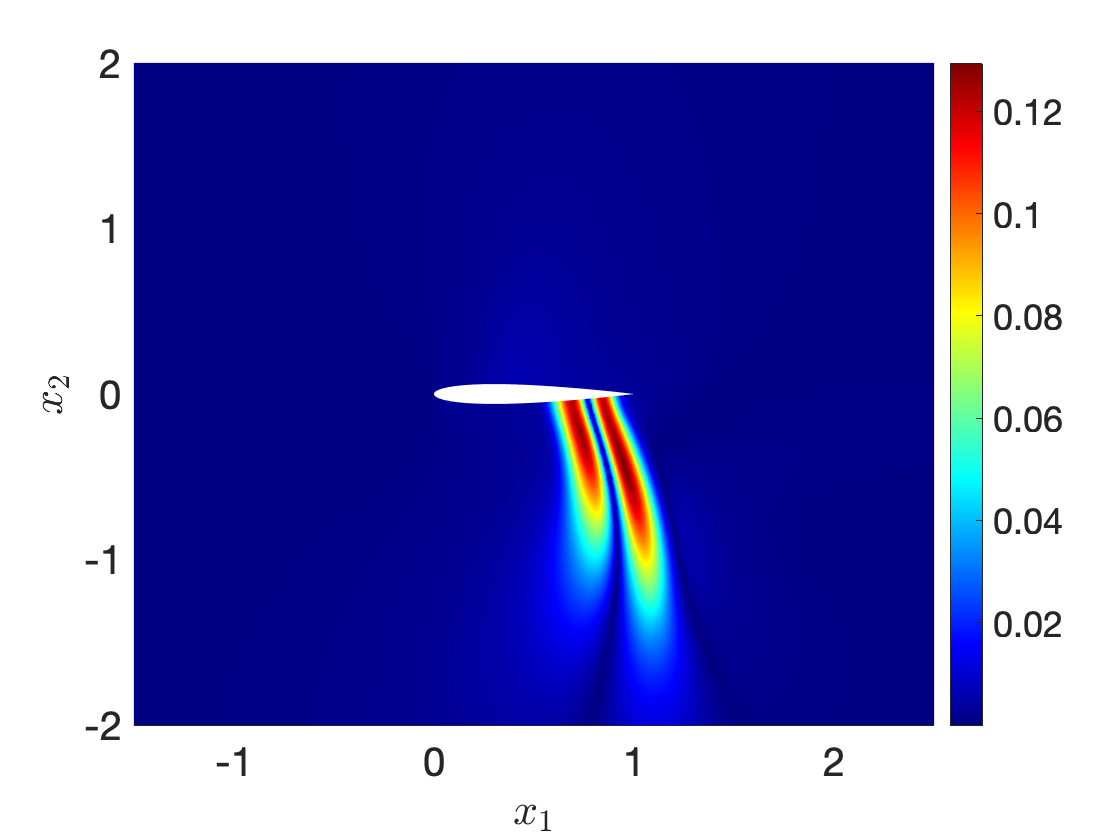}}
        
\caption{boundary-aware transportation of Gaussian models (angle of attack $4^o$). (a) behavior of the relative $L^2$ error for BA CDI.
b)-(c) behavior  of the error fields
$|   \widehat{\rho}(s,x)  -  {\rho} (x; M) |$ and
$|   \widehat{\rho}^{\rm co}(s,x)  -  {\rho} (x; M)  |$
for $M=0.8$ and $s=0.5$. 
 }
\label{fig:airfoil_density_bnd_aware_2}
\end{figure}

\begin{figure}[H]
\centering
\subfloat[$x_2=-0.1$]{
\includegraphics[width=.45\textwidth]{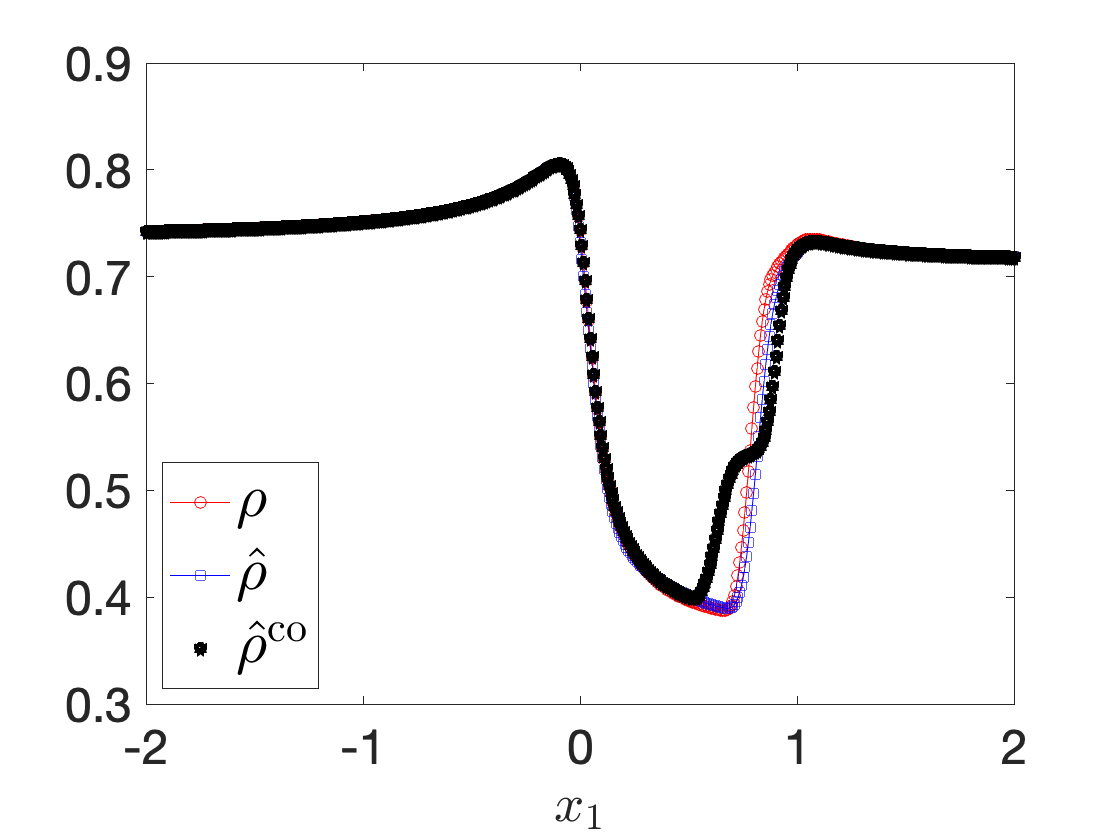}}
~~
\subfloat[$x_2= - 0.5$]{
\includegraphics[width=.45\textwidth]{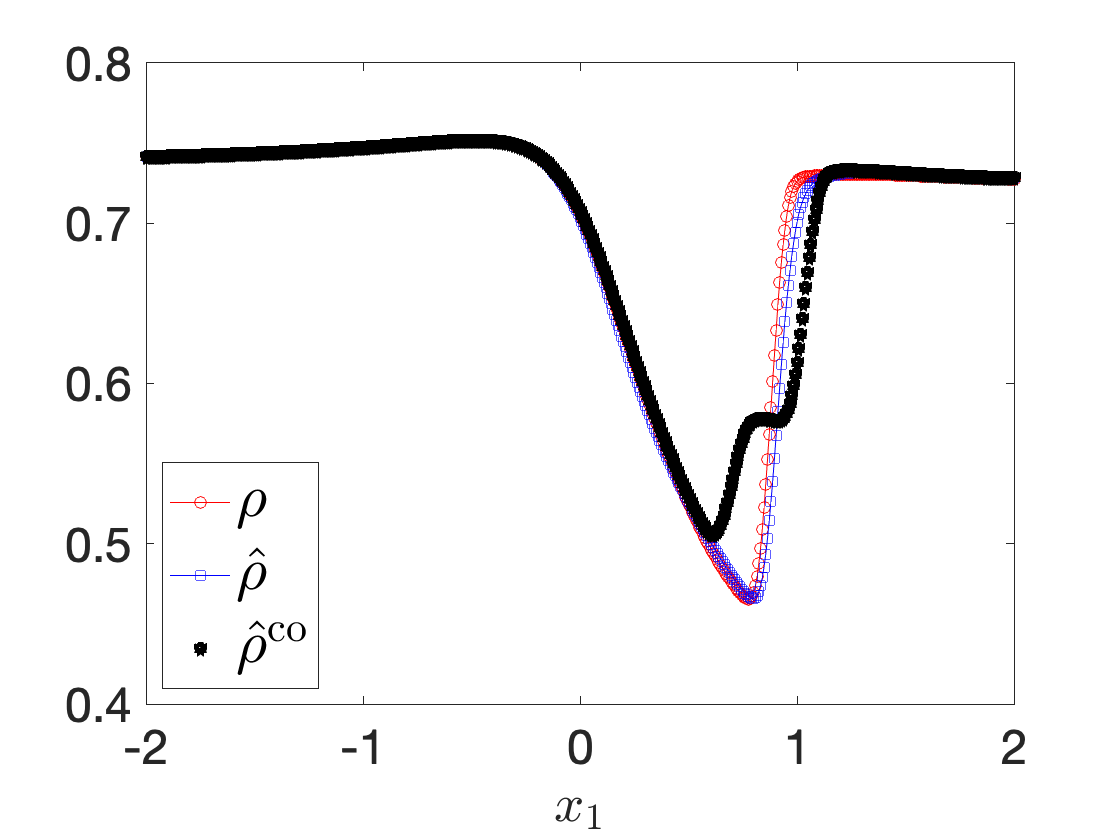}}
        
\caption{boundary-aware transportation of Gaussian models (angle of attack $4^o$). (a)-(b)  horizontal slices of  truth and predicted density profiles for $M=0.8$.
 }
\label{fig:airfoil_density_bnd_aware_3}
\end{figure}

Figure  \ref{fig:airfoil_density_bnd_aware_2}
compares performance of the nonlinear interpolation \eqref{eq:convex_displacement_interpolation_bnd_aware} with the linear convex interpolation for $M\in [0.77,0.83]$ --- we here consider 
$s_{\alpha} = \alpha$ for both linear and nonlinear interpolation.
Similarly, Figure \ref{fig:airfoil_density_bnd_aware_3} 
compares 
horizontal slices of  the truth  density profile
with horizontal slices of linear and nonlinear interpolation \eqref{eq:convex_displacement_interpolation_bnd_aware},
 for $M=0.8$.
We observe that nonlinear interpolation leads to more accurate performance in terms of relative $L^2$  error, and in particular is  more accurate in the proximity of the shock.

\subsubsection{Transonic flow at angle of attack $1^o$}

\begin{figure}[H]
\centering
\subfloat[$M=0.8$]{
\includegraphics[width=.48\textwidth]{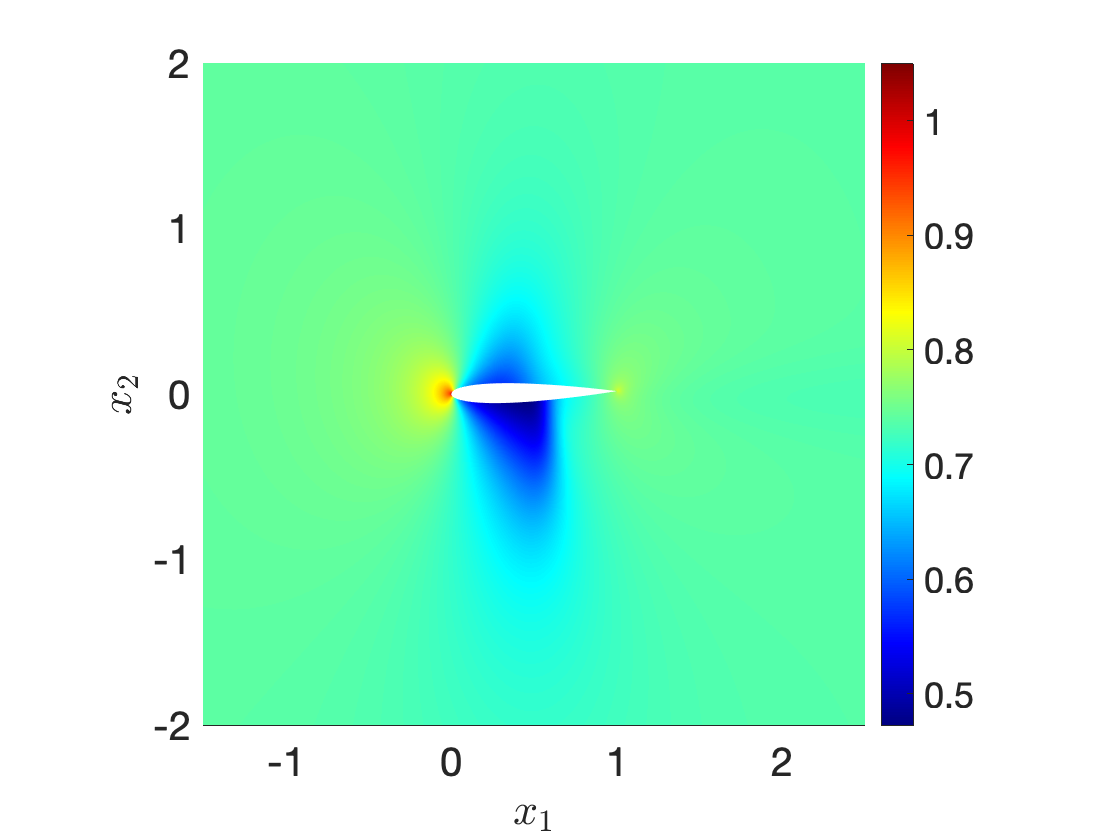}}
~~
\subfloat[$M=0.86$]{
\includegraphics[width=.48\textwidth]{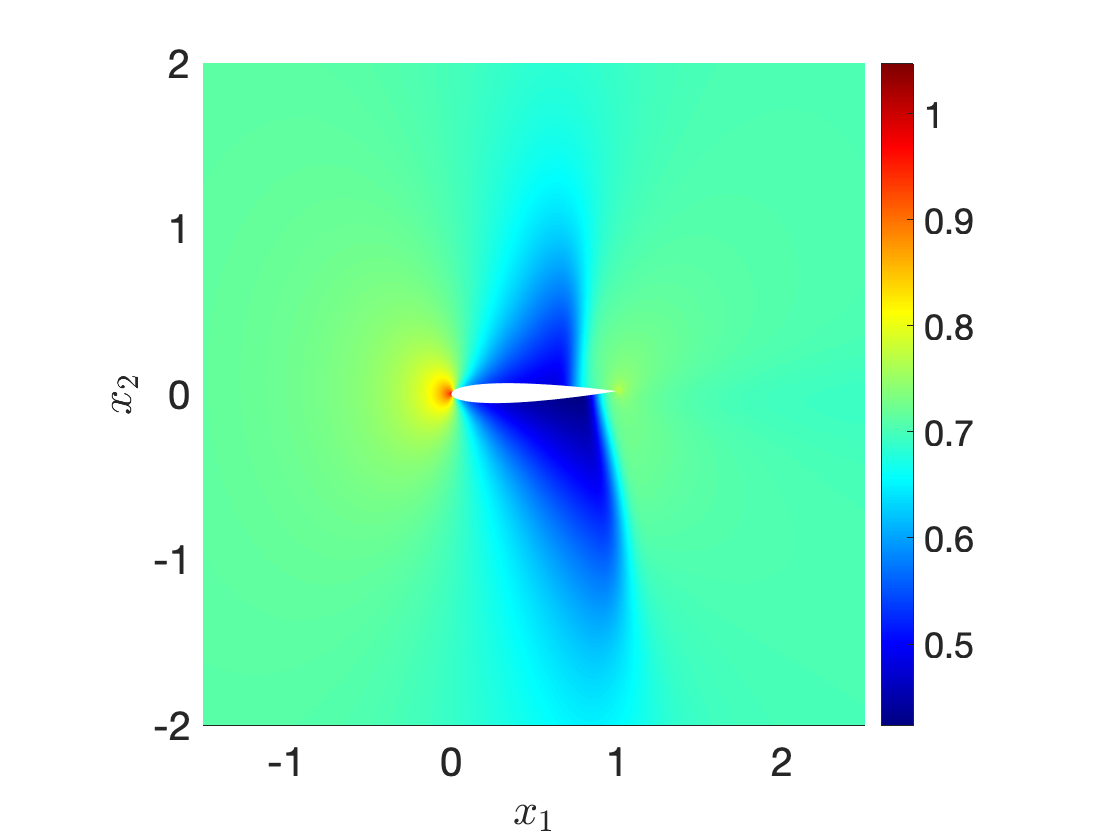}
}
        
\caption{transonic flow past an airfoil at angle of attack $1^o$. Flow density for two values of the Mach number.}
\label{fig:airfoil_density_1}
\end{figure}

We consider a two-dimensional transonic flow past a NACA 0012 airfoil at angle of attack $\alpha = - 1^o$; we let the solution vary with  respect  to the free-stream Mach number $M\in [0.8,0.86]$. As shown in 
Figure \ref{fig:airfoil_density_1}, the flow density exhibits two shocks that are very sensitive to the value of   the parameter.

\begin{figure}[H]
\centering
\subfloat[$M=0.8$]{
\includegraphics[width=.45\textwidth]{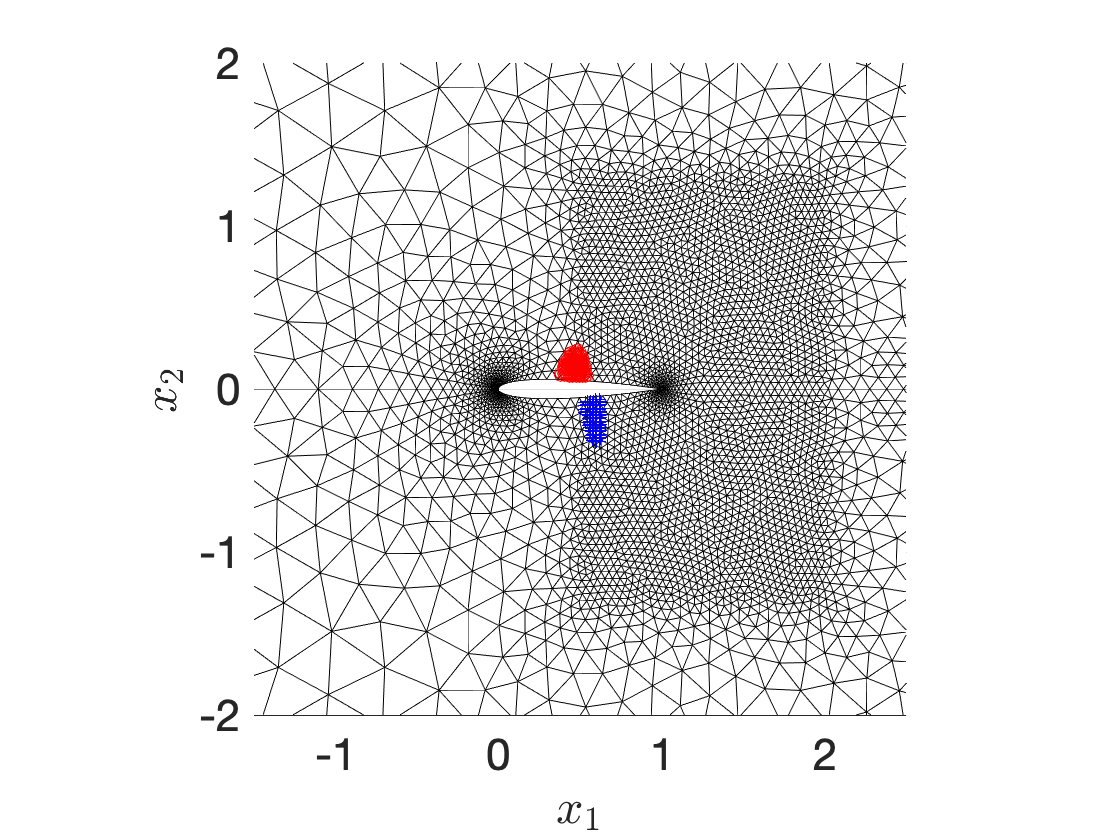}}
~~
\subfloat[ $M=0.86$]{
\includegraphics[width=.45\textwidth]{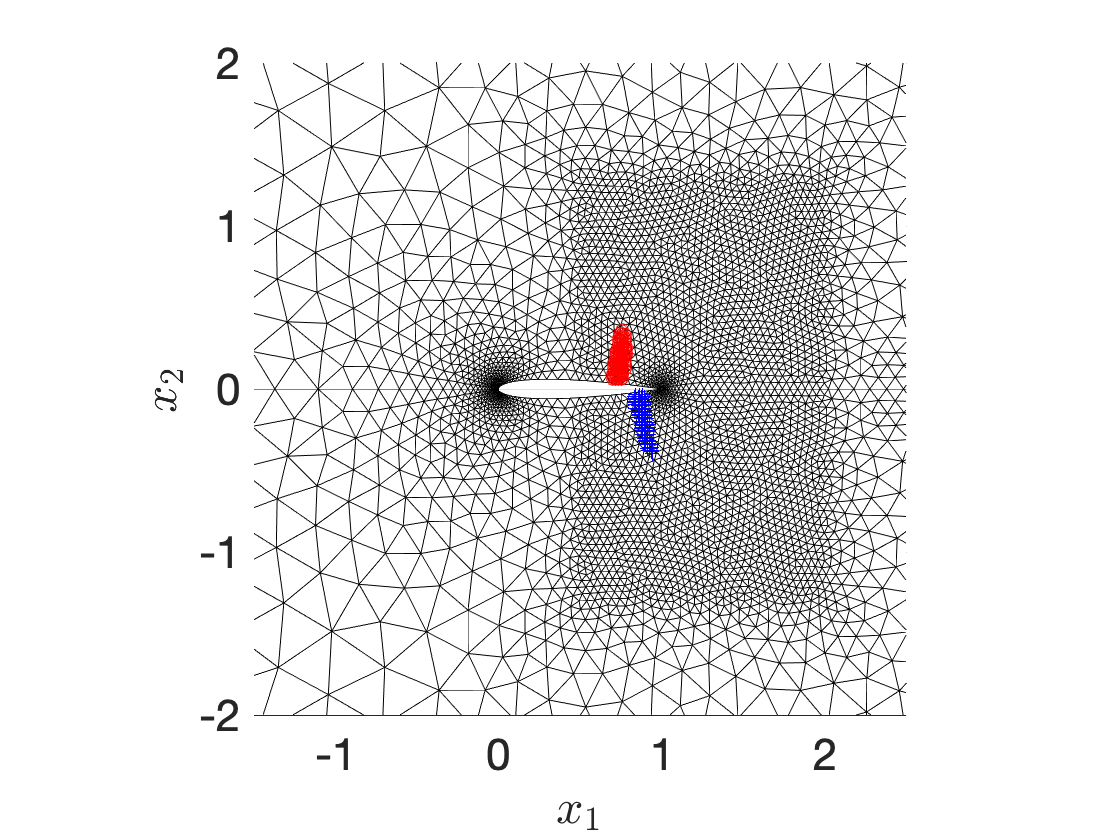}
}
        
\caption{transonic flow past an airfoil at angle of attack $1^o$.
Elements of $P_{\rm hf}^+$ for two values of the Mach number. }
\label{fig:airfoil_selectpoints_1}
\end{figure}

Figure \ref{fig:airfoil_selectpoints_1}  shows the selected points $P_{\rm hf}^+$ for two values of the Mach number.
We resort to the same indicator introduced in \eqref{eq:modified_ducros} to identify the set $P_{\rm hf}^+$. To facilitate the interpolation task we discard points outside $(0,1)\times \mathbb{R}$; furthermore, we separate \emph{a priori} the two clouds of points by discriminating between positive and negative heights. Note that the latter expedient  allows us to apply the  procedure  described in section \ref{sec:gauss_model_coherent} to build the Gaussian models for the upper and lower shocks, and ultimately robustifies the identification task. In the future, we wish to investigate performance of automated detection algorithms for Gaussian mixtures, \cite{mclachlan2019finite}.

\begin{figure}[H]
\centering
\subfloat[$s=0.25$]{
\includegraphics[width=.45\textwidth]{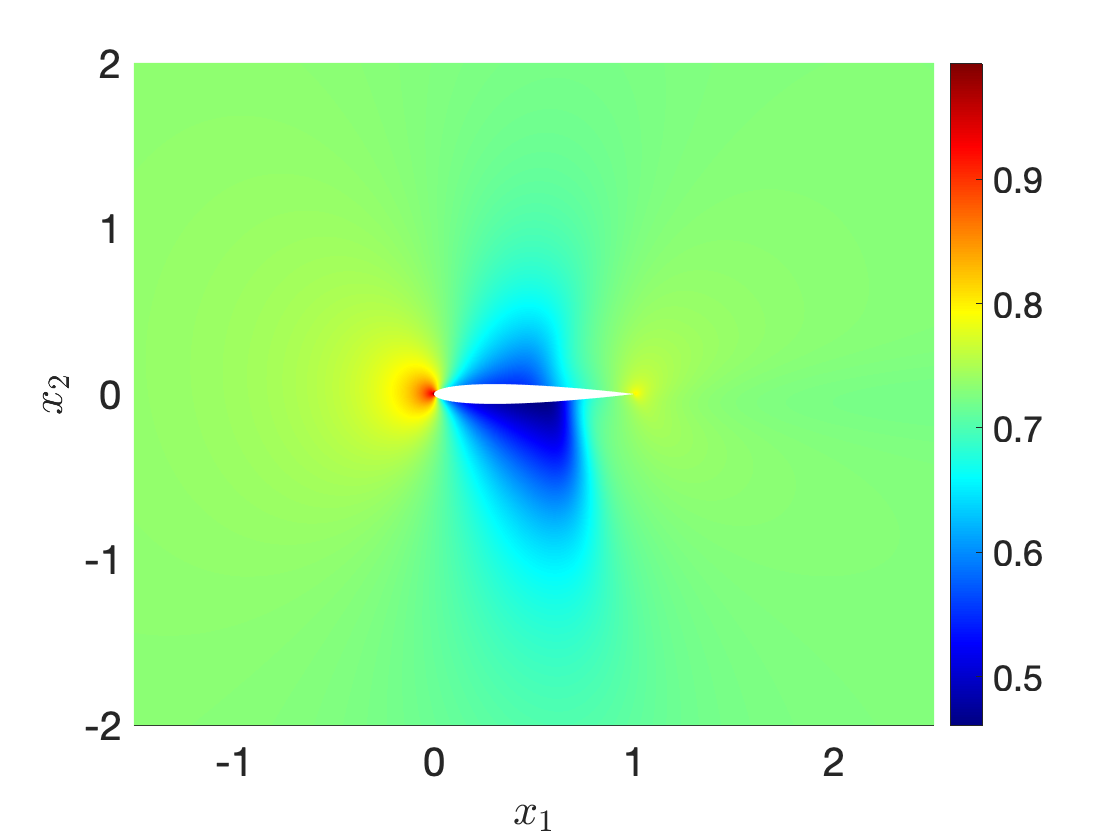}}
~~
\subfloat[$s=0.75$]{
\includegraphics[width=.45\textwidth]{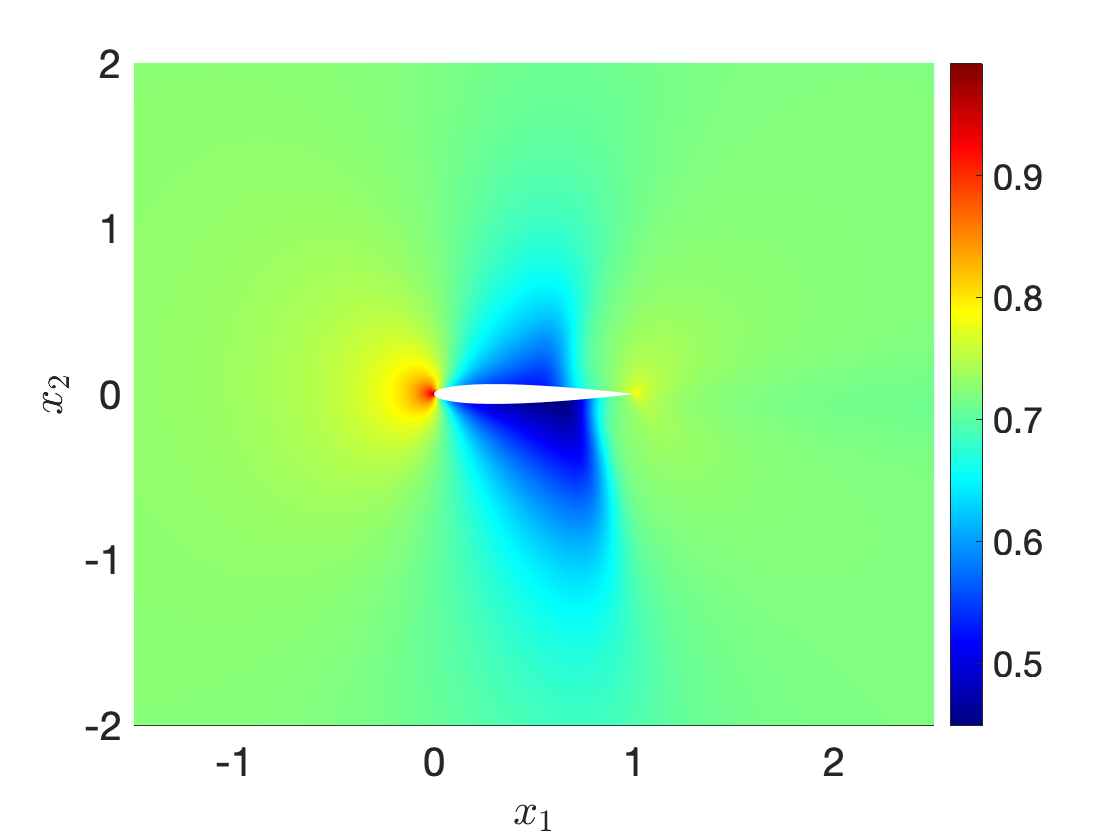}}
        
\caption{boundary-aware transportation of Gaussian models (angle of attack $1^o$). BA CDI   
\eqref{eq:convex_displacement_interpolation_bnd_aware}
for two values of $s$. }
\label{fig:airfoil_density_bnd_aware_angle1}
\end{figure}

\begin{figure}[H]
\centering
\subfloat[]{
\includegraphics[width=.45\textwidth]{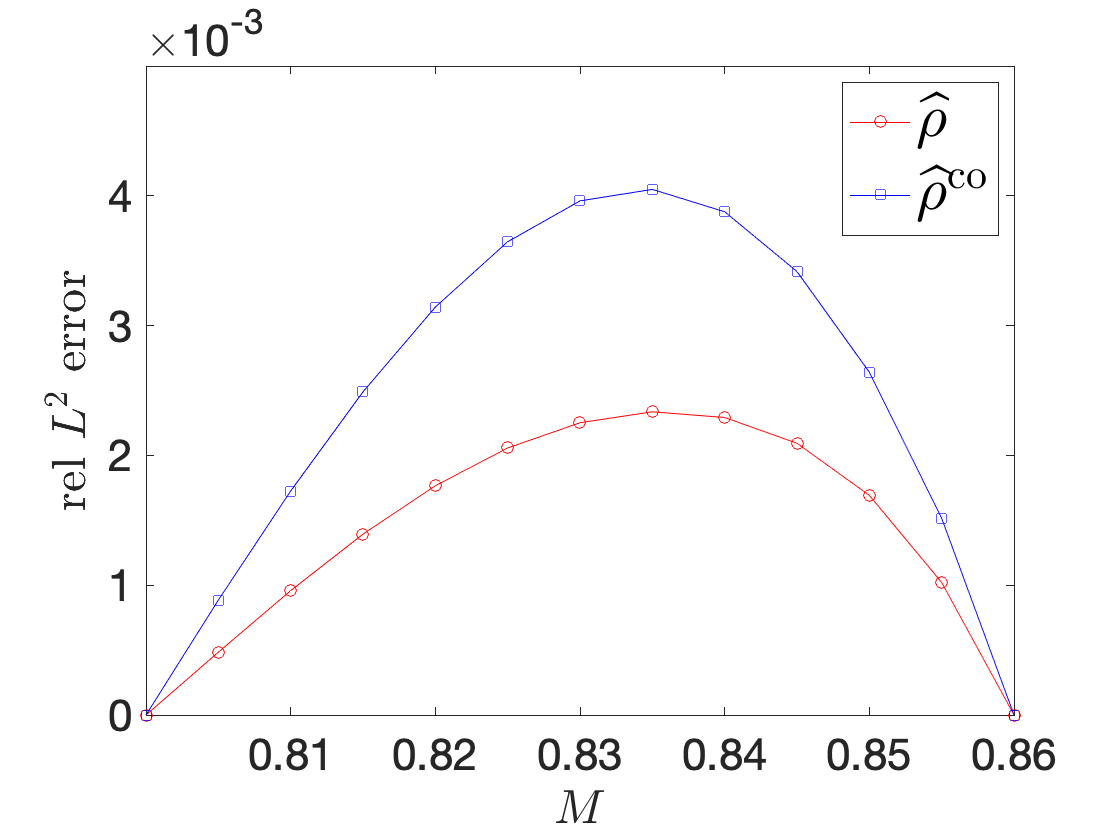}}

\subfloat[]{
\includegraphics[width=.45\textwidth]{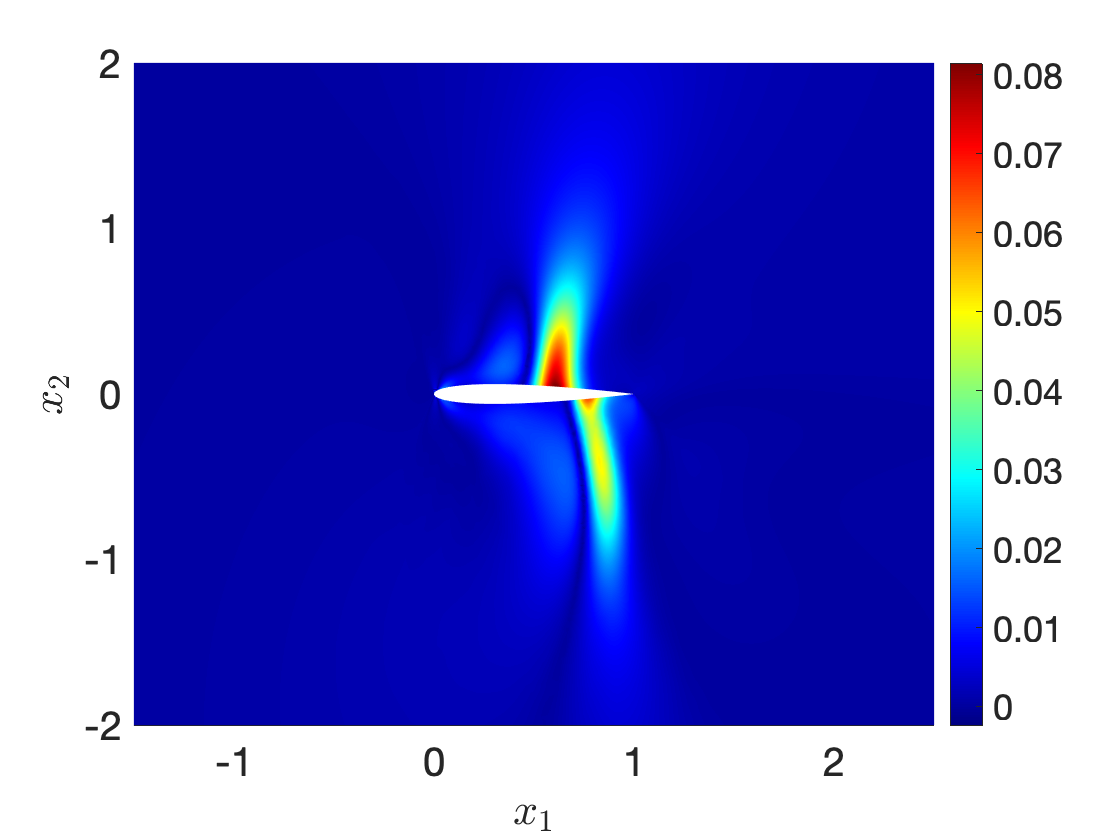}}
~~
\subfloat[]{
\includegraphics[width=.45\textwidth]{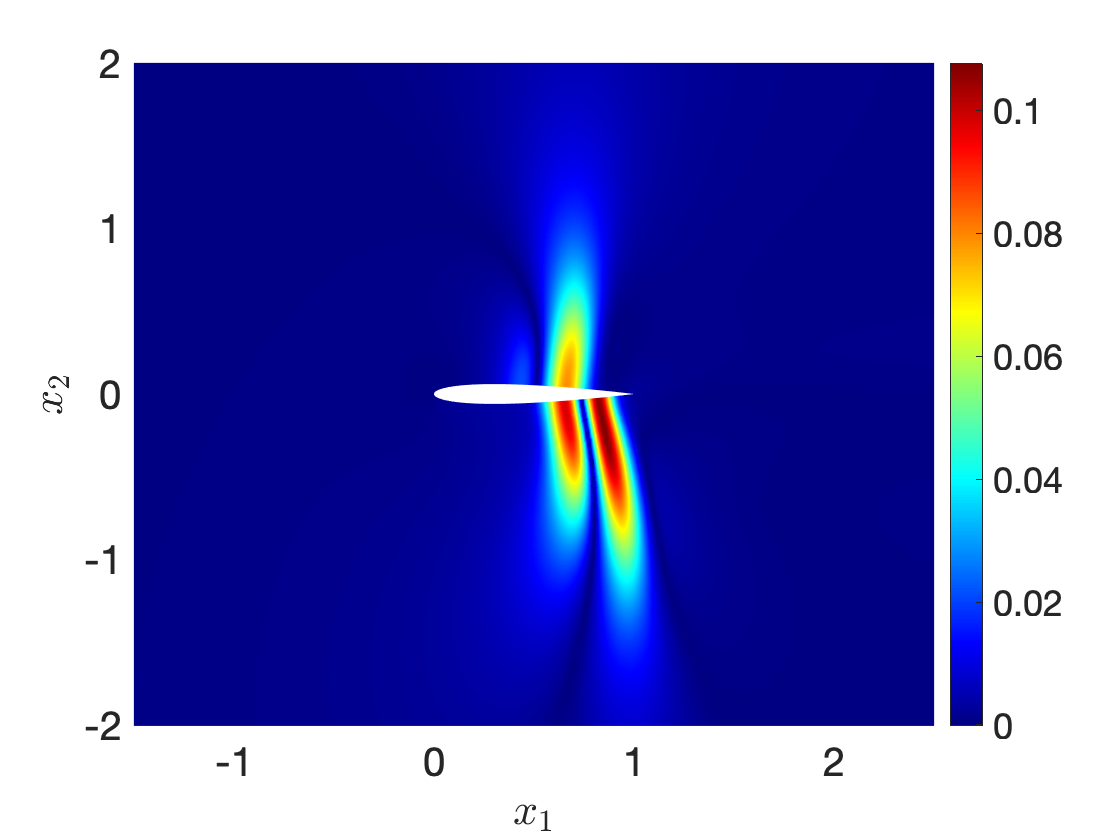}}
        
\caption{boundary-aware transportation of Gaussian models (angle of attack $1^o$). (a) behavior of the relative $L^2$ error for BA CDI.
(b)-(c) behavior  of the error fields
$|   \widehat{\rho}(s,x)  -  {\rho} (x; M) |$ and
$|   \widehat{\rho}^{\rm co}(s,x)  -  {\rho} (x; M)  |$
for $M=0.83$ and $s=0.5$. 
 }
\label{fig:airfoil_density_bnd_aware_2_angle1}
\end{figure}

\begin{figure}[H]
\centering
\subfloat[$x_2=0.1$]{
\includegraphics[width=.45\textwidth]{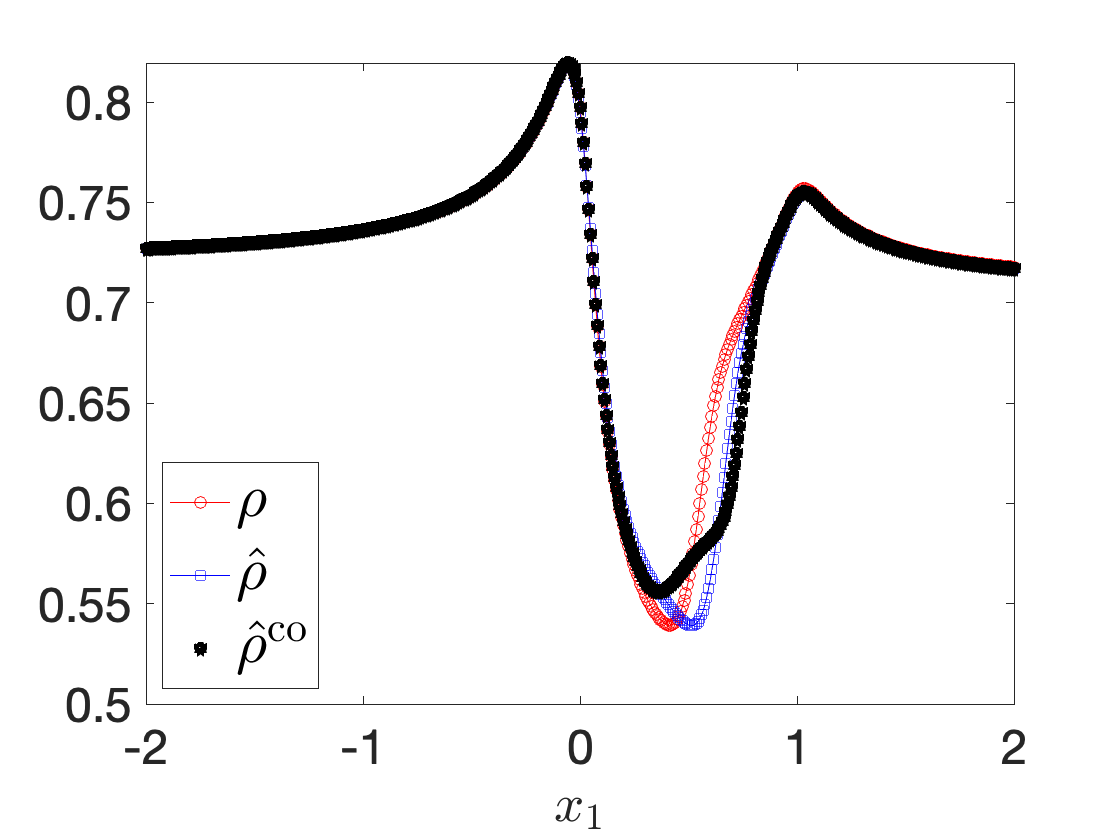}}
~~
\subfloat[$x_2= -0.1$]{
\includegraphics[width=.45\textwidth]{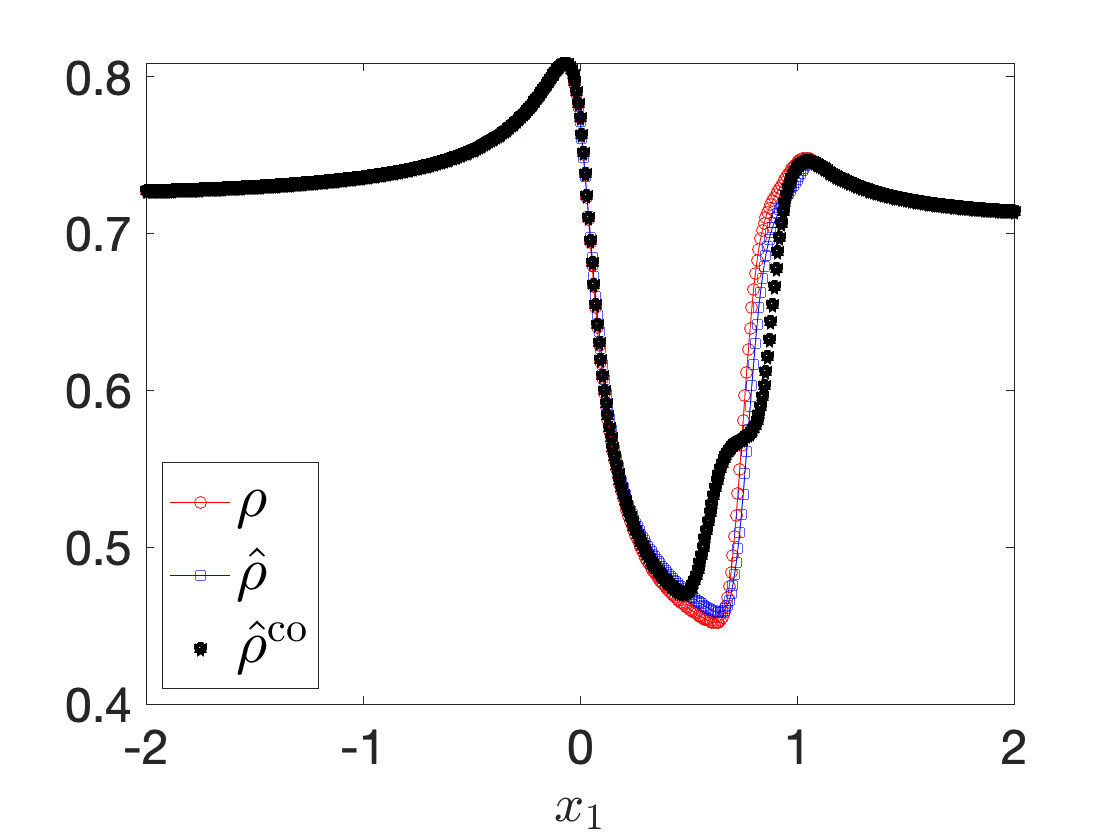}}
        
\caption{boundary-aware transportation of Gaussian models (angle of attack $1^o$). (a)-(b)  horizontal slices of  truth and predicted density profiles for $M=0.83$.
 }
\label{fig:airfoil_density_bnd_aware_3_angle1}
\end{figure}

Figures \ref{fig:airfoil_density_bnd_aware_angle1}, \ref{fig:airfoil_density_bnd_aware_2_angle1} and
\ref{fig:airfoil_density_bnd_aware_3_angle1} show performance of our nonlinear interpolation procedure.
Figure \ref{fig:airfoil_density_bnd_aware_angle1},
shows the BA  CDI  for two values of the parameter $s$: we observe that the interpolation procedure is able to generate physically-meaningful interpolations. Figure \ref{fig:airfoil_density_bnd_aware_2_angle1} compares the behavior of the relative $L^2$ error for BA CDI   and linear 
convex  interpolation: similarly,  Figure \ref{fig:airfoil_density_bnd_aware_3_angle1} shows horizontal slices of truth and predicted density profiles   for $M=0.83$ --- we set $s=1/2$ for both linear and nonlinear interpolation. Note that the shock on the upper part of the airfoil is not tracked as accurately as the lower shock by our nonlinear interpolation: this might be due to the inaccuracy of the Gaussian model and might also be due to the fact that the optimal value of $s$ is not necessarily a linear function of $\alpha = \frac{M-M_0}{M_1-M_0}$ (see discussion in 
section \ref{sec:motivating_examples} and results in 
Figure \ref{fig_sw4}).

\section{Conclusions}
\label{sec:conclusions}

We presented a general interpolation technique based on optimal transportation of Gaussian models for parametric advection-dominated problems. Application of optimal transportation to a Gaussian model of the  solution field, for which the transport map is known explicitly, simplifies the implementation of the method and allows to deal with fields that are neither scalar nor positive and that do not satisfy conservation of mass over the parameter domain.
We presented several examples to establish the connection between self-similarity and optimal transportation, which is at the foundation of the proposed technique. Furthermore, we presented several numerical investigations to illustrate the many features of the approach and assess strengths and weaknesses.

As  discussed in section \ref{sec:gauss_model_coherent} --- and shown numerically in the example of section  \ref{sec:transonic_flow} --- the choice of the Gaussian distribution model might not suffice to properly track relevant coherent structures of the flow; furthermore, the  approach is not suited to accurately represent the flow in the proximity of the boundaries. To address these issues, 
we proposed in section \ref{sec:extension} a more sophisticated interpolation procedure that combines Gaussian mixture models with a nonlinear registration  procedure.

The key elements  of our approach are (i) a scalar testing function $\mathcal{T}(\cdot; U): \mathbb{R}^n \to \mathbb{R}$, (ii) a (non-necessarily conforming) mapping technique for the construction of $T_g,W_g$, and (iii) a registration  (or mesh morphing) procedure to 
{project  the mappings $T_g,W_g$ onto a suitable subspace of bijective maps in $\Omega$}.
In this work, we proposed simple yet effective strategies based on (i)  physics-informed scalar testing functions, (ii) Gaussian models and optimal transportation maps between Gaussian distributions, and (if needed) (iii) the registration approach proposed in \cite{taddei2021registration}. In the future, we aim to design more accurate strategies for each of the three steps, and discuss the application  to  a broad class of problems in computational mechanics.

The aim of this work is to devise a nonlinear interpolation procedure  for continuum mechanics applications
that is simple to implement, interpretable,  and robust for small training sets.  
{In particular, we wish to apply our interpolation procedure to flow visualization problems.
Flow visualization methods such as particle image velocimetry
(PIV, \cite{adrian2011particle}) are of paramount importance in experimental fluid dynamics to investigate the flow behavior for relevant physical systems and ultimately inform the design and the assessment of engineering components.
These methods typically rely on piecewise-linear approximations of the form \eqref{eq:convex_interpolation} to estimate the flow field over a prescribed time interval:
if the acquisition frequency is small compared to the characteristic frequencies of the system, convex interpolations \eqref{eq:convex_interpolation} might be extremely inaccurate.  It is thus important to devise more advanced and physics-informed interpolations to achieve accurate predictions.
}

In the past decade, the spectacular successes of deep learning methods \cite{lecun2015deep} for data science applications have motivated the development of deep convolutational architectures in model reduction 
\cite{bhatnagar2019prediction,fresca2020deep,kim2022fast,lee2020model,mojgani2021low}:
in our experience, these approaches require large training sets and are difficult to interpret; furthermore, for small datasets, convergence to local minima might impact their robustness and generalization properties. In this regard, we observe that our approach might be interpreted as an Eulerian and non-intrusive counterpart of the 
registration-based approach in \cite{taddei2020registration,taddei2021registration}: we remark that the latter relies on the introduction of a template space and thus cannot deal with datasets of very modest size $n_{\rm train} = \mathcal{O}(1)$.
{In the future, we aim to devise strategies to optimally combine the many available linear and nonlinear reduction strategies for a wide range of offline computational budgets.
In this respect, 
similarly to \cite{bernard2018reduced}, 
we wish to apply the proposed technique  in the framework of  projection-based schemes, to augment the dataset of snapshots used to generate the reduced-order basis. 
 }

\section*{Acknowledgements}
The authors acknowledge the support by European Union’s Horizon 2020 research and innovation programme under the Marie Skłodowska-Curie Actions, grant agreement 872442 (ARIA).
Tommaso Taddei also acknowledges the support of IdEx Bordeaux (projet EMERGENCE 2019).

\appendix 

\section{Construction of the approximation map $\mathcal{N}$}
\label{sec:registration}

\begin{figure}[h!]
\centering
\includegraphics[width=0.5\textwidth]
 {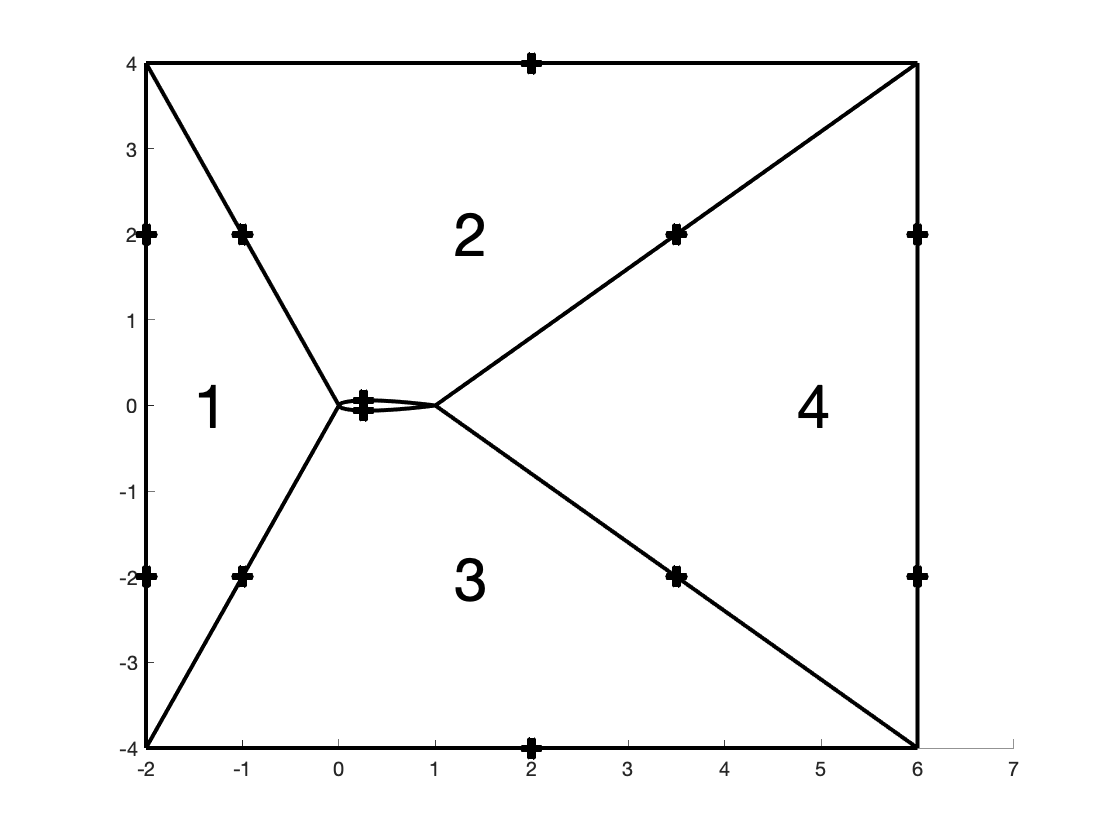}
  
 \caption{partition considered for the model problem  in section     \ref{sec:transonic_flow}.}
 \label{fig:vis_partitioned}
  \end{figure}

Let $\widehat{\Omega}=(0,1)^2$ be the  unit square,
let $\mathbb{P}_J$ be the space of one-dimensional  polynomials of degree lower or equal to $J$, and let 
 $\mathbb{Q}_J$ be the space of two-dimensional tensorized polynomials
  $$
  \mathbb{Q}_J ={\rm span}  \left \{ \varphi(x) = \ell_1(x_1) \ell_2(x_2) e_d :  \ell_1,\ell_2\in  \mathbb{P}_J, \;d \in \{ 1,2 \}
  \right \},
  $$ 
where  $e_1,e_2$ are the canonical basis of $\mathbb{R}^2$.
Given the domain $\Omega \subset \mathbb{R}^2$, we define the non-overlapping partition $\{ \Omega_q \}_{q=1}^{N_{\rm dd}}$  such that each element is isomorphic to the unit square; we denote by $\Psi_1,\ldots,\Psi_{N_{\rm dd}}$ Gordon-Hall maps from $\widehat{\Omega}$ to $\Omega_1,\ldots,\Omega_{N_{\rm dd}}$, respectively: we recall that Gordon-Hall maps are uniquely defined  based on the parameterizations of the partition interfaces.
To provide a concrete example, for the problem in section \ref{sec:transonic_flow} we consider the partition depicted in Figure  \ref{fig:vis_partitioned}.

Given the set of polynomials $\overrightarrow{{\varphi}} : = [{\varphi}_1,\ldots, {{\varphi}}_{N_{\rm dd}}] \in  \bigotimes_{q=1}^{N_{\rm dd}}   \mathbb{Q}_J$,  we define 
\begin{equation}
\label{eq:calN_gen}
\widetilde{\mathcal{N}} \left( x;  \overrightarrow{{\varphi}}  \right)  \, := \,
\sum_{q=1}^{N_{\rm dd}} \; 
\Psi_q \circ \Phi_q \circ \Psi_q^{-1}(x) \mathbbm{1}_{\Omega_q}(x) \;\; {\rm where} \;\;
\Phi_q(x) = x + \varphi_q(x) \; q=1,\ldots,N_{\rm dd}.
\end{equation}
It is possible to verify that the space 
\begin{equation}
\label{eq:calW_reg}
{\mathcal{W}}_0  = 
\left\{ 
\overrightarrow{{\varphi}} 
=
 [{\varphi}_1,\ldots, {{\varphi}}_{N_{\rm dd}}],
\in  \bigotimes_{q=1}^{N_{\rm dd}}   \mathbb{Q}_J
\, : \,
{\varphi}_q\cdot n|_{\partial \widehat{\Omega}} = 0,
\;
 \widetilde{\mathcal{N}} \left( \cdot ;  \overrightarrow{{\varphi}}  \right) \in C(\Omega)
\right\}
\end{equation}
is a linear space of size $M< 2 (J+1)^2 N_{\rm dd}$; we denote by 
$\left\{ \overrightarrow{{\varphi}}_m \right\}_{m=1}^M$ a basis of ${\mathcal{W}}_0$. Finally, we define the approximation map
${\mathcal{N}}:\Omega \times \mathbb{R}^M \to \mathbb{R}^2$ given by
\begin{equation}
\label{eq:calN_reg}
{\mathcal{N}} \left( x;  \mathbf{a} \right)  \, := \,
\widetilde{\mathcal{N}} 
\left( x;  \sum_{m=1}^M ( \mathbf{a} )_m  \overrightarrow{{\varphi}}_m  \right),
\quad
x \in \Omega, \quad
\mathbf{a} \in \mathbb{R}^M.
\end{equation}

\bibliographystyle{abbrv}   
\bibliography{all_refs}

\end{document}